\documentclass[11pt,a4paper]{article}
\usepackage{graphicx}

\usepackage{rotating}
\usepackage{subfigure}
\usepackage{caption} %must go after subfigure, wrapfigure, etc.
\usepackage{multirow}

\usepackage[T1]{fontenc}% pour les accents français
\usepackage{amsmath}
\usepackage{amssymb}
\usepackage{amsthm}
\usepackage{stmaryrd}
\usepackage{bbm}
\usepackage{delarray}
\usepackage{citesort}

\usepackage{textcomp}

\usepackage{fancyheadings}
\usepackage{a4}

{\theoremstyle{plain}
	\newtheorem{thm}{Theorem}[section]
	
	\newtheorem{prop}[thm]{Proposition}
	\newtheorem*{cor}{Corollary}
	\newtheorem*{con}{Conjecture}
}
{\theoremstyle{definition}
	\newtheorem{dfn}{Definition}
	\newtheorem*{ntn}{Notation}
}
{\theoremstyle{remark}
\newtheorem*{rem}{Remark}
}

\newcommand{\SP}[0]{\mbox{\,$\mathcal{S\!P}$}}
\newcommand{\UP}[0]{\mbox{\,$\mathcal{U\!P}$}}

\renewcommand{\O}{\emptyset}

\newcommand{\T}{\mathcal{I}}

\newcommand{\DefEmph}[0]{\sffamily\upshape \bf }

\newcommand{\mm}[0]{\emph{mutatis mutandis}}

\renewcommand{\O}{\emptyset}
\newcommand{\bi}[2]{{{#1}\choose{#2}}}
\newcommand{\dd}[1]{\frac{\delta}{\delta #1}}

\newcommand{\second}[0]{\(2^{\mbox{\tiny{nd}}}\) }

\newcommand{\nth}[1]{\({#1}^{\mbox{\tiny{th}}}\)}

\newcommand{\EE}[3]{E\bigg[\frac{#1}{#2}{#3}\bigg]}
\newcommand{\Ex}[4]{E_{#1}\bigg[\frac{#2}{#3}{#4}\bigg]}
\newcommand{\lr}[1]{\left(#1\right)}
\newcommand{\lrf}[2]{\left(\frac{#1}{#2}\right)}

\begin{document}

\title{		Families of $m$-convex polygons: $m = 2.$} 
\author{	W. R. G. James, I. Jensen and A. J. Guttmann,\\
			Department of Mathematics and Statistics,\\
			The University of Melbourne, Vic. 3010, Australia}
\date{\today}
\maketitle
\begin{abstract}
	Polygons are described as almost-convex if their perimeter differs
	from the perimeter of their minimum bounding rectangle by twice their
	`concavity index', $m$. Such polygons are called \emph{$m$-convex}
	polygons and are characterised by having up to $m$ indentations in
	the side. We use a `divide and conquer' approach, factorising 2-convex
	polygons by extending a line along the base of its indents. We then use
	the inclusion-exclusion principle, the Hadamard product and extensions to
	known methods to derive the generating functions for each case.
\end{abstract}

\section{Introduction} \label{s_intro}

This is the second in a series of papers that look at families of $m$-convex
self-avoiding polygons (SAPs). In the first \cite{wrgj} we outlined the 50-year
history of polygon enumeration on the square lattice before enumerating the
$m=1$ case. This began when Temperley \cite{Temperley:56} defined enumeration
problems involving self-avoiding walks (SAWs) that are necessarily closed,
forming SAPs. Exact results in this area have generally been difficult to
obtain, whereas asymptotic results have been more numerous. For example,
Hammersley \cite{Hammersley:62} showed that there exists a certain exponential
asymptotic growth in the number of SAWs, counted by their length, and SAPs,
counted by their perimeter, which is known to be the same for a given lattice.
Furthermore, for length $n$, it is believed that their asymptotic behaviour is
described by $\mu^n n^{\gamma-1}$, where $\mu$ is the growth constant and
$\gamma$ the critical exponent.

\begin{figure}
	\begin{center}
	\subfigure[A pyramid.]
				{\qquad\includegraphics{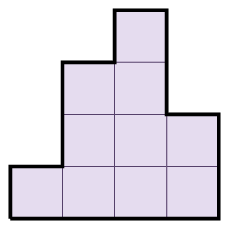}\qquad}
	\subfigure[A unimodal, or directed convex polygon.]
				{\qquad\includegraphics{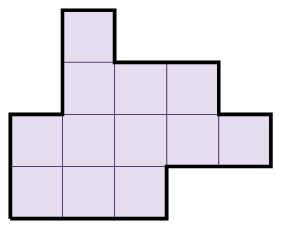}\qquad}
	\\
	\subfigure[A convex polygon.]
				{\qquad\includegraphics{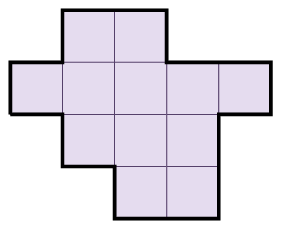}\qquad}
	\subfigure[A column-convex polygon.]
				{\qquad\includegraphics{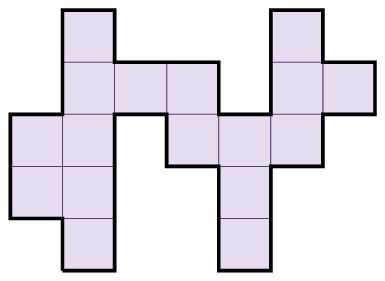}\qquad}
	\end{center}
	\caption{Some examples of convex polygons that are in some way convex.}
	\label{Fig:convex_examples}
\end{figure}

Exact results have so-far required the restriction of SAPs to subclasses that
are in some way {\em convex}. In two dimensions, convexity means that the
perimeter is equal in length to the minimum bounding rectangle (MBR).
Column-convexity means that any vertical cross-section may only intersect the
the polygon twice, such that all columns are connected. Examples of convex and
column-convex polygons can be seen in Figure~\ref{Fig:convex_examples}.
In 1997, Bousquet-M\'elou and Guttmann \cite{BMG:convex} (BMG) gave exact
results for convex SAPs of three dimensions and a method for their enumeration
in an arbitrary dimension \cite{BMG:convex}, which was based on the principle of
inclusion-exclusion. This became central to the methods used in the previous
paper \cite{wrgj}, and will be one of several methods used throughout this
paper.

\begin{figure}
	\begin{center}
	\subfigure[A 1-convex polygon.]
				{\qquad\includegraphics[scale=0.9]{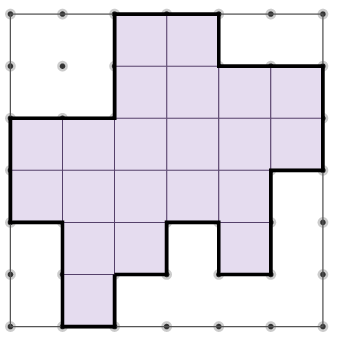}\qquad}
	\subfigure[A 2-convex polygon.]
				{\qquad\includegraphics[scale=0.9]{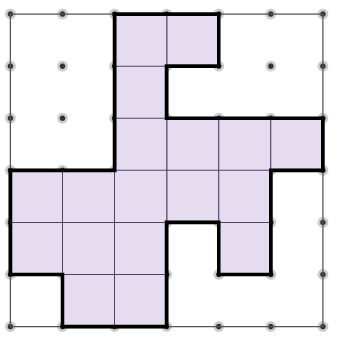}\qquad}
	\end{center}
	\caption{Some examples of almost-convex polygons, with the minimum bounding
	rectangle marked.}
	\label{Fig:almost_examples}
\end{figure}

Enting et al. \cite{Enting:92} described polygons as almost-convex if their
perimeter differs from the perimeter of their minimum bounding rectangle by
twice their `concavity index', $m$. Such polygons are called \emph{$m$-convex}
polygons and are characterised by having up to $m$ indentations in the side.
Examples of 1-convex and 2-convex polygons can be found in
Figure~\ref{Fig:almost_examples}.
Enting et al.  proceeded to derive the asymptotic behaviour of the
number of $m$-convex polygons according to their perimeter, $n$ for $m =
o(\sqrt{n})$. The results were confirmed for the case $m=0$ (i.e. convex
polygons) by the known perimeter generating function.  Soon after their paper was
submitted, Lin~\cite{Lin:1-convex} derived the exact generating function for
1-convex polygons, using a `divide and conquer' technique introduced to the
problem of convex animals (the interior of a convex SAP) by Klarner and Rivest
\cite{Klarner:74}. His result provided support for a conjecture in
\cite{Enting:92}, giving the next term in the asymptotic expansion.
The previous paper \cite{wrgj} then re-derived the generating functions for
1-convex polygons in an effort to generalise the methodology and extended the
results to osculating\footnote{Osculating polygons are those that may touch
themselves, but not cross.} and neighbour-avoiding\footnote{Neighbour-avoiding
polygons are those that may not occupy a neighbouring lattice vertex without
being connected by an arc.} polygons.

So why are we interested in almost-convex polygons? Column-convex polygons are
asymptotically more numerous than our almost-convex polygons. So if we are
trying to generalise to obtain SAPs, shouldn't that be a faster avenue of
investigation?
The answer to this question lies in the fact that we are interested in
understanding how the generating functions change as we gradually relax the
convexity restriction in all directions. For example, as deeper indentations are
allowed in the side of the polygons, we may start to see dense and complex
shapes that become fractal patterns in the limit as the hole depth grows
significantly fast with respect to the perimeter. Column-convex polygons, in
contrast, lose the convexity restriction in one direction altogether, and we
do not learn anything about the complex behaviour of those polygons that are the
least convex.
We would therefore like to predict how the form of the generating functions
changes as the concavity index grows. This will give us a means of understanding
what happens in the scaling limit as the concavity index grows in proportion to
the perimeter.

In this paper, we start be revisiting the factorisation used by Lin
\cite{Lin:1-convex} for almost-convex polygons. We then introduce the
methodology used in our results, followed by examples of its application in
re-deriving 1-convex generating functions. We then enumerate all the separate
building blocks required in the factorisation of 2-convex polygons, before
enumerating 2-staircase, 2-unimodal and 2-convex polygons using the `divide and
conquer' approach.\\

The Temperley method is central to the enumeration of partially convex polygons.
The so-called `functional-Temperley' method allowed Bousquet-Mélou
\cite{Bousquet:96} to enumerate classes of column-convex polygons.  It differs
from the Temperley method in that it can be used to concatenate several large
enumerable factors rather than individual columns, such that the building
blocks may be characteristically different.  A variation \cite{BoRe:02}
allowed for the enumeration of certain classes of animals, represented as heaps
of dimers. Rechnitzer \cite{Rechnitzer} identified these methods as equivalent,
the superiority of one over the other lying in its ease of use and
appropriateness to the recurrence relation underlying the problem.

The approach that has been used repeatedly by Lin
\cite{Lin:Chang88,Lin:90,Lin:91} in the construction of convex polygons is to
build them up vertically, block by block.  The functional-Temperley method can
therefore be used in this case.  One tool which can be used
to `join' polygonal blocks together is the Hadamard product (see
\cite{Rechnitzer}).  This is particularly useful when a few blocks need to be
joined in a non-recurring manner. It means that polygons with a side length
enumerated according to a certain parameter can be joined along those sides,
removing the arcs along the respective side perimeters and creating a larger
polygon. We will use the Hadamard product in calculating most of the generating
functions in this paper.

In enumerating 1-convex polygons, Lin ran a line along the bottom of the indent
of a 1-convex polygon, such that it is factored into three distinct parts.
Assuming that the indent occurs on the top, he divided all the possible polygons
into six cases, up to vertical symmetry.
In Section~\ref{s_bimodal} we extend these cases to apply to an indent of depth
$m$, which we call `bimodal' $m$-convex polygons, as the path has two modes in
the direction of the indent. We then calculate the generating function for the
cases $m=1,2$ and 3.

In Section~\ref{s_methodology}, we present the methodology required when
adopting a `divide and conquer' approach, such that we join together blocks
that are often enumerated using an inclusion-exclusion approach. We begin by
recalling the definition of the $E$ operator, which was central to the
inclusion-exclusion approach to enumerating 1-convex polygons. We then introduce
a notation that allows us to apply the inclusion-exclusion approach to distinct
blocks by making perimeter parameters independent in different blocks. In
Section~\ref{s_distinguished}, we show how we can sometimes avoid an
inclusion-exclusion argument by inserting indents at a distinguished vertex.
Next, Section~\ref{s_wrapping}, we then highlight how the block-by-block
consruction can be extended to `wrap' blocks that are attached to those
enumerated using the inclusion-exclusion approach. Then, in
Section~\ref{s_Hadamard}, we present the Hadamard product, which is required to
carry out the joining of such blocks. Finally, in Section~\ref{s_blocks}, we
enumerate the various staircase and unimodal blocks which we will need according
to their side perimeters.

Section~\ref{s_1-convex} serves as an example for the use of the above
techniques. We rederive the 1-convex polygon generating functions in a manner
that is far more concise than the pure inclusion-exclusion approach adopted in
\cite{wrgj}.  We then factorise 2-convex polygons by extending a line along the
base of its indents, and then derive the generating functions for 2-staircase,
2-unimodal and 2-convex polygons, in Sections~\ref{s_2-staircase},
\ref{s_2-unimodal} and \ref{s_2-convex} respectively.
We will see that as the examples get more complicated, requiring more subtle
inclusion-exclusion arguments, the combination of Lin's factorisation, the
Hadamard join and our general inclusion-exclusion arguments will allow for the
enumeration of any arbitrary family of SAPs that has a fixed factorisation. The
intermediary results, as well as much of the detail, are omitted for reasons of
conciseness and clarity. These may be found in \cite{thesis}, where the
presented results first appeared.

The result of this paper is that we now have all the pieces that we need to
enumerate $m$-convex polygons where each indentation has a pyramid shape, that
is, there are no turns or hooks. The next challenge, and the most interesting
outstanding problem, is to enumerate 3-convex polygons that have only one
indent, but an indent that forms a hook. It will be interesting to see what
happens to the denominator of the generating function in this case. This will
appear in the next paper in this series.
Going forward, it is not realistic to factorise almost-convex polygons
for high concavity indices as we have done here, as there will be a polynomial
growth in the number of cases to evaluate. It would be more sensible
to restrict the size of each indentation first, and then generalise these cases.
Eventually, defining operators that can add more and more complex indentations
in the side of convex polygons and looking at the effect of the asymptotic
growth of their number seems to be the path of least resistance in learning
more about how convex polygons become general SAPs.

\section{Bimodal polygons} \label{s_bimodal}

The generating function for pyramids, counting the width (by $x$) and the height
(by $y$), also given according to the base, $n$, of the polygons, was given as
\begin{equation} \label{Lin:P}
	P(x,y) = \sum_{n\geq1} P_n = \frac{xy(1-x)}{(1-x)^2-y},	
\end{equation}
where we write
\begin{equation} \label{Lin:P_n}
	P_n = E_{\{y\}}\bigg[ y^2 \left( \frac{x}{1-y} \right)^n \bigg]
\end{equation}
in terms of the operator $E$ given by Bouquet-Mélou and Guttmann
\cite{BMG:convex}.  The generating functions for unimodal and convex polygons
were respectively given as
\begin{equation} \label{Lin:H}
	H(x,y) = \sum_{n\geq1} H_n = \frac{xy}{\sqrt{\Delta}}
\end{equation}
and
\begin{equation} \label{Lin:C}
	C(x,y) = \sum_{n\geq1} C_n =
	\frac{xy A_c(x,y) }{\Delta^2} -\frac{4x^2y^2}{\Delta^{3/2}},
\end{equation}
where $A_c=1-3x-3y+3x^2+3y^2+5xy-x^3-y^3-xy(x+y+(x-y)^2)$ and $\Delta$ has the
usual definition, $1-2x-2y-2xy+x^2+y^2$.

We note that solutions that can be expressed purely in terms of rational
functions and the operator $E$ will be rational in terms of $u$ and $v$, where
\[
	x = u(1-v) \qquad \mbox{and} \qquad y = v(1-u).
\]
The generating functions of the latter two classes of polygon can therefore
be given according to their base as:
\begin{equation}
	H_n =  E_{\{y\}}\bigg[ \frac{y^2(1-y)(1-x-y)}{((1-y)^2-x)}
			\left( \frac{x}{1-y} \right)^n \bigg] - \frac{2xy^2u^n}{\Delta}
\end{equation}
and
\begin{equation}
	C_n = \frac{1}{\Delta^2}  E_{\{y\}}\bigg[ y^2(1-y)^2(1-x-y)^2((1+y)^2-x)^2
  \left( \frac{x}{1-y} \right)^n \bigg] - \frac{4x^2y^3u^n}{\SP\cdot\Delta^{3/2}},
\end{equation}
where $\SP=(1-x-y-\sqrt{\Delta})/2$ is the staircase polygon generating
function. We denote the generating function of
these functions of exact height $m$ to be \( X_{n,m} \), where $X$ is either
$P,H,$ or $C,$ while \( X^+_{n,m} = X_n - \sum_{i=1}^{m-1}  X_{n,m} \) is
that of polygons of \emph{at least} height $m$.\\

\begin{figure}
	\begin{center}
	\subfigure[Case 1.]{\includegraphics[scale=0.6]{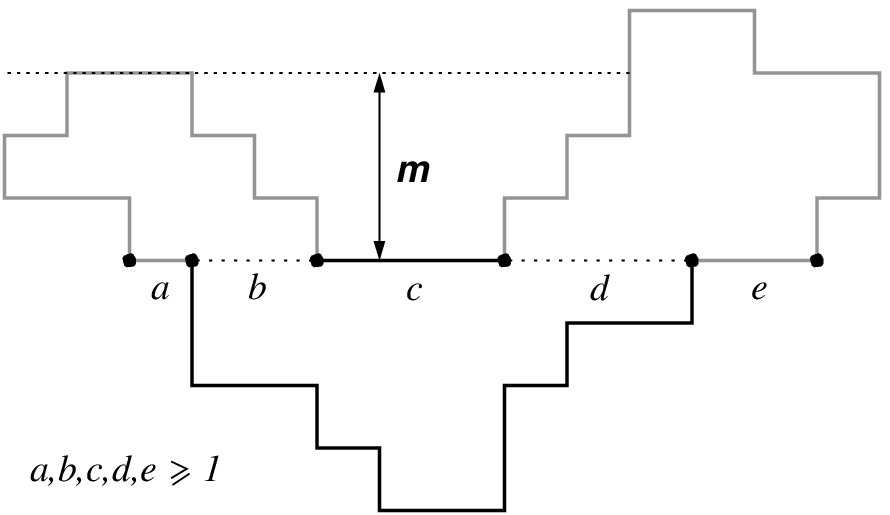}} \quad
	\subfigure[Case 2.]{\includegraphics[scale=0.6]{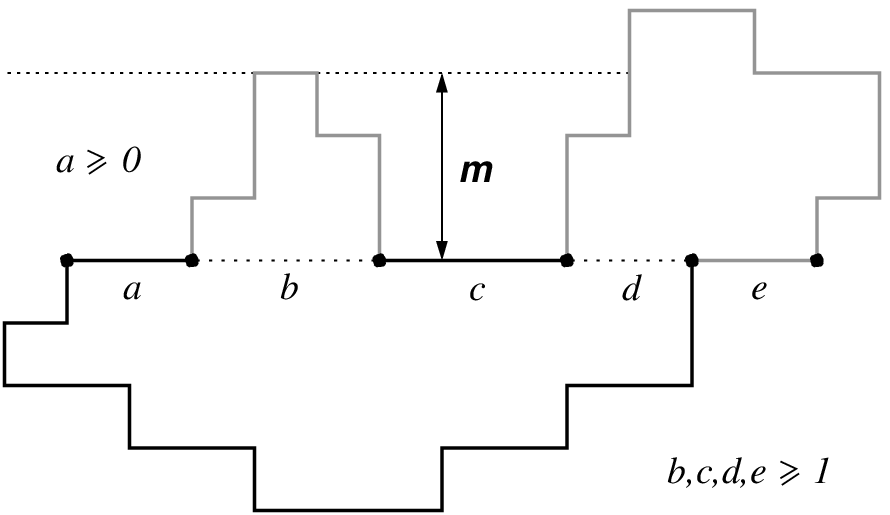}} \\
	\subfigure[Case 3.]{\includegraphics[scale=0.6]{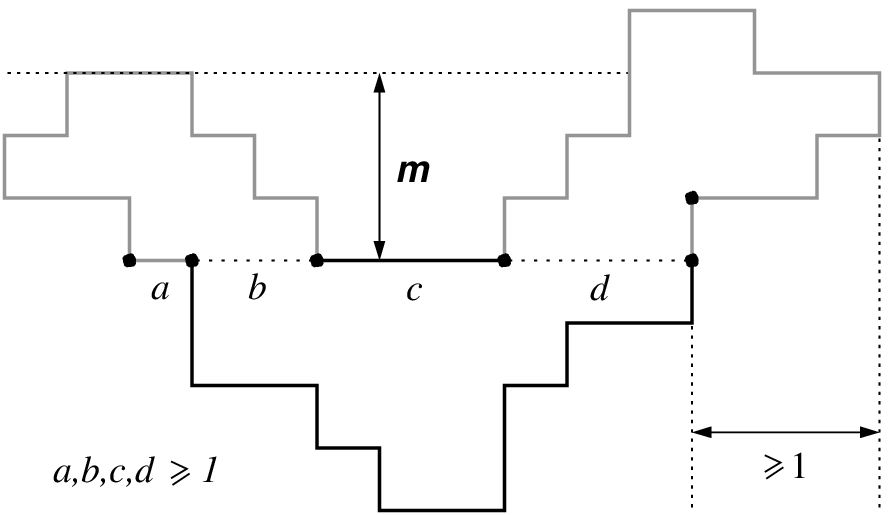}} \quad
	\subfigure[Case 4.]{\includegraphics[scale=0.6]{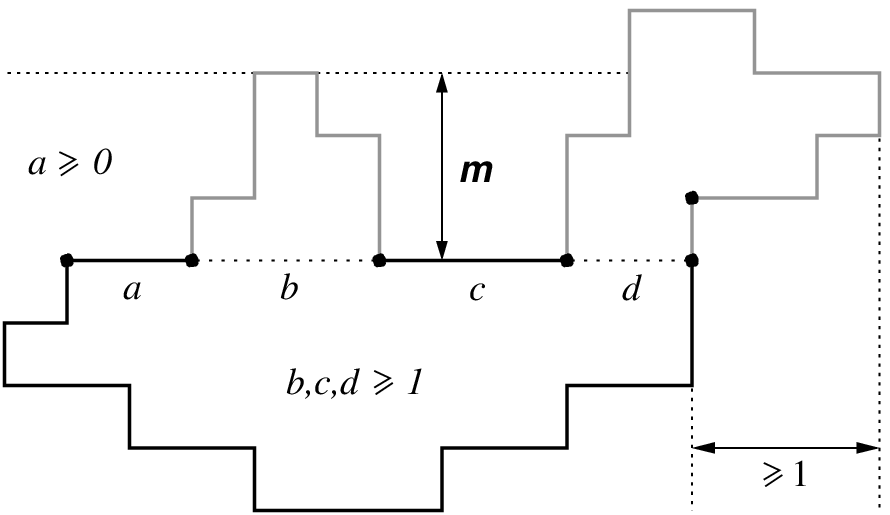}} \\
	\subfigure[Case 5.]{\includegraphics[scale=0.6]{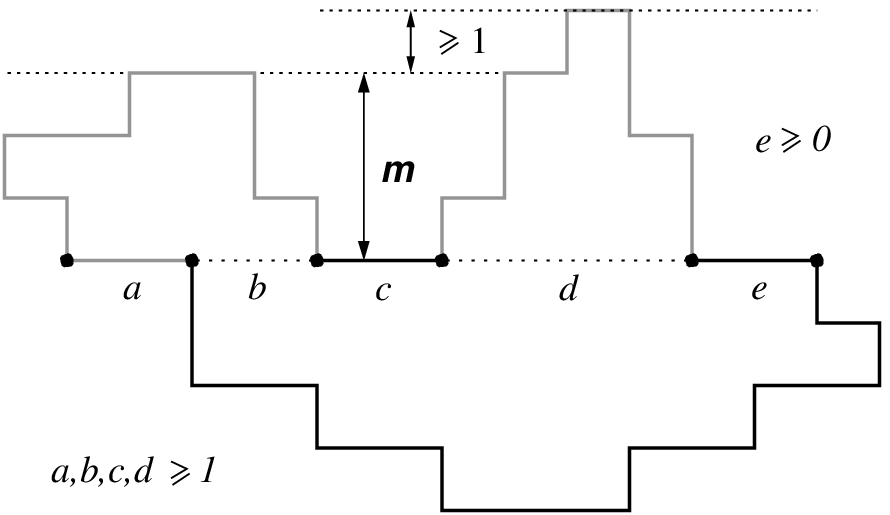}} \quad
	\subfigure[Case 6.]{\includegraphics[scale=0.6]{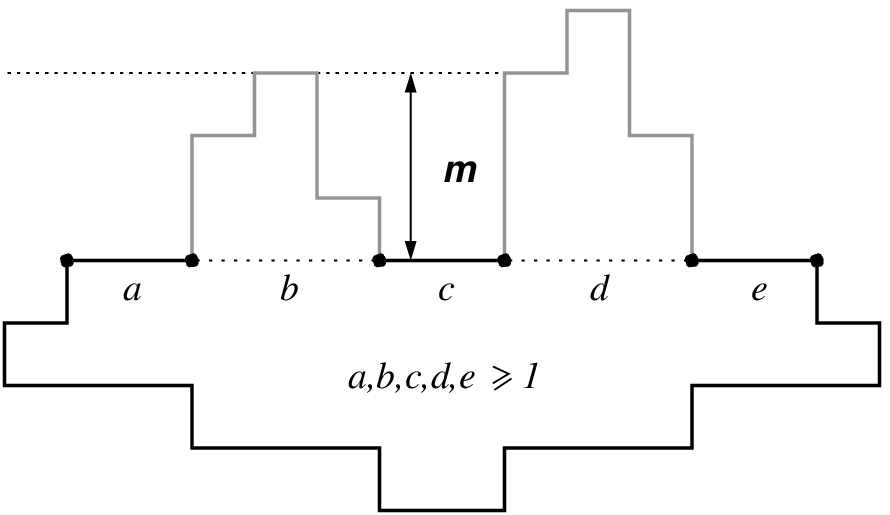}} \\
	\subfigure[Case 7.]{\includegraphics[scale=0.6]{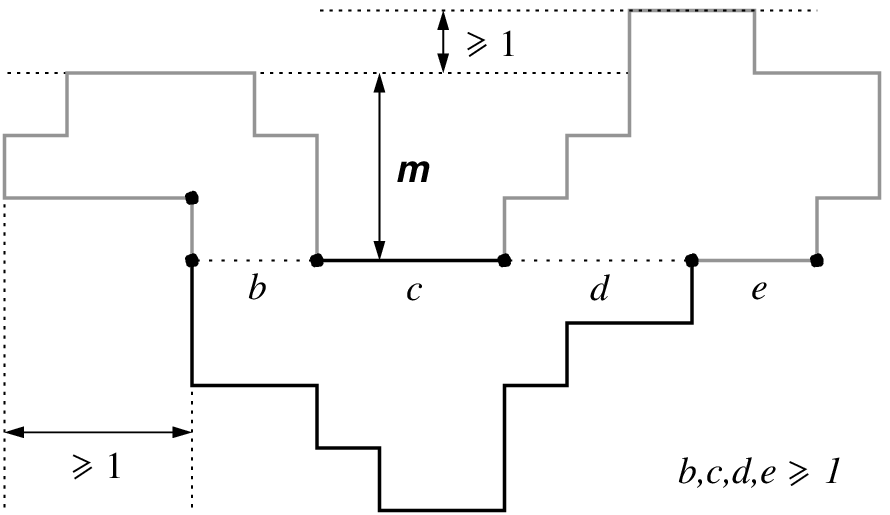}} \quad
	\subfigure[Case 8.]{\includegraphics[scale=0.6]{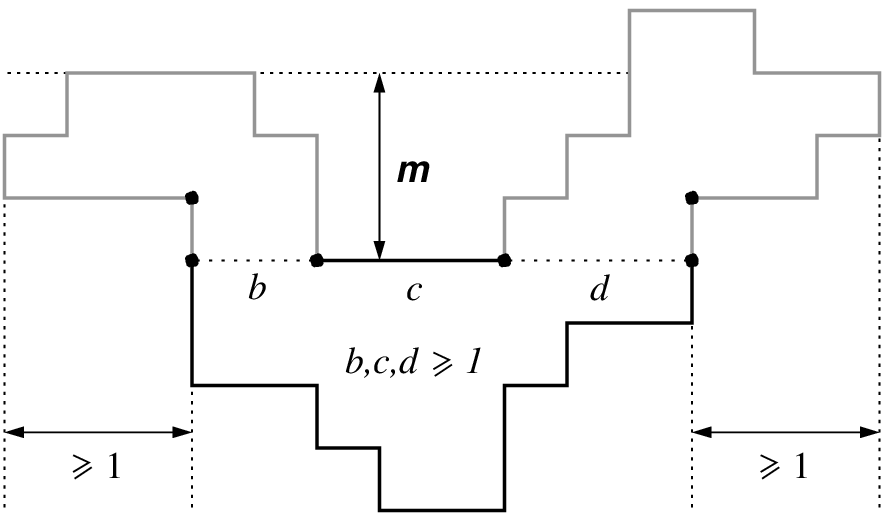}} \\
	\subfigure[Case 9.]{\includegraphics[scale=0.6]{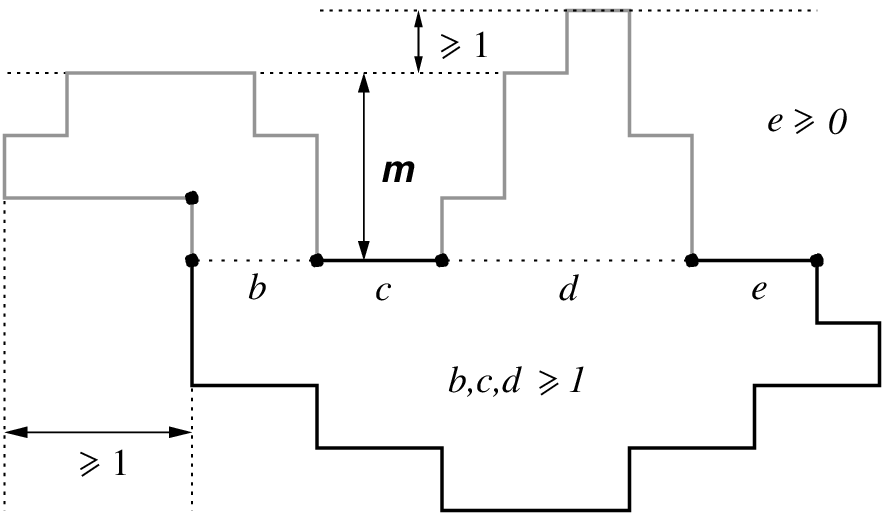}}
	\end{center}
	\caption{The cases used to construct bimodal $m$-convex polygons
			 following Lin's method.}
	\label{Fig:bimodal}
\end{figure}

These cases have a cardinality of two, due to the vertical symmetry.  Referring
to Figure~\ref{Fig:bimodal}, there are two symmetrical cases that are double
counted in cases 1 and 6, if the top factor is of height $m$. We also note the
symmetry of case 8. The generating functions of bimodal $m$-convex polygons,
which we denote as $G_i, i=1\ldots9$ for the nine classes given by Lin's
factorisation, are therefore given by
\begin{eqnarray*}
	G_1 &=& \sum_{a,b,c,d,e\geq1} x^{-b-d}
		H_{a+b,m} P_{b+c+d} ( 2 H^+_{d+e,m} - H_{d+e,m} ) \\
	G_2 &=& \sum_{a\geq0,b,c,d,e\geq1} 2 x^{-b-d}
		P_{b,m} H_{a+b+c+d} H^+_{d+e,m} \\
	G_3 &=& \sum_{a,b,c,d\geq1} 2 x^{-b-d}
		H_{a+b,m} P_{b+c+d} (H^+_{d,m} - P^+_{d,m}) \\
	G_4 &=& \sum_{a\geq0,b,c,d\geq1} 2 x^{-b-d}
		P_{b,m} H_{a+b+c+d} (H^+_{d,m} - P^+_{d,m}) \\
	G_5 &=& \sum_{a,b,c,d\geq1,e\geq0} 2 x^{-b-d}
		H_{a+b,m} H_{b+c+d+e} P^+_{d,m+1} \\
	G_6 &=& \sum_{a,e\geq0,b,c,d\geq1} x^{-b-d}
		P_{b,m} R_{a+b+c+d+e} ( 2 P^+_{d,m} - P_{d,m} ) \\
	G_7 &=& \sum_{b,c,d,e\geq1} 2 x^{-b-d}
		(H_{b,m} - P_{b,m}) P_{b+c+d} H^+_{d+e,m+1} \\
	G_8 &=& \sum_{b,c,d\geq1} x^{-b-d} (H_{b,m} - P_{b,m}) P_{b+c+d}
		(2(H^+_{d,m}- P^+_{d,m}) - (H_{d,m} - P_{d,m})) \\
	G_9 &=& \sum_{b,c,d\geq1,e\geq0} 2 x^{-b-d}
		(H_{b,m} - P_{b,m}) P_{b+c+d+e} P^+_{d,m+1}.
\end{eqnarray*}

Implementing the sum using any of the standard computer algebra programs then
verifies the result for 1-convex polygons.
The generating function for bimodal 2-convex polygons is
\begin{equation}\label{Eq:2-bimodal}
	\frac{2x^2 A^{(2)}_c}{(1-x)^3((1-x)^2-y)\Delta^{5/2}}
+	\frac{x^2 B^{(2)}_c}{(1-x)^5((1-x)^2-y)^3\Delta^3}
\end{equation}
where {\small
\begin{align*}
A^{(2)}_c =& -4(1-x)^{11} + 2 (1-x)^9 (10+9x) y - 2 (1-x)^7 (20+21+15x^2) y^2 \\&
			+ 2 (1-x)^5 (40+17x+18x^2+21x^3) y^3 - (1-x)^3 (20-26x-31x^2 \\&
			-14x^3+3x^4) y^4
			+ (1-x) (2+x) (2-13x+5x^2-7x^3-3x^4) y^5 \\&
			+ x (4-11x-3x^2-5x^3-x^4) y^6 + x(1+x)^2 y^7
\end{align*}
\begin{align*}
B^{(2)}_c =& 8(1-x)^{18} - 4 (1-x)^{16}(16+11x) y + 8(1-x)^{14}(28+29x+12x^2) y^2 \\&
	- 2 (1-x)^{12} (224+233 x+158 x^2+51 x^3) y^3 + 2 (1-x)^{10} (280+189 x \\&
	+144 x^2+95 x^3+24 x^4) y^4 - (1-x)^8 (448-57 x-56 x^2+36 x^4+x^5) y^5 \\&
	- (1-x)^6 (-224+345 x+182 x^2+336 x^3+292 x^4+75 x^5+6 x^6) y^6 \\&
	+ (1-x)^4 (-64+265 x+30 x^2+128 x^3-496 x^4-267 x^5+66 x^6+2 x^7) y^7 \\&
	- (1-x)^2 (-8+89 x-132 x^2-634 x^3-558 x^4+593 x^5-248 x^6-64 x^7 \\& +2 x^8) y^8
	+  x (7-236 x-476 x^2+962 x^3-388 x^4-296 x^5+248 x^6-78 x^7 \\& +x^8) y^9
	+ x (13+188 x+99 x^2-288 x^3+247 x^4-84 x^5+17 x^6) y^{10} \\&
	- x (11+58 x-28 x^2+6 x^3+x^4) y^{11}-(-3+x) x (1+x) y^{12}
\end{align*}
}The generating function for bimodal 3-convex polygons was obtained and its form is
\[
	\frac{x A}{(1-x)^5((1-x)^2-y)^5y^4\Delta^{5/2}}
+	\frac{x B}{(1-x)^9((1-x)^2-y)^{11}y^4\Delta^3}.
\]
Generating functions for higher concavity indices can be calculated, but we
shall see later that there is a much more efficient way of
calculating the generating functions, whilst still adopting the general
methodology.

For reference, we point out that the form of the 1-convex generating function
was \[
	\frac{2x^2y A}{(1-x) \Delta^{5/2}}
+	\frac{x^2y B}{(1-x)((1-x)^2-y) \Delta^3}.
\]
We note that the dominant factor asymptotically is that of $\Delta$ in the
denominator, whose index is not changing. This implies that this class of
polygon, as a proportion of the total number of $m$-convex polygons, is
asymptotically small. For example, if there are $m$ separate height one
indents, for large $n$ there is approximately $\bi nm$ ways of choosing the
location of these indents, which implies an \nth m derivative of the convex
polygon generating function. This will increase the index of $\Delta$ by $m$ and
there will be $\sim n^m$ times as many as the bimodal polygons in the limit.

\section{Methodology} \label{s_methodology}

Let us recall the $E$ operator that was originally defined in \cite{BMG:convex}
and was central to the inclusion-exclusion approach to enumerating 1-convex
polygons in the previous paper.

\begin{dfn}[$E_I$]
	For $I \subset \{1, \ldots, d\}$, $E_I$ acts on a Laurent series in the
	variables $x_1,\ldots,x_d$ with real coefficients. Let $f(x_1,\ldots,x_d) =
	\sum a_{n_1,\ldots,n_d} x_1^{n_1} \cdots x_d^{n_d}$. Then
	\begin{equation*}
		E_I(f(x_1,\ldots,x_d)) = \begin{cases}
			f(x_1,\ldots,x_d),
			&	\mbox{if } I = \O, \\
			\sum a_{n_1,\ldots,n_{i-1},2n_i,n_{i+1},\ldots,n_d}
				 x_1^{n_1} \cdots x_d^{n_d},
			&	\mbox{if } I = \{i\}, \\
			E_{I \backslash \{i\}} \left[ E_{\{i\}} (f(x_1,\ldots,x_d)) \right],
			&	\mbox{if } |I| > 1 \mbox{ and } i \in I.
		\end{cases}
	\end{equation*}
	Equivalently,
	\begin{multline*}
		E_{\{i\}} [f(x_1,\ldots,x_d)]
		= \frac{1}{2} \big[
			f(x_1,\ldots,x_{i-1}, \sqrt{x_i},x_{i+1},\ldots,x_d) \\
			+ f(x_1,\ldots,x_{i-1},-\sqrt{x_i},x_{i+1},\ldots,x_d)
		\big].
	\end{multline*}
\end{dfn}

\begin{ntn}[E]
	For our 2-dimensional case, $E_{[1,2]}\left[ f(x,y) \right]$ is denoted
	$E\left[ f(x,y) \right]$.
\end{ntn}

We now introduce a new notation, denoted ``star'' (*) notation, which designates
when we are holding enumeration variables constant while using the $E$ operator.
This allows us to apply the inclusion-exclusion approach to distinct blocks by
making perimeter parameters independent in different blocks. As summation is
commutative with $E$, this allows us to evaluate series expressions more simply.

\begin{ntn}[$*$]
	Any variable that is acted on by $E_I$ which is asterisked is assumed to be
	\emph{independent} of the variables in $I$. These variables are replaced by
	their original value outside the operator. For example,
	\[
		\EE 1{(1-x-y^*)(1-x-y)}{} = \frac 1{2(1-x-y)} \left(
			\frac 1{(1-x-\sqrt y)} + \frac 1{(1-x+\sqrt y)}
		\right).
	\]
\end{ntn}

For notational convenience, we denote the staircase and unimodal polygon generating
functions to be respectively
\begin{equation} \label{Eq:SP-UP}
	\SP = (1-x-y-\sqrt{\Delta})/2 (= u - x = v - y) \qquad \mbox{and} \qquad
	\UP = x y Z,
\end{equation}
where $Z$ is the generating function for pairs of directed walks, such that
\[
	Z = 1/\sqrt{\Delta} = 1/(1-u-v),
\]
and $\Delta$ has the usual definition, $1-2x-2y-2xy+x^2+y^2$. We define the
1-deep indent generating function to be
\begin{equation} \label{Eq:I}
	\T = \frac{u^2}{(1-u)^2} = \frac{\SP^2}{y^2}, \qquad \mbox{with} \qquad
	\bar{\T} = \T(y,x) = \frac{v^2}{(1-v)^2},
\end{equation}
and, for $m>1$, the $m$-deep indent generating function is
\begin{equation} \label{Eq:Im}
\T_m \equiv \T_m(x,y) = \frac{u^2}{(1-u)^{2m}} = \frac{v^2 \SP^{2m-2}}{y^{2m}},
	\qquad \mbox{with} \qquad \bar{\T}_m = \T_m(y,x).
\end{equation}

\begin{ntn}[Edge, side of a polygon]
	In \cite{BMG:convex}, staircase polygons were defined as those that had two
	roots, one in each corner of the MBR, and factored as two directed walks in
	opposite directions. Unimodal polygons were defined similarly,
	requiring that they be rooted in one corner and that there be a
	factorisation for each direction such that there are either only positive
	or only negative steps in that direction for each factor. Finally, convex
	polygons were defined as those for whom there is a unimodal factorisation
	for a given direction, taking the root to be different vertex for each
	direction. Each of these factors, plus any steps along the edge of the MBR)
	we refer to as an {\bf edge}. For example, the walk running from
	the bottom to the top of the MBR on the left (resp. right) side of a convex
	polygon (including the bottom and top perimeters) is called the left (resp.
	right) edge. We note that the top and bottom edges are defined similarly.

	Now, each edge can be factored into two (possibly empty) directed walks. We
	refer to these factors, plus any steps running along the MBR, as a
	{\bf side} of the convex polygon. That is, a convex polygon has four
	directed walks running between the sides of the MBR, and we call them sides
	of the polygon. We note that edges are equivalent to sides for staircase
	polygons. Also, for unimodal polygons rooted in the bottom left corner, the
	bottom edge is the bottom-right side, and the left edge is the top-left
	side.
\end{ntn}

Without loss of generality, we will usually assume that there is a vertical indent
on the top-left side of the polygon, and that staircase and unimodal polygons
are rooted to the bottom left. We can therefore classify 2-convex polygons
according to the relative direction and position of their indents. They are
either bimodal, or belong to one of four distinct cases:
\begin{enumerate}
\item	in the same direction, on the same edge;
\item	in the same direction, on opposite edges;
\item	in different directions, on the same side; and
\item	in different directions, on opposite edges.
\end{enumerate}
Moreover, we note that when the indents are on edges of unimodal or convex
polygons that have two distinct sides, then they can be broken-down further into
either adjacent or non-adjacent sides.

\subsection{Distinguishing vertices to insert indents} \label{s_distinguished}

The first application of Lin's factorisation enables us to insert indents by
first distinguishing a height at which the indent is to be inserted. We are able
to do this when the indents are joined to a staircase factor. We shall see that
this is possible because the factor of the generating function that corresponds
to the indent is independent of its position.

\subsubsection{Distinguishing a vertical step to insert the indent}
Let us consider bimodal $m$-staircase polygons. So as to have a
visual guide, we refer to Figure~\ref{Fig:bimodal}(b) and (d), which show their
form, except for their being rooted in the bottom-left corner. And so, measuring
the width of the bottom and indent factors from the left edge of the top factor,
their generating function is
\(
	\T_m (u/x)^d v,
\)
and that of the top factor is $u^d v$. Combining these sums gives
\begin{equation}
	\T_m v^2 \sum_{d\geq 1} \left(\frac{u}{1-v}\right)^d = \T_m u v^2 Z.
\end{equation}
We see that the factor of $\T_m$ is
independent of the height at which the indentation is found. 

\begin{prop} \label{Prop:s.bimodal}
	The generating function for bimodal $m$-staircase polygons whose single
	$m$-deep indent is vertical is
	\[
		\T_m y^2 \dd y \frac{\SP}y.
	\]
\end{prop}

\begin{proof}
	We saw above that the generating function that enumerates the indent is
	independent of the height of the indent. This is due to the simple form of
	the staircase polygon generating function with $n$ fixed horizontal steps,
	which is $u^nv$. The generating function can therefore be factored as
	the indent generating function, $\T_m$, and the generating function for
	staircase polygons with a distinguished left-most vertical step (other than
	the bottom one), which is given by $y^2 \dd y \SP/y$.
\end{proof}

\subsubsection{2-unimodal polygons with both indents on the left edge}

\begin{figure}
	\begin{center}
		\includegraphics[scale=0.5]{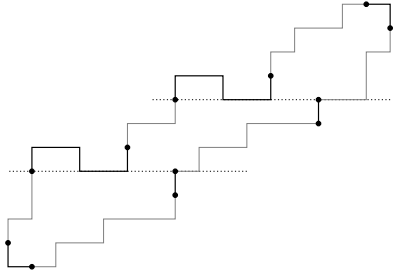}
	\end{center}
	\caption{The form of 2-staircase polygons with their indent on the same side
		and in the same direction.}
	\label{Fig:2.stair.same}
\end{figure}

This argument can be extended to cases where there are multiple indents on the
one side, giving staircase factors. (See Figure~\ref{Fig:2.stair.same}, for
example.) As each indent is attached to a staircase factor and is independent of
the height or length of each join, we can enumerate them by distinguishing a
vertical step and multiplying by $\T$. We note that this is true, no matter the
shape of the top factor. Therefore, this approach works for both 2-staircase and
2-unimodal polygons with their indents on the left.

\begin{figure}
	\begin{center}
	\subfigure[The two possible configurations, depending on whether the first
	height selected is above or below the second.]{
		\includegraphics[scale=0.5]{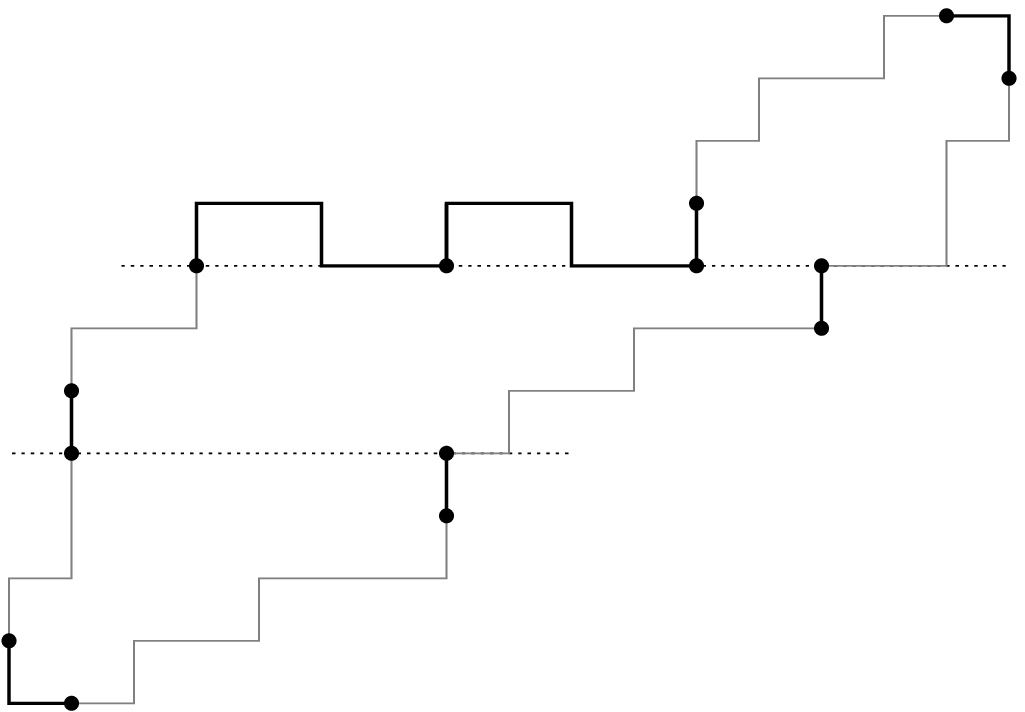}\qquad
		\includegraphics[scale=0.5]{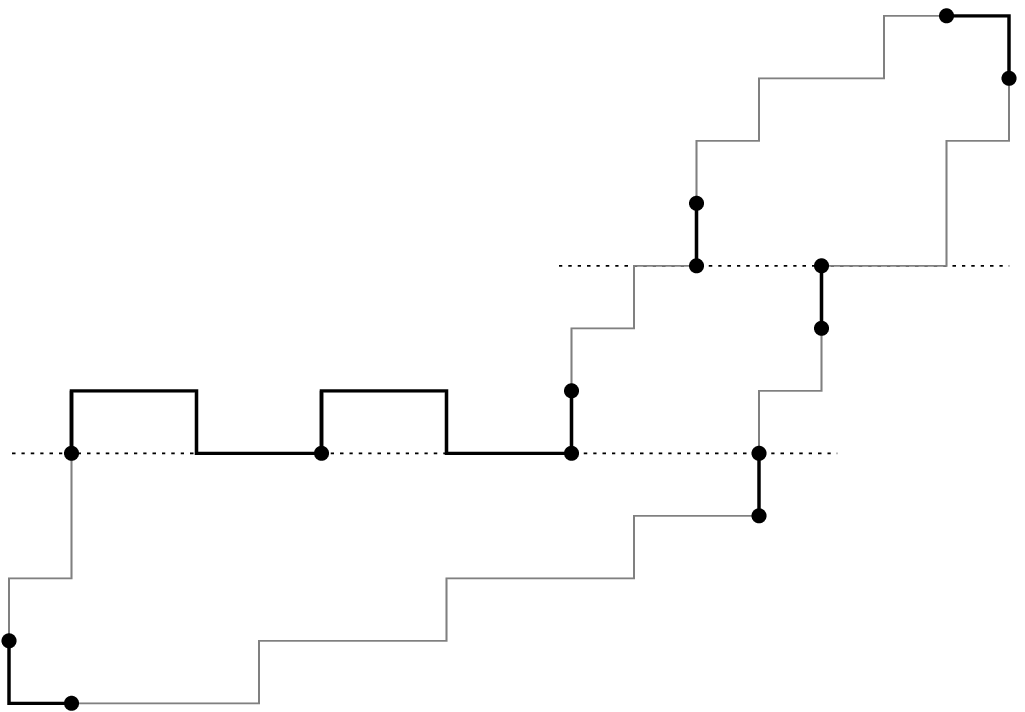}}
	\subfigure[The combined form. In both instances, the bottom walk can go
	below the two boxes formed. Only in one can it go inside.]{
\hspace{24pt} \includegraphics[scale=0.7]{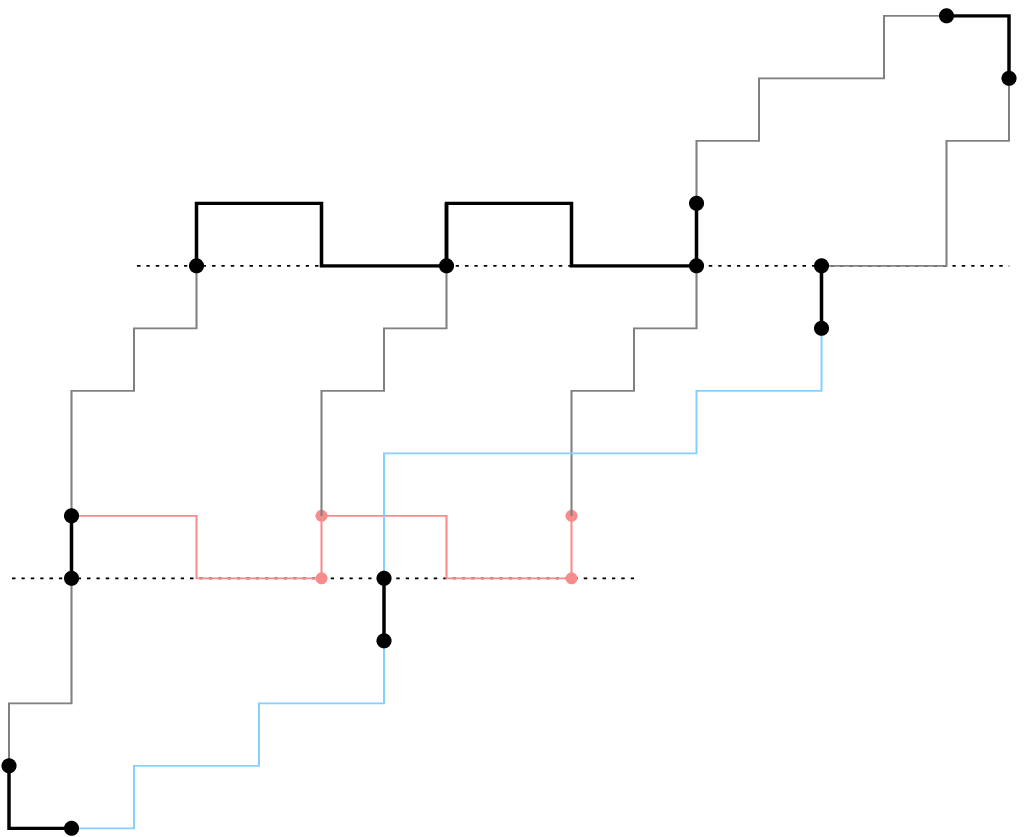} \hspace{24pt}}
%\hspace{48pt} \includegraphics[scale=0.5]{2-stair.total.eps} \hspace{48pt}}
	\end{center}
	\caption{The construction of 2-staircase polygons. Two vertical steps are
	selected as the heights of the indents. Both indents are placed at the first
	choice.}
	\label{Fig:2-stair}
\end{figure}

In the 2-staircase example, we select a vertical step on the
left side that is not the bottom-most to designate the height at which we will
place the indents. We then multiply by $\T^2$ to place the two indents.
Therefore, the case where the indents are next to each other is enumerated by
\(
	\T^2 y^2 \dd y \SP/y.
\)
We then take the derivative again to designate the location of the second
indent. This gives another factorisation line, either above or below the
indents, as shown in Figure~\ref{Fig:2-stair}(a). For a given pair of heights,
there are therefore two ways that is occurs: with the indents either at the
upper or lower height.
This creates the combined form shown in Figure~\ref{Fig:2-stair}(b). If we take
the walk between the two distinguished vertical steps that indicate the
designated heights and translate it horizontally to the edges of the indents,
we form two boxes. When the indents are at the higher of the two designated
heights, the bottom walk can go inside these boxes. Otherwise, it cannot go into
either of them.

We want to enumerate the case where the left occurrence of the two indents is at
the lower of the two heights and the right one is higher. This implies that the
lower walk must be able to go into the right box, no matter the order in which
we select the two heights, but never into the left box. We will now show that
the generating function for this is equal to the case where one can go into both
boxes if the upper height is selected first, and one cannot go into either box
otherwise.

\begin{prop}
	Let us consider a unimodal polygon with two staircase regions with identical
	sides, but possibly differing widths, identified as in
	Figure~\ref{Fig:2-stair}(b).  The generating function of those polygons
	whose bottom perimeter intersects the right box, but not the left, is equal
	to that of the polygons whose bottom perimeter intersects both boxes.
\end{prop}

\begin{figure}
	\begin{center}
	\subfigure[The bottom perimeter intersects the right box.]{
		\includegraphics[scale=0.5]{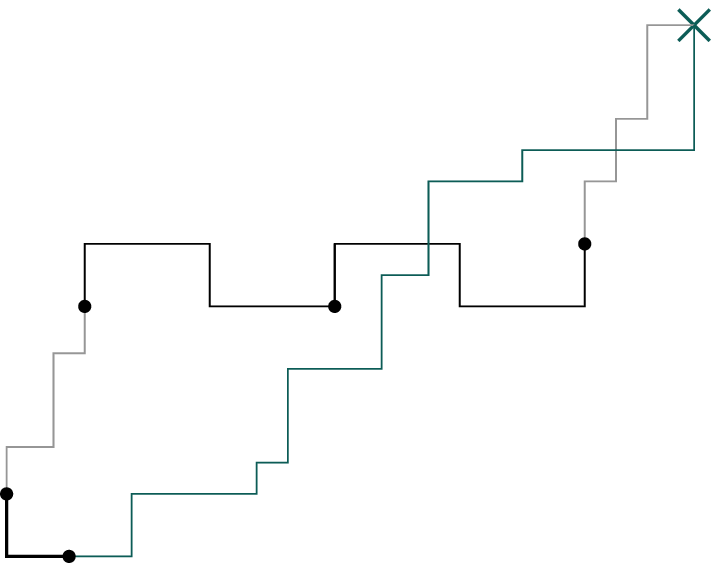}
		\includegraphics[scale=0.5]{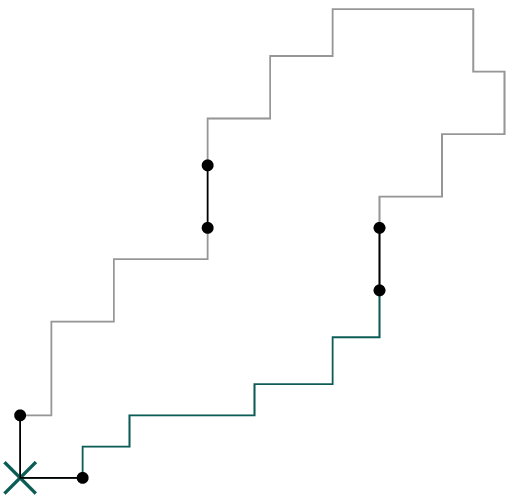}} \qquad
	\subfigure[The bottom perimeter intersects the left box.]{
		\includegraphics[scale=0.5]{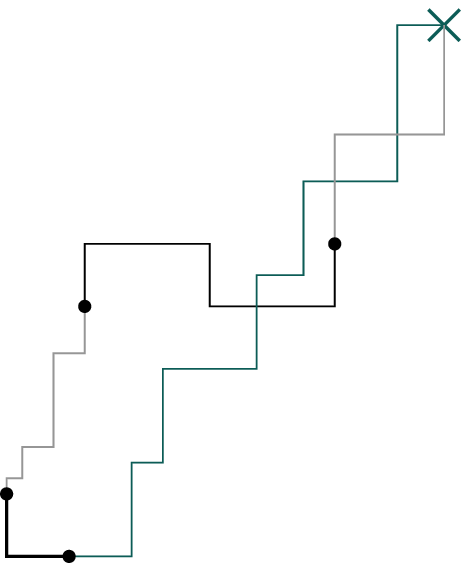}
		\includegraphics[scale=0.5]{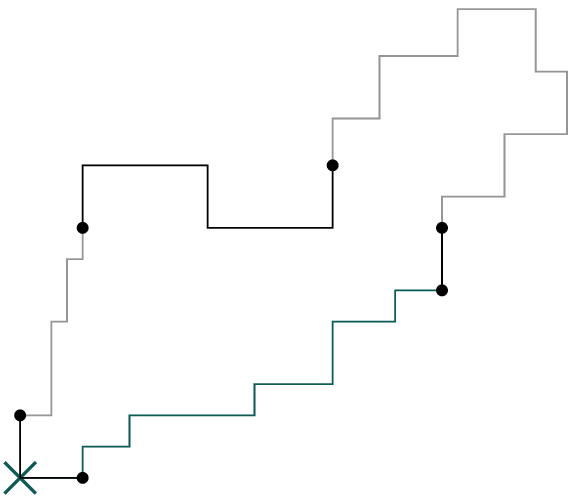}}
	\end{center}
	\caption{The 2-unimodal polygon can be factored where it intersects the
	boxes that delimit the staircase region between the possible indent
	positions.}
	\label{Fig:2-stair.factor}
\end{figure}

\begin{proof}
	Let us first consider the case where the bottom perimeter intersects the
	right box. This is illustrated in Figure~\ref{Fig:2-stair.factor}(a). If we
	factorise there, the bottom and top factors are respectively enumerated by
	\[
		\dd y \T^2 = 2\T \dd y \T \qquad \mbox{and} \qquad
		y^2 \dd y \frac{\UP} y.
	\]
	We then consider the case where the bottom perimeter intersects the
	left box. This is illustrated in Figure~\ref{Fig:2-stair.factor}(b).
	Following Proposition~\ref{Prop:s.bimodal} \mm, the bottom and top factors
	are respectively enumerated by
	\[	
		\dd y \T \qquad \mbox{and} \qquad \T \cdot y^2 \dd y \frac{\UP}{y}.
	\]
	The first generating function is twice the second, hence the difference is
	equal to that of the latter case, where the bottom perimeter intersects the
	left box (and thus has to intersect the right one as well). The difference
	between the classes of polygons generated consists of those polygons whose
	bottom perimeter intersects the right box, but not the left.
\end{proof}

\begin{con}
	Considering the polygons described in the above proposition, there exists
	some perimeter preserving bijection\footnote{That is, the horiztonal and
	vertical perimeters remain constant.} between those polygons whose bottom
	perimeter intersects the right box, but not the left, and those whose bottom
	perimeter intersects both boxes.
\end{con}

\begin{cor}
	The generating function for 2-staircase (resp. 2-unimodal) polygons with
	vertical indents on the left side is
	\begin{equation} \label{Eq:2-unimodal.same}
		\frac{y^3}{2} \dd y \T^2 \dd y \frac{\mathcal{Q}}{y},
	\end{equation}
	where $\T = (\SP/y)^2$ and $\mathcal{Q} = \SP/y$ (resp. $\UP/y$).
\end{cor}

\begin{rem}
	In Appendix~\ref{s_m-same}, we extend the above result to $m$-unimodal
	polygons with all their indents being vertical and on the left side.
\end{rem}

\subsection{Wrapping} \label{s_wrapping}

The new star notation relates to the $E$ operator. When applied to a
directed path generating function, $E$ converts the perimeter generating
function into the half-perimeter generating function. Graphically, this is
equivalent to folding the directed path at the half-way point in each direction.
Indeed, the pyramid generating function is given Section~\ref{s_bimodal} as a
folded directed walk, and the unimodal generating function is given as a folded
pyramid.

This leads us to another new notion: that of ``wrapping''. To explain this, let
us consider 1-unimodal polygons formed by joining a bottom, staircase factor to
a top, unimodal factor and an indent. The total height of the polygon is given
by the sum of the heights of the top and bottom factors. The total width of the
polygon is measured by summing the width of the top factor to the width of the
bottom factor \emph{that lies to the left of the top factor}. (In the 1-unimodal
case, the generating function for the bottom factor is therefore $\T v(u/x)^d$.)

\begin{figure}
	\begin{center}
	\subfigure[The form of the top factor.]
				{\qquad\includegraphics[scale=0.8]{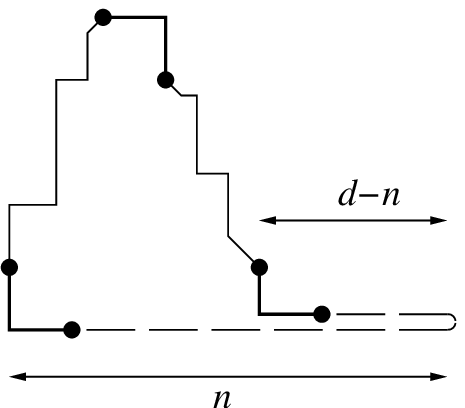}\qquad} \qquad
	\subfigure[The form of the resulting polygon.]
			{\raisebox{10pt}
			{\hspace{30pt}
				 \includegraphics[scale=0.4]{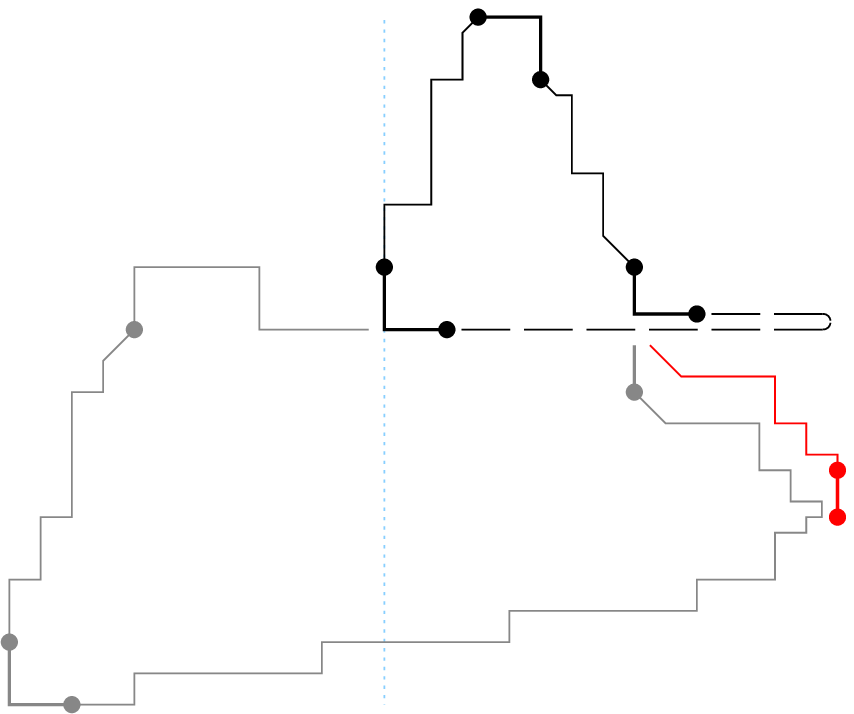}\hspace{40pt}}}
				 \\
	\subfigure[Flipping the top factor produces a polygon with the indent in the corner.]
				{\includegraphics[scale=0.4]{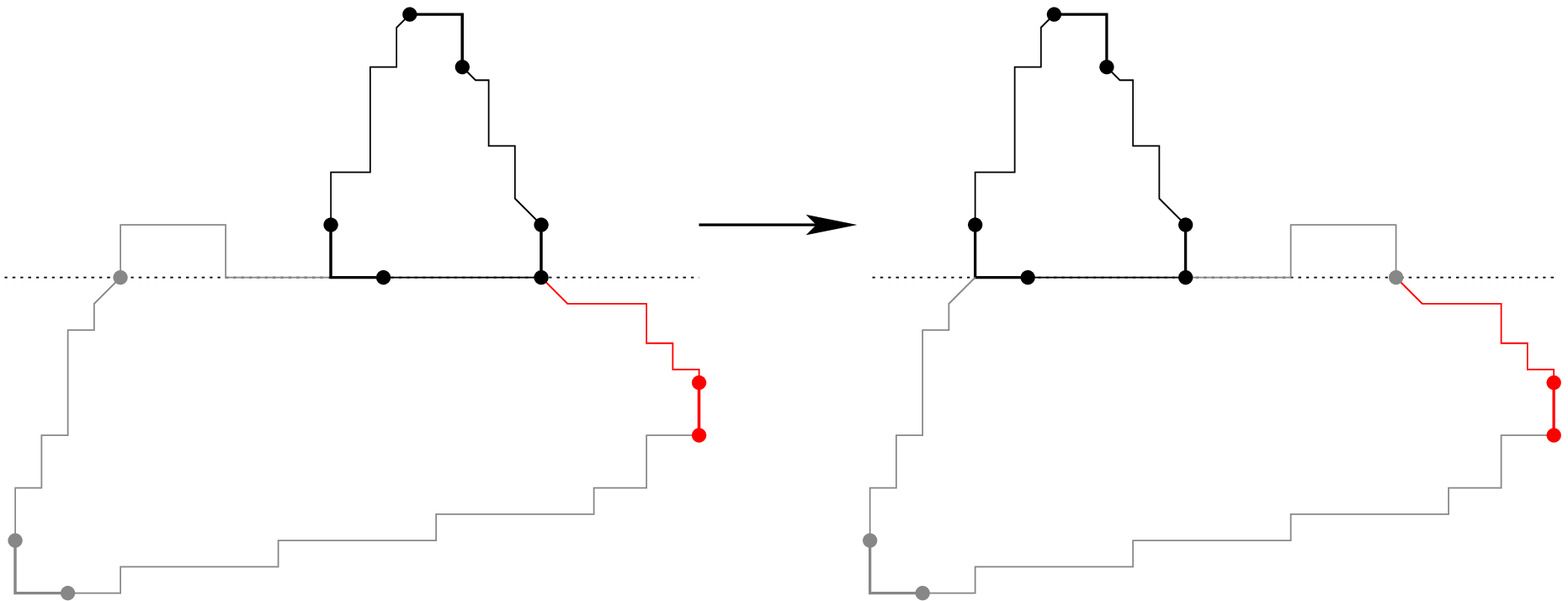}}
	\end{center}
	\caption{The form of the polygons when the top factor is narrower than $d$.}
	\label{Fig:wrapping}
\end{figure}

Now, if the join is of length $d$, then there are $d$ fixed horizontal steps in
each of the top and bottom factors that are identified, but are then removed and
do not form part of the polygon.
If there are fewer than $d$ horizontal steps in the rest of the top
unimodal factor, as shown in Figure~\ref{Fig:wrapping}(a), then the contribution
to the polygon is a pyramid of width $2n-d$, with a weight of $x^n$.
Importantly, although those fixed steps are not part of the polygon itself, they
do contribute to the weight. This is because the width of the polygon here is
given by the top factor. When the $E$ operator folds the walk, using the
standard inclusion-exclusion argument, the fixed steps are still identified with
those of the bottom factor, which is therefore folded as well. We call the
effect of one factor folding another along the join as ``wrapping''.

\begin{ntn}[Wrapping]
	The action of the $E$ operator is equivalent to the folding of a walk
	or polygon. When this forces the folding of a joined polygonal factor, as
	described abve, we refer to this as \DefEmph{wrapping}.
\end{ntn}

\subsection{The Hadamard product} \label{s_Hadamard}

\begin{figure}
	\begin{center}
		\includegraphics{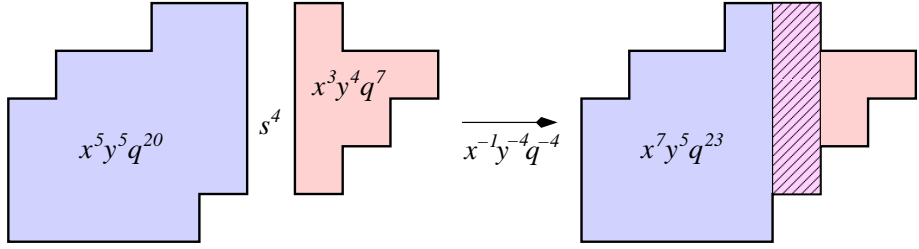}
	\end{center}
	\caption{An example of the action of the Hadamard product.}
	\label{Fig:2.hadamard}
\end{figure}

The Hadamard product is an operator which `joins' generating functions.  This
means that polygons with a side length enumerated according to a certain
parameter can be joined along those sides. For example, if we enumerate
staircase and stack polygons according to their right perimeter, total
perimeter and area, we can join them, making the neighbouring columns overlap.
Making the transformation $s\mapsto s/yq$ and dividing by $x$, so that the
overlapping column is not double-weighted, we form unimodal polygons, as in
Figure~\ref{Fig:2.hadamard}.

\begin{dfn}[Hadamard product / join]
	Let us consider two series, $f(t)=\sum_n f_n t^n$ and $g(t)=\sum_n g_n t^n$.
	We denote as $\odot_t$ the \emph{Hadamard product} with respect to $t$. We
	define
	\begin{equation} \label{Eq:hadamard}
		f(t) \odot_t g(t) = \sum_n f_n g_n t^n.
	\end{equation}
	The \emph{restricted Hadamard product} with respect to $t$ is defined as
	\begin{equation} \label{Eq:r-hadamard}
		f(t) \odot_t g(t) \arrowvert_{t=1} = \sum_n f_n g_n
		= \frac{1}{2\pi i} \oint f(t) g(1/t) \frac{\text{d}t}{t}.
	\end{equation}
	For notational convenience, we will refer to this as a \emph{Hadamard
	join} (over $t$), or simply a `join'.
\end{dfn}

\subsubsection{Hadamard arithmetic}

It is straight-forward to show the following properties of restricted Hadamard
products. (Alternatively, see \cite{Rechnitzer} for proofs.)

\paragraph{It is distributive:}
\begin{equation}
	f(t) \odot_t (g(t)+h(t)) = f(t) \odot_t g(t) + f(t) \odot_t h(t)
\end{equation}

\paragraph{It follows the product rule:}
\begin{equation}
	\dd s (f(s,t) \odot_t g(s,t)) = ( \dd s f(s,t)) \odot_t g(s,t)
	+ f(s,t) \odot_t ( \dd s g(s,t) )
\end{equation}

\paragraph{It evaluates simply at poles:}
\begin{gather}
	f(t) \odot_t \frac{1}{1-\alpha t} = f(\alpha)
\\	f(t) \odot_t \frac{t^k k!}{(1-\alpha t)^{k+1}} = \left(\dd t\right)^k f(t)
	\bigg|_{t=\alpha}
\end{gather}

Another useful observation is that the join with a function of the
form $(1-\alpha t)^{-k}$ can be simplified by the use of partial fractions. Thus,
\begin{equation}
	f(t) \odot_t \frac{1}{(1-\alpha t)^{k+1}}
	= \sum_{k=0}^n \frac{a^k}{k!} \left(\dd t\right)^k f(t) \bigg|_{t=\alpha}\\
\end{equation}

\noindent We can therefore simplify joins using partial fractions. For example,
we will often see joins of the form 
\begin{multline} \label{Eq:pfs}
	\frac{\alpha s^3}{(1-s)^3(1-\alpha s)} \\
=	\frac{\alpha}{2(1-\alpha)} \frac{2s^2}{(1-s)^3}
	- \frac{\alpha}{(1-\alpha)^2} \frac{s}{(1-s)^2}
	+ \frac{\alpha}{(1-\alpha)^3} \lr{\frac{1}{1-s}-\frac{1}{1-\alpha s}}.
\end{multline}
If, say, $\alpha=1/(1-v)=u/x$, we therefore have $\alpha/(1-\alpha)=-1/v$ and
$1/(1-\alpha)=-x/\SP$, which allows us to write the join as
\[
	\frac{vs^3}{(1-s)^3(1-v-s)} \odot_s f(s)
=	\frac{x^2(f(u/x)-f(1))}{\SP^2} - \bigg[ \frac{xf'(s)}{\SP}
	+ \frac {f''(s)}2 \bigg]_{s=1} .
\]

\subsection{Enumerating the building blocks} \label{s_blocks}

\subsubsection{Staircase and unimodal blocks}

When enumerating 2-convex polygons using Lin's factorisation, we will need the
generating functions for staircase and unimodal blocks according to their top
and bottom perimeters, as well as staircases according to their base and
right perimeter. The generating function for staircase polygons that counts the
left and right perimeters by $s$ and $t$ was given by Bousquet-Mélou
\cite{Bousquet:96}. She enumerated according to area, defining generalised
versions of modified $q$-Bessel functions. Setting $q=1$, we obtain
\[
	J_0(s) = \frac{(1-s)(1-sy)}{(1-s)(1-sy)+xs} \qquad\mbox{and}\qquad
	J_1(s,t) = \frac{xsJ_0(s)}{1-sty}.
\]
Denoting $T(s,t)\equiv T(x,y,s,t)$ as the generating function according to the
width (by $x$), height (by $y$), and right and left perimeters (by $s$ and $t$
respectively), we solve the recurrence relation given in \cite{Bousquet:96} and
obtain
\begin{equation} \label{Eq:Tst}
T(s,t) = x \bigg( J_0(s) \frac{sty}{1-sty} + (1-J_0(s)) \frac{vt}{1-vt} \bigg).
\end{equation}
For notational convenience, we define $\bar{T}(s,t)\equiv T(y,x,s,t)$.

We similarly define the generating function for unimodal polygons according to
their left and right perimeter and, through a slight abuse of notation, denote
it as $U(s,t)\equiv U(x,y,s,t)$. From \cite{Bousquet:96} one obtains
\begin{equation} \label{Eq:Ust}
U(s,t) = x \bigg( J_0(s) \frac{sty}{1-sty} \frac{(1-sy)((1-sy)-x)}{(1-sy)^2-x}
+ (1-J_0(s)) \frac{vt}{1-vt} (1+\SP Z) \bigg)
\end{equation}
with $U(1,t)=xvt(1+\SP Z)/(1-vt)$.\\

\begin{figure}
	\begin{center}
		\includegraphics[scale=0.5]{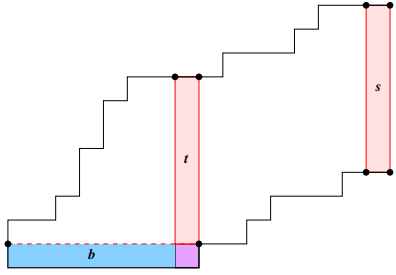}
	\end{center}
	\caption{Joining a staircase polygon to a Ferrers diagram.}
	\label{Fig:stair.side}
\end{figure}

Now, we will also need to enumerate staircases according to their base and right
perimeter, which we do by joining a Ferrers diagram to a staircase polygon, as
in Figure~\ref{Fig:stair.side}. Denoting the generating function
$F(b,s)\equiv F(x,y,b,s)$, we have
\begin{equation} \label{Eq:Fbs}
F(b,s)	= by (1-J_0(s)) \lr{\frac{(1-s)(1-x-sy)}{1-bx-sy} + \frac{u}{1-bu}},
\end{equation}
which is in a simple form to join over $b$.

%
%%
%%%
%%%% Pyramids with indents
%%%
%%
%

\subsubsection{Pyramids with indents}

\begin{figure}
	\begin{center}
	\subfigure[When the indent is not the highest factor.]{ 
		\includegraphics[scale=0.5]{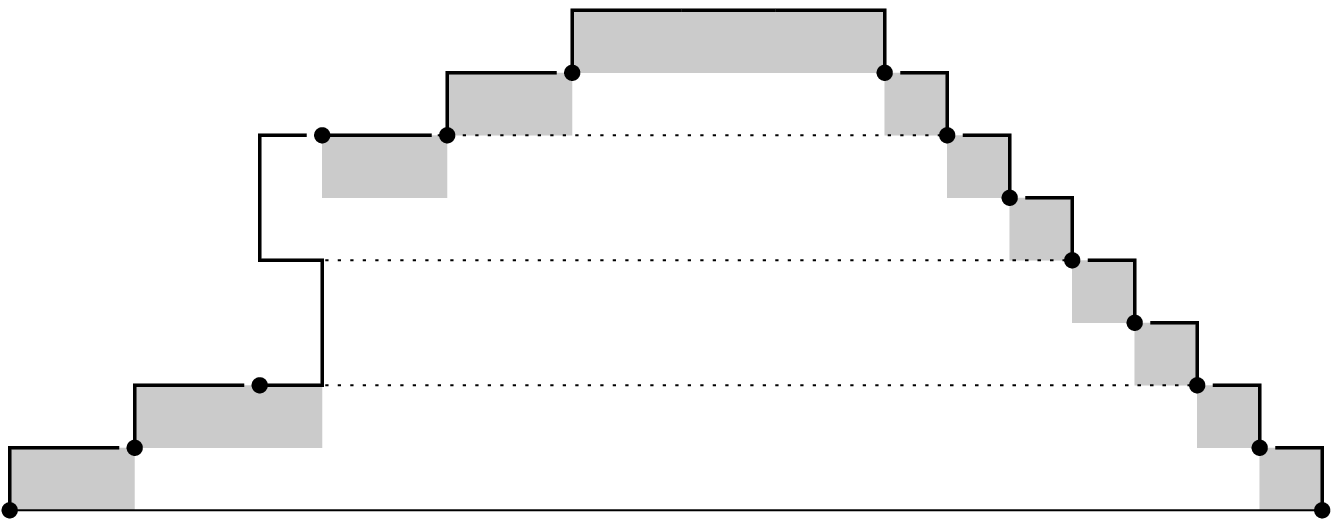}}\qquad
	\subfigure[When the indent is the highest factor.]{\quad 
	\includegraphics[scale=0.5]{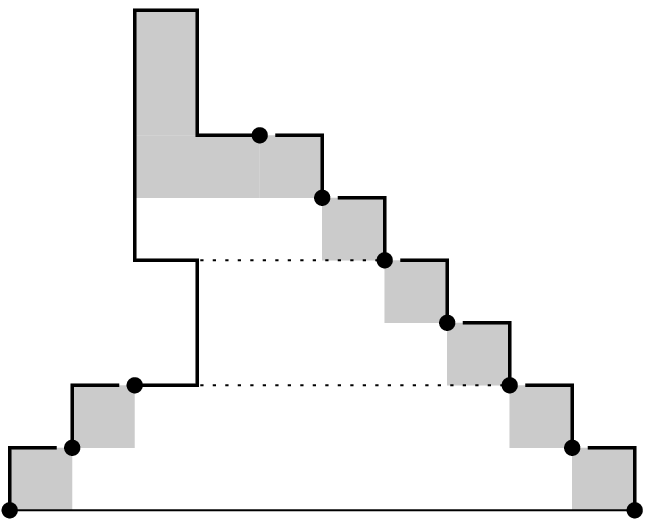} \quad}
	\end{center}
	\caption{The factorisation of a pyramid with a horizontal indent.}
	\label{Fig:2-u.pyramid}
\end{figure}

In the enumeration of 2-unimodal polygons, we will need to enumerate pyramid
factors with indents. A vertical indent we
can enumerate by simply taking the derivative to distinguish a height, and
multiplying by $x^2/(1-x)^2$, and so we only need to enumerate those with a
horizontal indent. In Figure~\ref{Fig:2-u.pyramid}, we illustrate the two cases
of such a pyramid: when the indent factor either does or does not touch the top
of the MBR. The generating function, which we denote $P'(x,y)$, is therefore
\begin{equation} \label{Eq:P'}
	P'(x,y) = \frac{xy}{(1-x)^2-y} \lr{\frac{y}{1-x-y}}^2
	\lr{ \frac{x(1-x)^2}{(1-x)^2-y} + \frac{xy}{1-y}}.
\end{equation}

%
%%
%%%
%%%% Convex: Top and bottom.
%%%
%%
%
\subsubsection{Enumerating convex polygons by their base and top perimeter}
\label{ss_convex_top_bottom}

Extending the inclusion-exclusion enumeration of convex polygons of
Bousquet-Mélou and Guttmann \cite{BMG:convex} \mm, the generating function for
convex polygons according to the left-perimeter is
\begin{multline} \label{Eq:Cs}
	\bar C(s) = \EE {xy(1-x)^2(1-y)^2}{2(1-x-y)^2} {
	\lr{ \frac{1-y}{x^2} - 1 } \frac{x^*s}{1-x^*s-y} } \\
	- Z U^2 \lr{\frac{1}{u}-1} \frac{us}{1-us}.
\end{multline}

\begin{figure}
	\begin{center}
		\subfigure[Non-pyramid unimodal polygons.]
		{\includegraphics[scale=0.6]{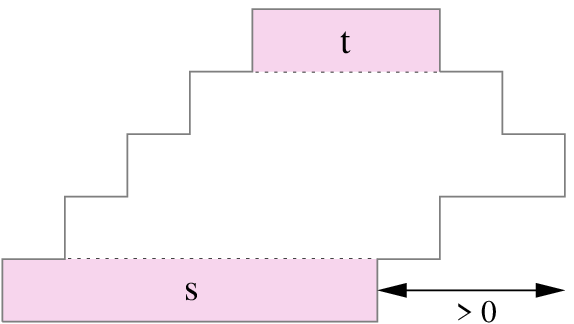}}
		\quad
		\subfigure[Non-unimodal convex polygons that touch left first.]
		{\quad\includegraphics[scale=0.6]{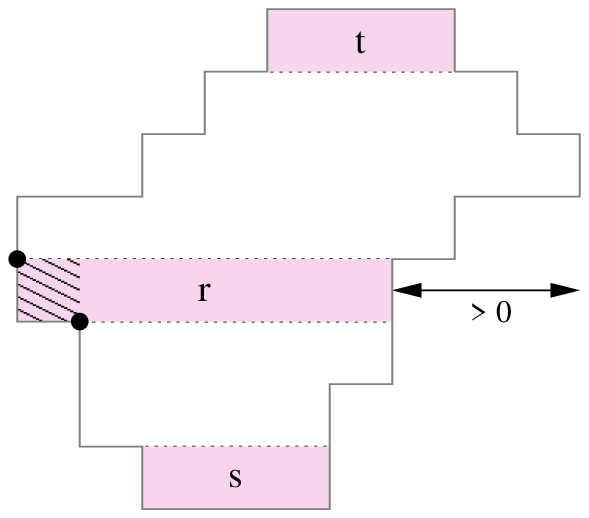}\quad}
		\quad
		\subfigure[Non-unimodal convex polygons that touch both sides at the
		same height.]
		{\qquad
		\includegraphics[scale=0.6]{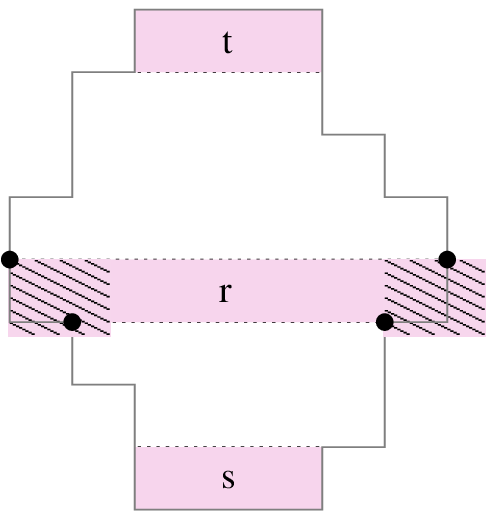}\qquad}
	\end{center}
	\caption{A decomposition of the form of convex polygons that allows their
			 enumeration by the base and top perimeter.}
	\label{Fig:2-c.tb}
\end{figure}

Now let us consider convex polygons according to their top and bottom
perimeters, whose generating function we denote $C(s,t)$.
We can decompose convex polygons into three classes: unimodal,
non-uni{\-}modal that touch the sides of the MBR at different heights, and those
that do not.
The first of these classes contains polygon that can be classified as rooted in
either of the bottom corners. If include both of these classes, we double-count
the pyramids. And so we enumerate unimodal polygons rooted in one corner, and
non-pyramid unimodal polygons rooted in the other. The latter class is depicted
in Figure~\ref{Fig:2-c.tb}(a), and the total generating function is
\(
	2 \bar U(s,t) - P(s,t).
\)

The next class are convex, but not unimodal, and include only those that touch
the side of the MBR at different heights. We assume, without loss of generality,
that they touch the left side first, as per the form shown in
Figure~\ref{Fig:2-c.tb}(b). We form these polygons by joining non-pyramid
unimodal polygons (as shown in part (a) of the figure) with upside-down
pyramids, overlapping the bottom rows, except for a single cell on the left. The
latter conditions determines that the polygon touches the right side at a
greater height. The generating function for this case is given by
\[
	\frac{2}{rs} \lr{\bar U(rs,t) - P(rs,t)} \odot_r
	\Ex y{rs}{1-rs}{\cdot \frac{y}{1-r-y}}.
\]

Finally, the case where the polygon touches each side at the same height
forms a diamond-like polygon, as shown in Figure~\ref{Fig:2-c.tb}(c). Similar to
the previous result, the generating function is simply
\[
	\frac{1}{(rs)^2} \lr{P(rs,t) - \frac{xrsty}{1-y}} \odot_r
	\Ex y{rs}{1-rs}{\cdot \frac{y}{1-r-y}}.
\]

The three cases then sum to give the desired generating function:
\begin{multline} \label{Eq:Cst}
	C(s,t)= \frac{1}{(1-s(1+(1-s)x-y)) (1-t(1+(1-t)x-y)) \Delta^{3/2}} \times \\
	\bigg( -stx^2y^2(1-x+y-s(1-x-y))(1-x+y-t(1-x-y)) \\ 
	 + \frac{stxy A(s,t)} {(1-stx) ((1-sx)^2-y) ((1-tx)^2-y) \sqrt{\Delta}} \bigg),
\end{multline}
{\small
\noindent
\begin{minipage}[t]{\linewidth}$
\mbox{\normalsize{where }} A(s,t) =
(1-s) (1-t) (1-x)^4 (1-s x)^3 (1-t x)^3 + (1-x)^2 (1-s x)^2 (1-t x)^2 
   (-5 - x - 2 x^2 - t^2 x (5 + 2 x + x^2) + t (6 + 4 x + 5 x^2 + x^3) + 
    s^2 x (-5 - 2 x - x^2 + t^2 x (-4 - 5 x + x^2) + t (6 + 4 x + 7 x^2 - x^3)) + 
    s (6 + 4 x + 5 x^2 + x^3 + t^2 x (6 + 4 x + 7 x^2 - x^3) + t (-7 - 9 x - 8 x^2 - 9 x^3 + x^4))) y + 
  (10 + 3 x - 3 x^2 + 5 x^3 + x^4 + 2 t^4 x^3 (5 + x^2 + 2 x^3) - t^3 x^2 (35 + 13 x^2 + 10 x^3 + 6 x^4) + 
    t^2 x (40 + 12 x + 17 x^2 + 11 x^3 + 15 x^4 + x^5) - t (15 + 23 x + 5 x^2 + 8 x^3 + 10 x^4 + 3 x^5) + 
    s^4 x^3 (2 t^4 x^3 (3 + x + 5 x^2 - x^3) + 2 (5 + x^2 + 2 x^3) + 2 t^2 x (17 + 5 x + 9 x^2 + 19 x^3 - 2 x^4) + 
      t^3 x^2 (-26 - 2 x - 27 x^2 - 12 x^3 + 3 x^4) - t (15 + 16 x + 7 x^2 + 22 x^3 + 4 x^4)) + 
    s^3 x^2 (-35 - 13 x^2 - 10 x^3 - 6 x^4 + t^4 x^3 (-26 - 2 x - 27 x^2 - 12 x^3 + 3 x^4) + 
      t^3 x^2 (100 + 5 x + 79 x^2 + 63 x^3 + 13 x^4 - 4 x^5) + t^2 x (-123 - 40 x - 71 x^2 -
$\end{minipage}
\noindent
\begin{minipage}[t]{\linewidth}$
	  105 x^3 - 50 x^4 + 5 x^5) + 
      t (51 + 66 x + 23 x^2 + 71 x^3 + 38 x^4 + 7 x^5)) + s^2 x (40 + 12 x + 17 x^2 + 11 x^3 + 15 x^4 + x^5 + 
      2 t^4 x^3 (17 + 5 x + 9 x^2 + 19 x^3 - 2 x^4) + t^2 x (145 + 92 x + 93 x^2 + 113 x^3 + 130 x^4 + 3 x^5) + 
      t^3 x^2 (-123 - 40 x - 71 x^2 - 105 x^3 - 50 x^4 + 5 x^5) - t (57 + 102 x + 43 x^2 + 79 x^3 + 78 x^4 + 23 x^5 + 2 x^6)) - 
    s (15 + 23 x + 5 x^2 + 8 x^3 + 10 x^4 + 3 x^5 + t^4 x^3 (15 + 16 x + 7 x^2 + 22 x^3 + 4 x^4) - 
      t^3 x^2 (51 + 66 x + 23 x^2 + 71 x^3 + 38 x^4 + 7 x^5) + t^2 x (57 + 102 x + 43 x^2 + 79 x^3 + 78 x^4 + 23 x^5 + 2 x^6) - 
      t (21 + 69 x + 44 x^2 + 39 x^3 + 56 x^4 + 20 x^5 + 7 x^6))) y^2 + 
  (-10 - 22 x - x^2 - 3 x^3 - 2 t^4 x^3 (5 + 4 x + 3 x^2) + t^3 x^2 (40 + 34 x + 29 x^2 + 9 x^3) - 
    t^2 x (50 + 62 x + 49 x^2 + 25 x^3 + 2 x^4) + t (20 + 52 x + 38 x^2 + 20 x^3 + 5 x^4 + x^5) + 
    s^4 x^3 (t^4 x^3 (-4 - x + x^2) - 2 (5 + 4 x + 3 x^2) + t^3 x^2 (24 + 10 x + 9 x^2 - 3 x^3) + 
      t^2 x (-36 - 28 x - 31 x^2 + 3 x^3) + t (20 + 26 x + 23 x^2 + 11 x^3)) + 
    s^3 x^2 (40 + 34 x + 29 x^2 + 9 x^3 + t^4 x^3 (24 + 10 x + 9 x^2 - 3 x^3) + t^2 x (152 + 130 x + 151 x^2 + 62 x^3 - 7 x^4) + 
      t^3 x^2 (-110 - 61 x - 70 x^2 - 13 x^3 + 6 x^4) - t (75 + 118 x + 109 x^2 + 72 x^3 + 18 x^4)) + 
    s^2 x (-50 - 62 x - 49 x^2 - 25 x^3 - 2 x^4 + t^4 x^3 (-36 - 28 x - 31 x^2 + 3 x^3) + 
      t^3 x^2 (152 + 130 x + 151 x^2 + 62 x^3 - 7 x^4) - t^2 x (200 + 234 x + 261 x^2 + 184 x^3 + 9 x^4) + 
      t (90 + 190 x + 192 x^2 + 147 x^3 + 54 x^4 + 7 x^5)) + s (20 + 52 x + 38 x^2 + 20 x^3 + 5 x^4 + x^5 + 
      t^4 x^3 (20 + 26 x + 23 x^2 + 11 x^3) - t^3 x^2 (75 + 118 x + 109 x^2 + 72 x^3 + 18 x^4) + 
      t^2 x (90 + 190 x + 192 x^2 + 147 x^3 + 54 x^4 + 7 x^5) - t (35 + 130 x + 147 x^2 + 121 x^3 + 47 x^4 + 23 x^5 + x^6))) y^3 + 
  (5 + 18 x + x^2 - x^3 + t^4 x^3 (5 + 4 x) - t^3 x^2 (25 + 27 x + 9 x^2) + t^2 x (35 + 59 x + 29 x^2 + 3 x^3) - 
    t (15 + 53 x + 32 x^2 + 6 x^3 + 3 x^4) + s^4 x^3 (5 + 4 x + t^4 x^3 + t^3 x^2 (-11 - 5 x + x^2) + t^2 x (19 + 15 x + 2 x^2) - 
      t (15 + 17 x + 11 x^2)) + s^3 x^2 (-25 - 27 x - 9 x^2 + t^4 x^3 (-11 - 5 x + x^2) + t^3 x^2 (65 + 45 x + 12 x^2 - 4 x^3) + 
      t^2 x (-103 - 102 x - 56 x^2 + x^3) + t (65 + 98 x + 80 x^2 + 23 x^3)) + 
    s^2 x (35 + 59 x + 29 x^2 + 3 x^3 + t^4 x^3 (19 + 15 x + 2 x^2) + t^3 x^2 (-103 - 102 x - 56 x^2 + x^3) + 
      t^2 x (155 + 210 x + 168 x^2
$\end{minipage}
\noindent
\begin{minipage}[t]{\linewidth}$
	  + 29 x^3) - t (85 + 180 x + 179 x^2 + 81 x^3 + 11 x^4)) - 
    s (15 + 53 x + 32 x^2 + 6 x^3 + 3 x^4 + t^4 x^3 (15 + 17 x + 11 x^2) - t^3 x^2 (65 + 98 x + 80 x^2 + 23 x^3) + 
      t^2 x (85 + 180 x + 179 x^2 + 81 x^3 + 11 x^4) - t (35 + 140 x + 163 x^2 + 93 x^3 + 35 x^4 + 4 x^5))) y^4 + 
  (-1 - 3 x + 2 x^2 - t^4 x^3 + 2 t^3 x^2 (4 + 3 x) - t^2 x (13 + 21 x + 4 x^2) + t (6 + 26 x + 9 x^2 + 3 x^3) + 
    s^4 x^3 (-1 + 2 t^3 x^2 - 2 t^2 x (2 + x) + t (6 + 4 x)) + s^3 x^2 (8 + 6 x + 2 t^4 x^3 + t^3 x^2 (-19 - 8 x + x^2) + 
      2 t^2 x (18 + 15 x + 2 x^2) - t (33 + 39 x + 16 x^2)) + s^2 x (-13 - 21 x - 4 x^2 - 2 t^4 x^3 (2 + x) + 
      2 t^3 x^2 (18 + 15 x + 2 x^2) - t^2 x (64 + 83 x + 33 x^2) + t (48 + 90 x + 63 x^2 + 11 x^3)) + 
    s (6 + 26 x + 9 x^2 + 3 x^3 + 2 t^4 x^3 (3 + 2 x) - t^3 x^2 (33 + 39 x + 16 x^2) + t^2 x (48 + 90 x + 63 x^2 + 11 x^3) - 
      t (21 + 87 x + 83 x^2 + 33 x^3 + 6 x^4))) y^5 + (-x - t^3 x^2 - s^4 t x^3 + 2 t^2 x (1 + x) - t (1 + 5 x + x^2) + 
    s^3 x^2 (-1 + 2 t^3 x^2 - t^2 x (5 + 2 x) + t (9 + 6 x)) + s^2 x (2 (1 + x) - t^3 x^2 (5 + 2 x) + t^2 x (11 + 12 x) - 
      t (15 + 22 x + 7 x^2)) - s (1 + 5 x + x^2 + t^4 x^3 - 3 t^3 x^2 (3 + 2 x) + t^2 x (15 + 22 x + 7 x^2) - 
      t (7 + 29 x + 18 x^2 + 4 x^3))) y^6 - s t (1 - 2 (-2 + s + t) x + (-2 s + s^2 + (1-t)^2) x^2) y^7
$\end{minipage}
}

\section{Enumerating 1-convex polygons using Lin's factorisation and the
inclusion-exclusion principle} \label{s_1-convex}

In section we use Lin's factorisation before the inclusion-exclusion principle
to rederive the 1-unimodal and 1-convex polygon generating functions. (The
1-staircase case was enumerated in the previous section by distinguishing a
vertical step.) This contrasts greatly with the more complicated derivation in
\cite{wrgj} that generated all polygons as a single case of inclusion-exclusion.
These results will serve as an example for the general approach we adopt later
in this paper.

\subsection{Enumerating 1-unimodal polygons}

\subsubsection{Indent on left.}
This case was also enumerated in the previous section. Adopting the same
approach \mm, the generating function is
\(
	\T \cdot y \dd y \UP.
\)
This confirms the following generating function, obtained via inclusion-exclusion:
\begin{multline} \label{1-unimodal-sum}
	\sum_{d\geq 1} \frac{u^d v}{x^d} \lr{
		\EE{x^dy(1-x)(1-y)}{1-x-y}{} - u^d v \EE{xy}{1-x-y}{} } \\
	= v \EE{x y^2}{1-x-y} {\sum_{d\geq 1} \left(\frac{x}{1-v^*}\right)^d}
	- v^2 {\sum_{d\geq 1} \lr{\frac u{1-v}}^d} \EE{x y}{1-x-y}{} \\
	= v x y \EE 1{(1-x-v^*)(1-x-y)}{} - u v^2 Z \cdot \frac {2xy}\Delta.
\end{multline}
Above, the inclusion term includes polygons of width less than $d$, such that
the polygon wraps, to form those shown in Figure~\ref{Fig:wrapping}(b).
We can see how this allows us to enumerate both pyramid and unimodal top factors
using the inclusion-exclusion method, without restricting the direction of the
path as it passes through Lin's factorisation line.

\subsubsection{Indent in corner.} \label{ss_unimodal_corner}

Now, in the previous enumeration, when the top factor is a pyramid the indent
occurs higher than all of its right edge steps. In other words, it touches the
right edge below the
factorisation line, as depicted in Figure~\ref{Fig:wrapping}(b). This means
that we can swap the order of the indent and the top pyramid factor, as per part
(c) of the figure, creating an indent in the top-right corner. Noting that
we must include the case where the indent factor touches the right edge of the
MBR and enumerating them as corner-staircase polygons from \cite{wrgj}, we
obtain \[
	\T y \SP Z \lr{ \frac{1}{1-x} + uZ }.
\]

One can check that this is the same result as that which is obtained by adopting
the inclusion-exclusion method. To enumerate the case where the
indent is not to the right of the bottom factor, and taking the first term of
equation~\ref{1-unimodal-sum} (n.b. there is no possibility of intersection),
we exclude the cases of width $\geq d$ by writing
\begin{align*}
\frac{u^2}{(1-u)^2}& \sum_{d\geq 1} \frac{u^d}{x^d} v y E\bigg[\frac{x^{d+1}}
{1-x-y} \bigg(1-\bigg(\frac{x}{1-y}\bigg)^{d-1} \bigg) \bigg] \\ &=
\frac{\SP^2}{1-u} \sum_{d\geq 1} u^d E\bigg[\frac{x^{-(d-1)}-(1-y)^{-(d-1)}}
{1-x-y} \bigg]\\
&=  \frac{\SP^2}{1-u} E\bigg[ \bigg(\frac{u^*}{1-u^*/x}-\frac{u^*}{1-u^*/(1-y)}
\bigg) \bigg/ (1-x-y) \bigg].
\end{align*}
This is equivalent to distinguishing a vertical step at a height above the right
edge of the polygon.

\subsubsection{Indent on bottom.} \label{ss_unimodal_bottom}

Using our above methodology \mm, the generating function for the case where the
indent is on the bottom edge of the polygon is
\begin{equation}
 x^2yv E_{\{y\}}\bigg[ \frac{u^*}{(1-u^*-y)(1-x-y)^2}
 \left( \frac{(1-y)^2}{(1-y)^2-x} + \frac{x}{1-x} \right)
 \bigg] - \frac{2xv\SP^3Z}{y\Delta}.
\end{equation}
Whenever indents are attached to the base of a unimodal factor, we will see
generating functions of this form.

\subsection{Enumerating 1-convex polygons}

In the case of 1-convex polygons, not having to root the polygon in one corner,
like we did for 1-unimodal polygons, returns certain symmetries. This means that
we no longer have to
enumerate three asymmetric cases; we only have to do the case where the indent
is on the top-left side of the polygon. We note, however, that when the top
factor is of height one, there is a vertical symmetry. And so, when the indent
is at the top of the polygon, it will be double counted. This results in the
term $(1+y/(1-x))$ appearing in the generating function.

Following the above methodology \emph{mutatis mutandis}, the generating function
for 1-convex polygons is
\begin{gather*} \hspace{-20pt} \EE{x(1-x)y^*}{(1-x)^2-y^*} {
	\lr{1+\frac{y^*}{(1-x)^2}}\lr{\frac{x^*}{1-x^*-y}}^2
	\lr{\frac{y^2(1-y)^2}{(1-y)^2-x^*}+\frac{x^*y^2}{1-x^*}} \frac{x}{1-x-y}
}\\ - \frac{4xyv}{\Delta} \EE{}{}{ \lr{\frac{x^*y}{1-x^*-y}}^2
	\lr{\frac{(1-y)^2}{(1-y)^2-x^*}+\frac{x^*}{1-x^*}} \frac{u^*}{1-u^*-y}
}\\ - \frac{2xyu^2v}{(1-u)^2\Delta} \EE{x(1-x)y^*}{(1-x)^2-y^*}{
	\lr{1+\frac{y^*}{(1-x)^2}} \frac{x}{1-x-v^*}
}\\ + 2 v \SP Z \lr{ \frac{2x\SP}{\Delta}}^2
- \frac{2\SP^3Z}{y}\left(\frac{1}{1-x}+uZ\right) \left(1+\frac{v}{1-u}\right)Z\UP,
\end{gather*}
where the last term comes from the possibility of intersection along the NW
diagonal, factoring into a unimodal polygon and a 1-unimodal SAP with the indent
in the corner.

\section{Enumerating 2-staircase polygons} \label{s_2-staircase}

In Section~\ref{s_distinguished}, we enumerated bimodal 2-staircase polygons.
We will now enumerate the four cases when the indents are distinct.

\subsection{Case 1: indents in the same direction, on the same edge}

In Section~\ref{s_distinguished}, we also enumerated 2-staircase polygons with
both indents on the left edge. The case when the indents are at the same height
follows from the 1-staircase result. These give the desired result.

\subsection{Case 2: indents in the same direction, on opposite edges}

If there is one indent each on the top and bottom edges, then there are
three different cases to enumerate: if the top one is above, level or below the
bottom one. If they are level, the top one can be to the left (enumerated by
$\T^2 v \SP Z$) or the right ($(x/(1-x))^3(\SP/(1-u))^2$). The other cases are
enumerated by
\[
	\frac{wvs}{(1-w)(1-s)(1-v-s)} \odot_s \bar{T}(s,t) \odot_t
	\frac{wvt}{(1-w)(1-t)(1-v-t)},
\]
where $w = u$ if the top indent is above the bottom indent, and $w = s$ or $t$
if it is below. To evaluate such expressions, we refer to the factorisation of
joins using partial fractions in \eqref{Eq:pfs} and let $\alpha=1/(1-v)=u/x$.

\subsection{Case 3: indents in different directions, on the same edge}

\begin{figure}
	\begin{center}
	\subfigure[The indents form a convex factor.]{
		\includegraphics[scale=0.4]{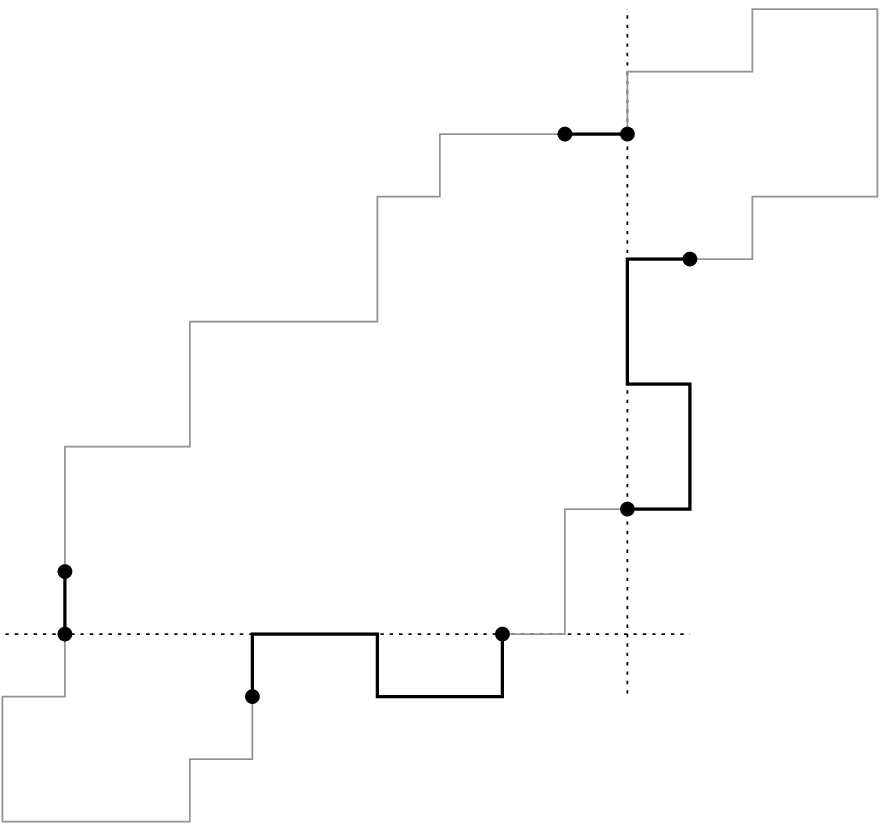}\qquad
		\includegraphics[scale=0.4]{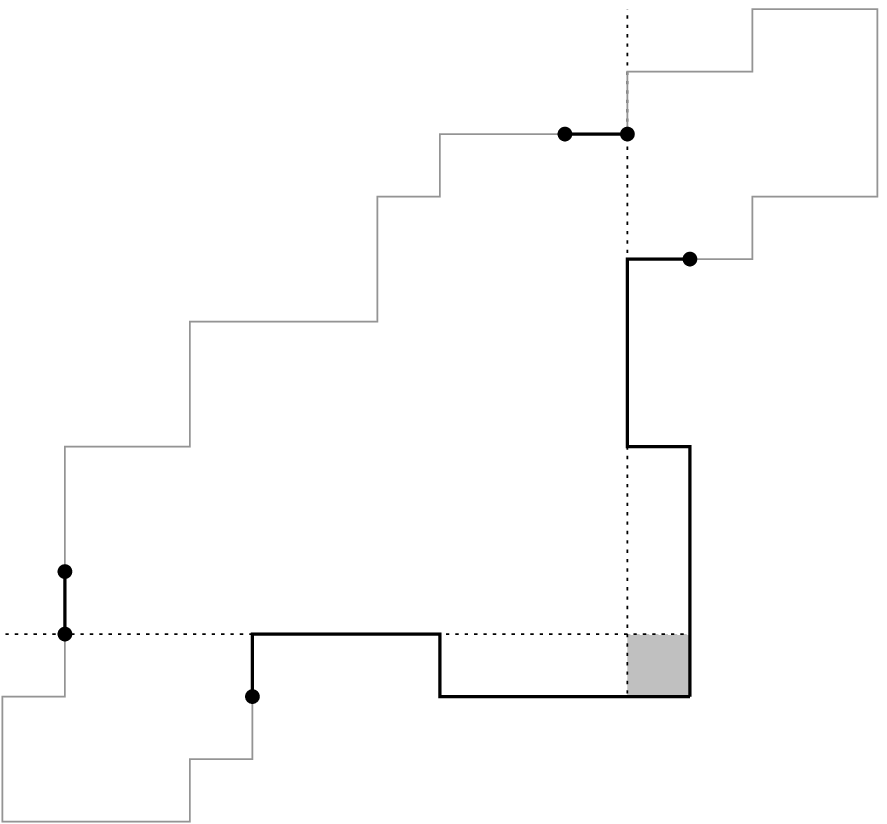}} \qquad
	\subfigure[The indents form a concave factor.]{
		\includegraphics[scale=0.4]{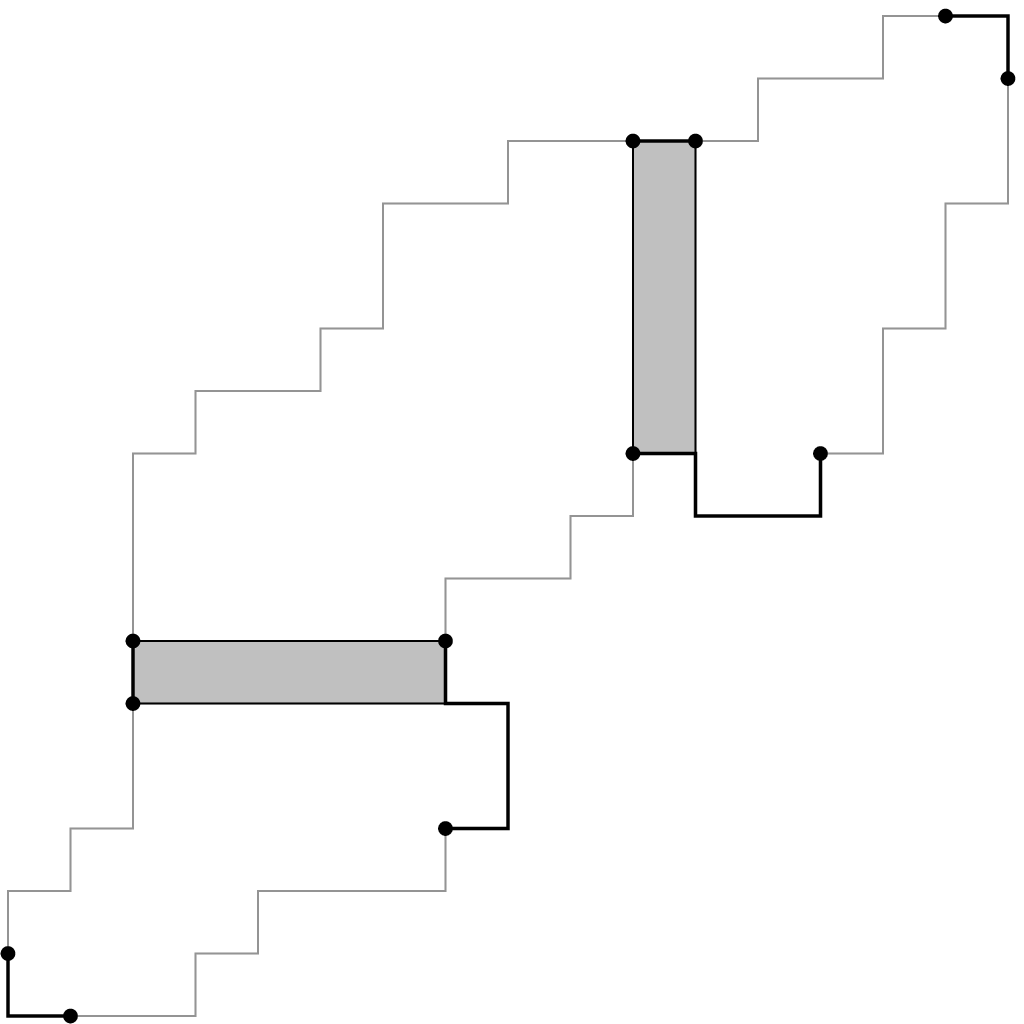}}
	\end{center}
	\caption[The relative position of indents on adjacent sides can cause
			 a convex or concave shape.]
			 {The relative position of indents on adjacent sides can cause
			 a convex or concave shape (in the usual non-lattice sense of
			 convexity).}
	\label{Fig:2-stair.adjacent}
\end{figure}

When the indents of a 2-staircase polygon are on the same edge, they can form
a middle factor that is either convex or concave (in the usual Euclidean sense
of the words). This is shown in Figure~\ref{Fig:2-stair.adjacent}. The
generating functions can be obtained by using the techniques of
Section~\ref{s_methodology} \mm.

\subsection{Case 4: indents in different directions, on opposite edges}

\begin{figure}
	\begin{center}
	\subfigure[The symmetric case, with the left indent above.]{
		\includegraphics[scale=0.5]{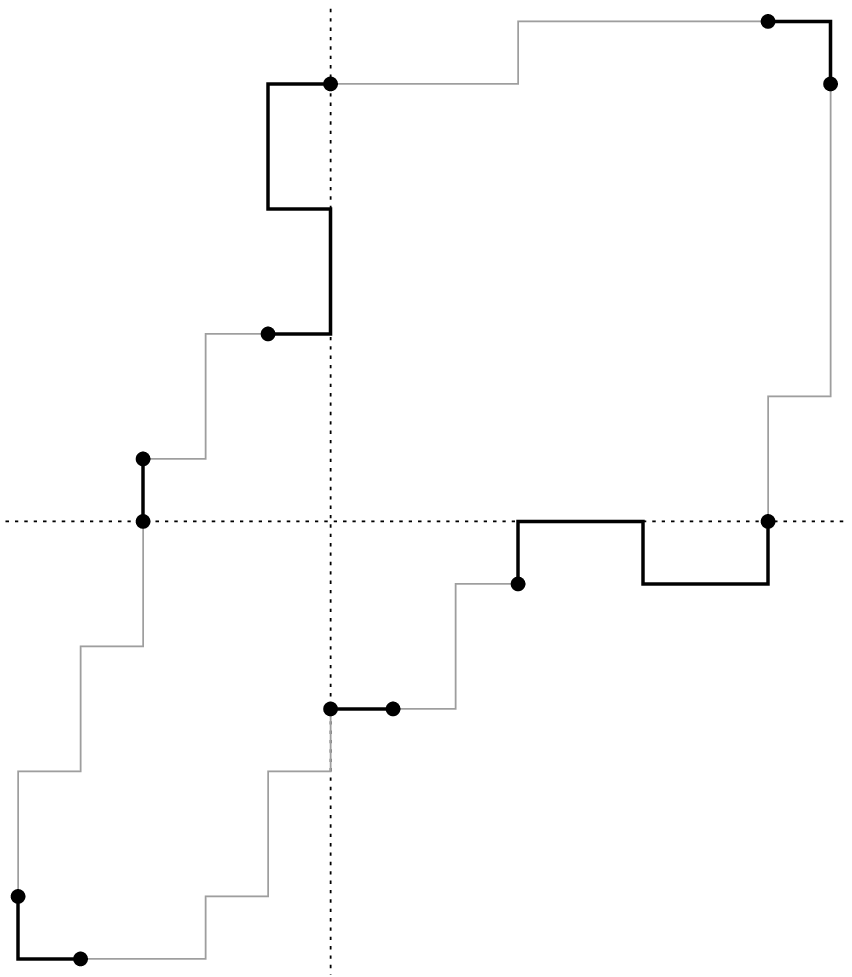}} \qquad
	\subfigure[The left indent is below the bottom indent.]{
		\includegraphics[scale=0.5]{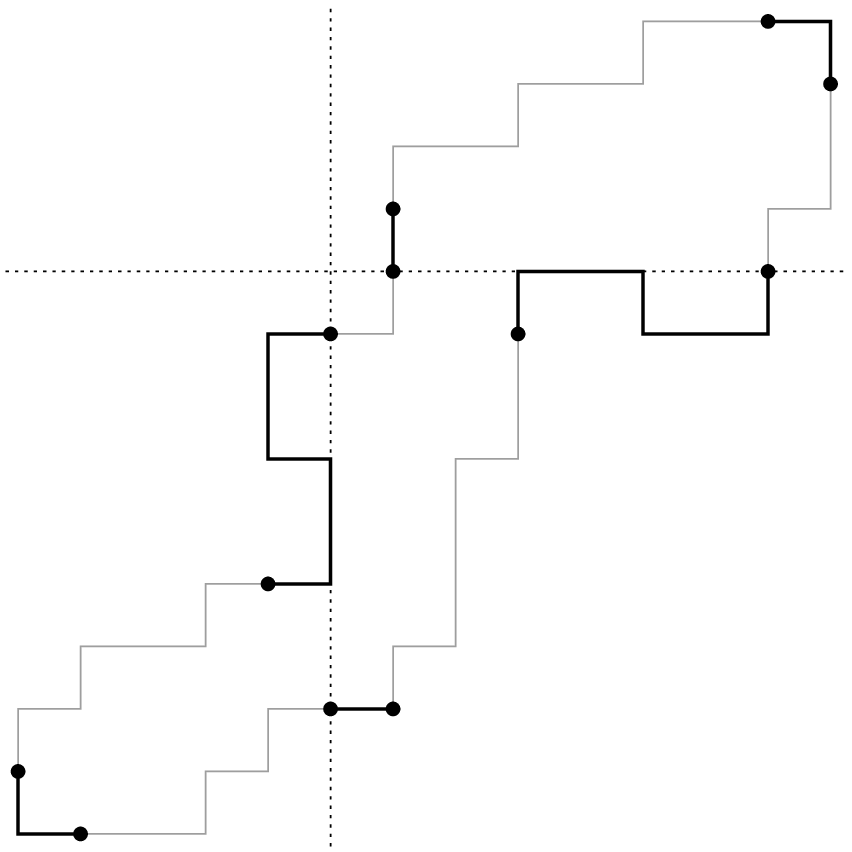}}
	\subfigure[The left indent is next-to the bottom indent.]{
		\includegraphics[scale=0.5]{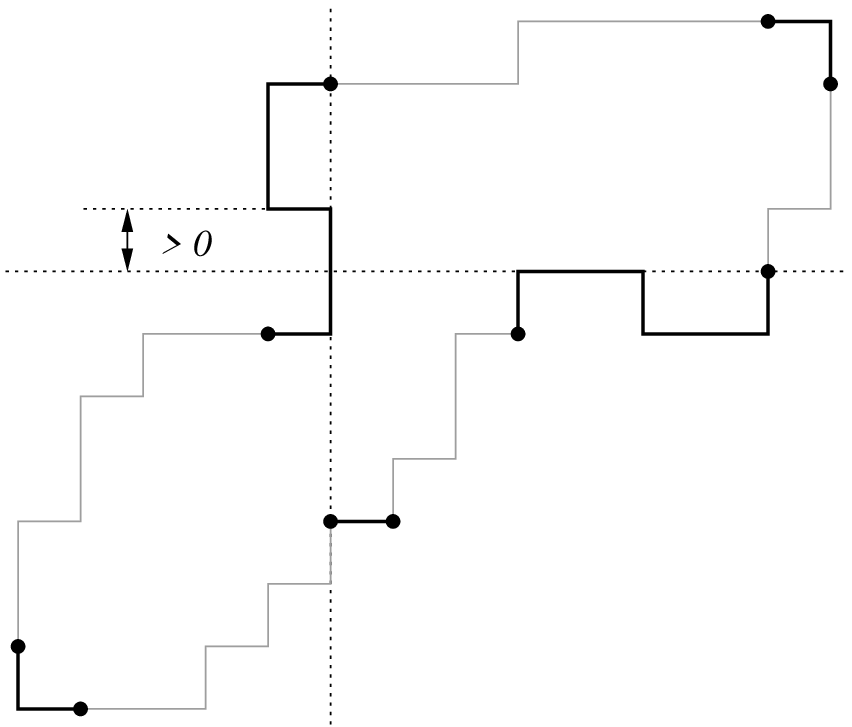}  \qquad
		\includegraphics[scale=0.5]{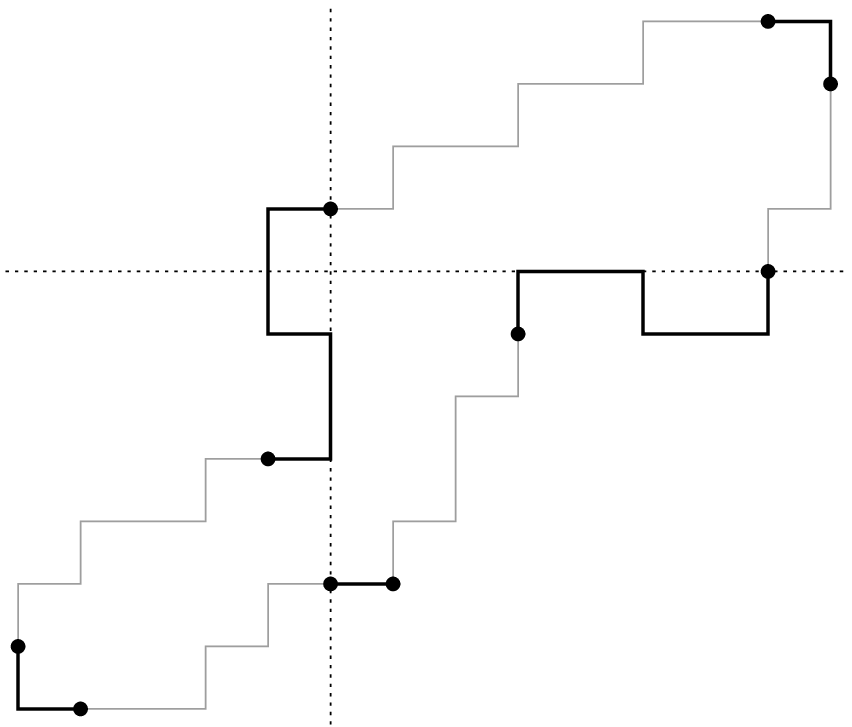}}
	\end{center}
	\caption{The relative position of indents on adjacent sides in opposite
	halves puts it in one of three cases.}
	\label{Fig:2-stair.axis}
\end{figure}

When the indents are the bottom and left, there are three possible cases in
terms of their relative position. These are shown in
Figure~\ref{Fig:2-stair.axis}. The first case is symmetrical about the $x=y$
diagonal. The other two cases have an equivalent class after reflection.

The generating function for the first case can be given on inspection, from the
corner-staircase generating function from \cite{wrgj}.
It is $\T \bar\T \SP^2 Z^2.$
Similarly, for both cases where the top indent is next-to the bottom indent, as
in Figure~\ref{Fig:2-stair.axis}(c), one can write down their generating
functions on inspection (which, incidentally, are equal).
Finally, when the top indent is below the bottom one, as in
Figure~\ref{Fig:2-stair.axis}(b), the generating function is given by simply
joining the staircase and indent factors.

\subsection{The generating function for 2-staircase polygons }

Combining all of the above cases gives the following generating function for
2-staircase polygons:
\begin{equation}\label{Eq:2-staircase}
	\frac{-A^{(s)}_2}{2x^4y^4(1-x)^3(1-y)^3(1-x-y)} + \frac{-B^{(s)}_2}{2x^4y^4(1-x)(1-y)\Delta^{3/2}}
\end{equation}
where {\small
\begin{align*}
A^{(s)}_2 =& (1-x)^9 x^4
- 3 (1-x)^7 x^4 (3 - 6 x + x^2) y
+ (1-x)^5 x^2 (4 - 8 x + 33 x^2 \\& - 109 x^3 + 112 x^4 - 47 x^5 + 3 x^6) y^2
- (1-x)^3 x^2 (28 - 96 x + 169 x^2  \\&- 321 x^3 + 471 x^4 - 395 x^5 + 184 x^6 - 39 x^7 + x^8) y^3
- (1-x) (-1 + 8 x \\& - 105 x^2 + 436 x^3 - 890 x^4 + 1198 x^5 - 1292 x^6 + 1120 x^7 - 701 x^8 + 301 x^9 \\& - 81 x^{10} + 9 x^{11}) y^4
+ (-9 + 81 x - 394 x^2 + 1144 x^3 - 2088 x^4 + 2562 x^5 \\& - 2238 x^6 + 1410 x^7 - 588 x^8 + 139 x^9 - 20 x^{10} + 2 x^{11} + x^{12}) y^5
+ (36 \\& - 318 x + 1087 x^2 - 2037 x^3 + 2490 x^4 - 2238 x^5 + 1530 x^6 - 738 x^7 \\& + 192 x^8 - 6 x^9 - 3 x^{10} - x^{11}) y^6
+ (-84 + 714 x - 2071 x^2 + 2940 x^3 \\& - 2412 x^4 + 1410 x^5 - 738 x^6 + 312 x^7 - 70 x^8 + 6 x^9 - x^{10}) y^7
+ (126 \\& - 1008 x + 2621 x^2 - 3103 x^3 + 1821 x^4 - 588 x^5 + 192 x^6 - 70 x^7 + 6 x^8 \\& + x^9) y^8
- (1-x)^3 (126 - 546 x + 167 x^2 + 18 x^3 + 9 x^4 + x^5) y^9 \\&
+ (1-x)^3 (84 - 294 x + 35 x^2 + 6 x^3 + x^4) y^{10}
+ (1-x)^3 (-36 + 90 x \\& + x^2 + x^3) y^{11}
- (1-x)^3 (-9 + 12 x + x^2) y^{12}
- (1-x)^3 y^{13} \\
B^{(s)}_2 = & (1-x)^9 x^4
- (1-x)^8 (9-x) x^4 y
- (1-x)^5 x^2 (-4 + 8 x - 33 x^2 + 51 x^3 \\& + 13 x^4 + 7 x^5) y^2
- (1-x)^3 x^2 (28 - 64 x + 97 x^2 - 147 x^3 + 10 x^4 - x^5 \\& - 34 x^6 + x^7) y^3
+ (1-x) (1 - 8 x + 105 x^2 - 268 x^3 + 326 x^4 - 338 x^5 \\& + 194 x^6 - 55 x^7 + 55 x^8 - 49 x^9 + 5 x^{10}) y^4
+ (-9 + 73 x - 336 x^2 + 658 x^3 \\& - 664 x^4 + 478 x^5 - 232 x^6 + 35 x^7 + 35 x^8 - 16 x^9 + 10 x^{10}) y^5
+ (36 - 260 x \\& + 672 x^2 - 806 x^3 + 532 x^4 - 232 x^5 + 64 x^6 - 11 x^7 - 33 x^8 - 10 x^9) y^6 \\&
+ (-84 + 532 x - 826 x^2 + 569 x^3 - 249 x^4 + 35 x^5 - 11 x^6 + 48 x^7 + 4 x^8) y^7 \\&
+ (126 - 686 x + 588 x^2 - 146 x^3 + 110 x^4 + 35 x^5 - 33 x^6 + 4 x^7) y^8 \\&
+ 2 (1-x) (-63 + 224 x + 110 x^2 + 65 x^3 + 13 x^4 + 5 x^5) y^9
- 2 (1-x) (-42 \\& + 112 x + 84 x^2 + 32 x^3 + 5 x^4) y^{10}
+ (1-x) (-36 + 64 x + 42 x^2 + 5 x^3) y^{11} \\&
+ (1-x)^2 (9 + x) y^{12}
- (1-x) y^{13}
\end{align*}
}

\section{Enumerating 2-unimodal polygons} \label{s_2-unimodal}

\subsection{Bimodal 2-unimodal polygons}

\subsubsection{Indent on the left}

Following Section~\ref{s_distinguished} \mm, we obtain the following generating
function by distinguishing a vertical step:
\begin{equation*}
	\T_2 y^4 \dd y \lr{\frac{\UP - xy/(1-x) - x(1+x)y^2/(1-x)^3}{y^3}}.
\end{equation*}

\subsubsection{Indent in corner}

Following Section~\ref{ss_unimodal_corner} \mm, the generating
function for bimodal 2-unimodal polygons with their indent in the corner, such
that the bottom factor is further to the right than the indent factor, is
$\T \SP^2 Z (Z-1/(1-u))/(1-u)$.
When the indent factor is further to the right, the generating function of the
indent factor for base $g$ is
$ \sum_{g\geq 1} g(g^2(1-x)^2+3g(1-x^2)+2(1+4x+x^2)) x^g/6(1-x)^3, $
whereas the other factors are enumerated straight-forwardly as $v\SP^2Z(u/x)^g$.
Summation then gives the required result.

\subsubsection{Indent on bottom}

\begin{figure}
	\begin{center}
	\subfigure[The top factor extends farthest to the right.]{
		\includegraphics[scale=0.5]{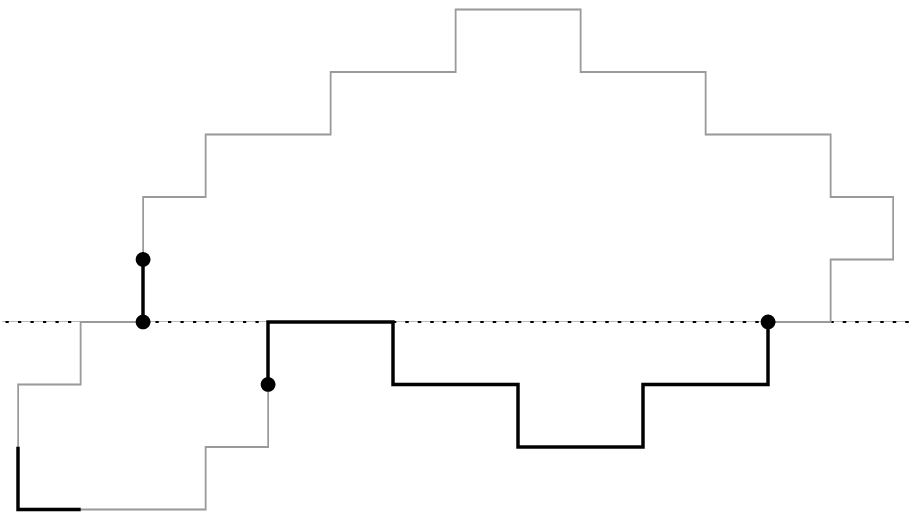}} \qquad
	\subfigure[The indent factor extends farthest to the right.]{
		\includegraphics[scale=0.5]{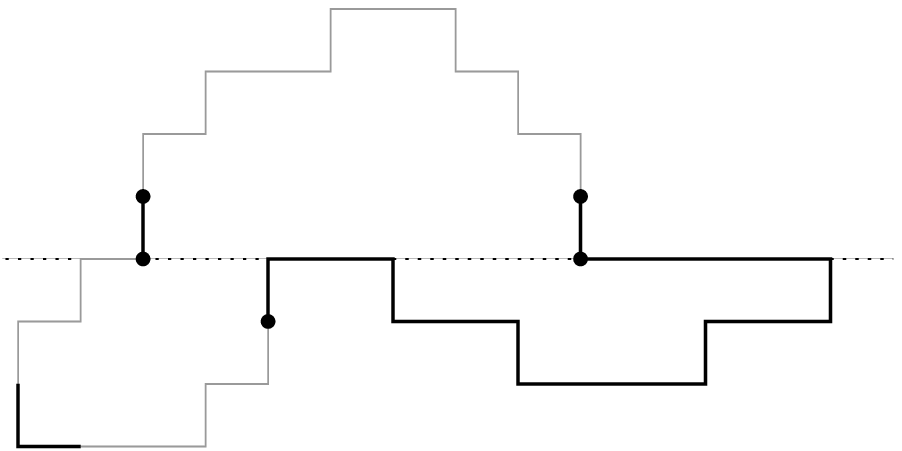}}
	\end{center}
	\caption{The form of a bimodal 2-unimodal polygon with the indent on the
			 bottom side.}
	\label{Fig:2-u.bim.bottom}
\end{figure}

If the indent of a bimodal 2-unimodal polygon is on the bottom, its form is one
of the two shown in Figure~\ref{Fig:2-u.bim.bottom}. The case where the top
factor extends farthest to the right is given by the inclusion-exclusion method
as
\begin{multline*}
	\Ex y{v^*x^2y^2(1-y)^4}{((1-y)^2-x)(1-x-y)^4}
	{\lr{\frac{u^*}{1-u^*-y}- \frac{x(1-u^*)}{(1-x)(1-x-y)} }} \\
	- v \SP \, \T_2 \lr{Z-\frac{1}{1-x}}\frac{2xy}{\Delta}.
\end{multline*}

If the indent is farthest to the right, we note that it is slightly different to
the case when the indent is in the corner, as we have already counted the case
where it is a pyramid. The generating function is
\begin{multline*}
	\sum_{g,h\geq 1} \lr{u^hv - \frac{x^hy}{1-x}}
	\frac{x(g^2(1-x)^2+2(1+x)+g(5-6x+x^2))}{2(1-x)^3} \\ \cdot
	\Ex y{xy^2}{1-x-y}{\lr{\frac{x}{1-y}}^g\lr{\frac{1}{1-y}}^h}.
\end{multline*}

%
%%
%%% Case 1
%%
%

\subsection{Case 1: indents in the same direction, on the same edge}
\label{ss_case1}

\subsubsection{Indents on top}

\begin{figure}
	\begin{center}
	\subfigure[The left-side indent is above the corner indent.]{
		\includegraphics[scale=0.45]{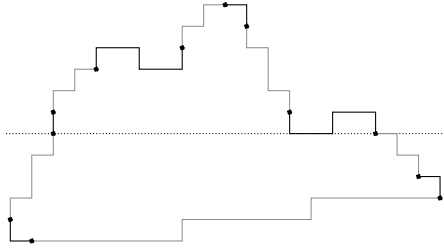}} \qquad
	\subfigure[Both indents are in the corner.]{
		\includegraphics[scale=0.45]{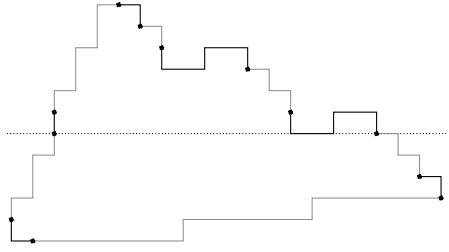}} \qquad
	\subfigure[The indents can be put at the same height. Then, the left indent
		can be moved into a position determined analytically by taking the
		derivative.]{
		\includegraphics[scale=0.45]{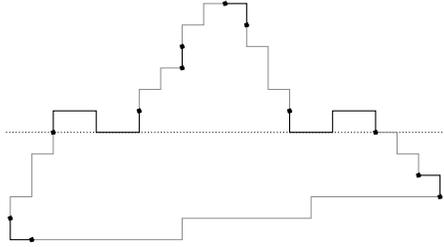}} \qquad
	\subfigure[The indent factor extends farthest to the right.]{
		\includegraphics[scale=0.45]{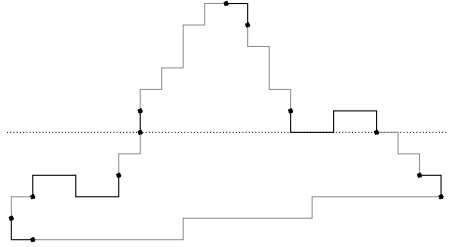}}
	\end{center}
	\caption{The form of a 2-unimodal polygon with its indents on the top.}
	\label{Fig:2-u.same}
\end{figure}

We enumerated the case where $m$-unimodal polygons have all their vertical
indents on the left side in Section~\ref{s_distinguished}. The case where the
indents are at the same height follow from Section~\ref{s_1-convex} \mm. We
therefore have two cases remaining to enumerate: when one indent is in the
corner and the other is either on the left or in the corner as well. Their
form is shown in Figure~\ref{Fig:2-u.same}.
We combine parts (a), (c) and (d) of the figure and enumerate it by
distinguishing a height of a 1-unimodal polygon with its indent in the corner,
giving the generating function
\[
	\T y^2 \dd y \T \SP Z \lr{u Z + \frac{1}{1-x}}.
\]
The case where both indents are in the corner can be generated  by taking the
derivative with respect to height of only the pyramid factor, or simply by
inspection, due its corner-staircase form. Its generating function is
\[
	\frac{u^4 v \SP Z^2}{(1-u)^3} \lr{\SP Z+\frac{v}{1-x}}.
\]

\subsubsection{Indents on bottom}

\begin{figure}
	\begin{center}
		\includegraphics[scale=0.55]{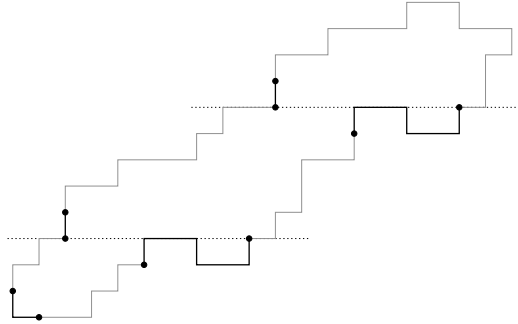}
	\end{center}
	\caption{The form of a 2-unimodal polygon with both indents on the bottom.}
	\label{Fig:2-u.both.bottom}
\end{figure}

If both indents are on the bottom, then the polygon is of the form illustrated in
Figure~\ref{Fig:2-u.both.bottom}. The generating function follows from
Sections~\ref{s_methodology} and \ref{ss_unimodal_bottom} \mm.

%
%%
%%% Case 2
%%
%

\subsection{Case 2: indents in the same direction, on different edges}
\label{ss_2-unimodal_opp}

\begin{figure}
	\begin{center}
	\begin{tabular}{ccc}
		\subfigure[The bottom indent is below the top indent.]
		{\includegraphics[scale=0.45]{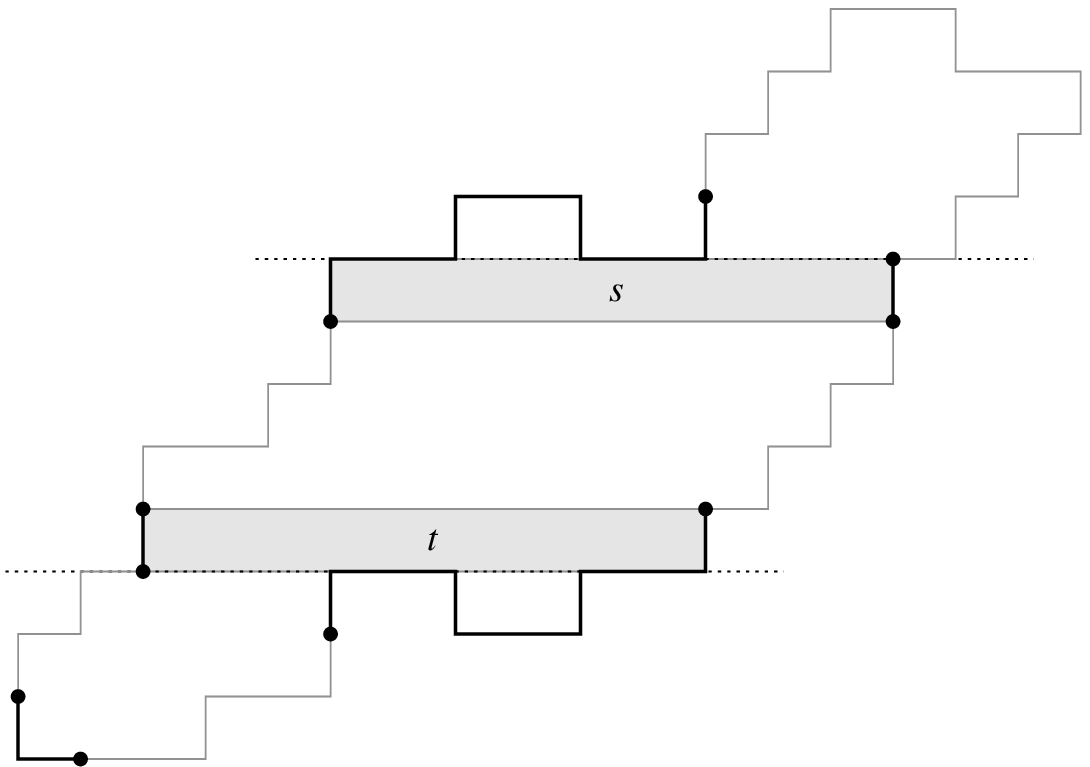}} &\qquad&
		\subfigure[The middle factor is further to the right than the top factor.]
		{\raisebox{27pt}
		 {\includegraphics[scale=0.45]{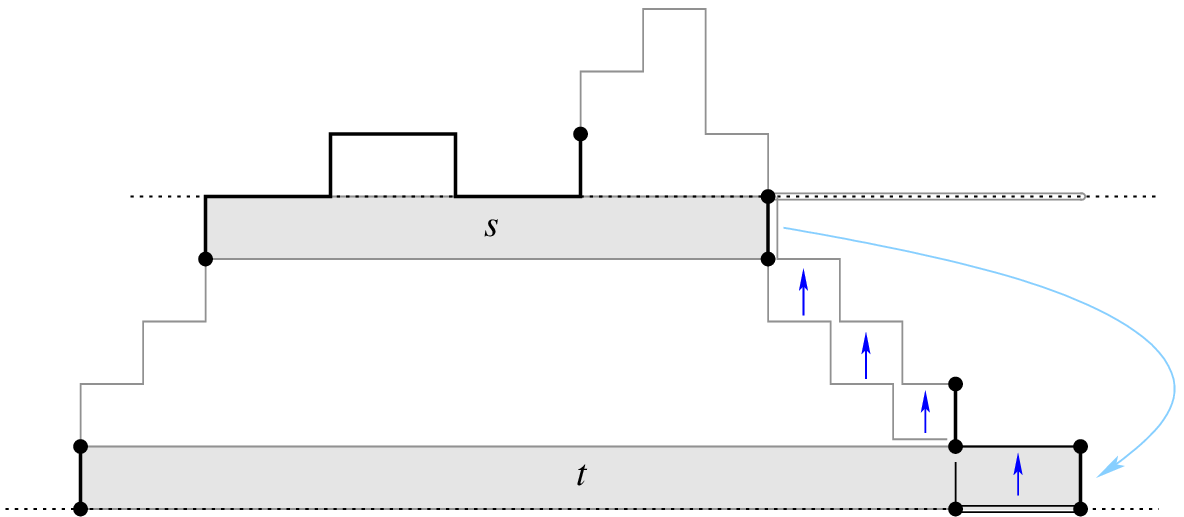}}}\\
		\subfigure[The indent factor is farthest to the right.]
		{\includegraphics[scale=0.45]{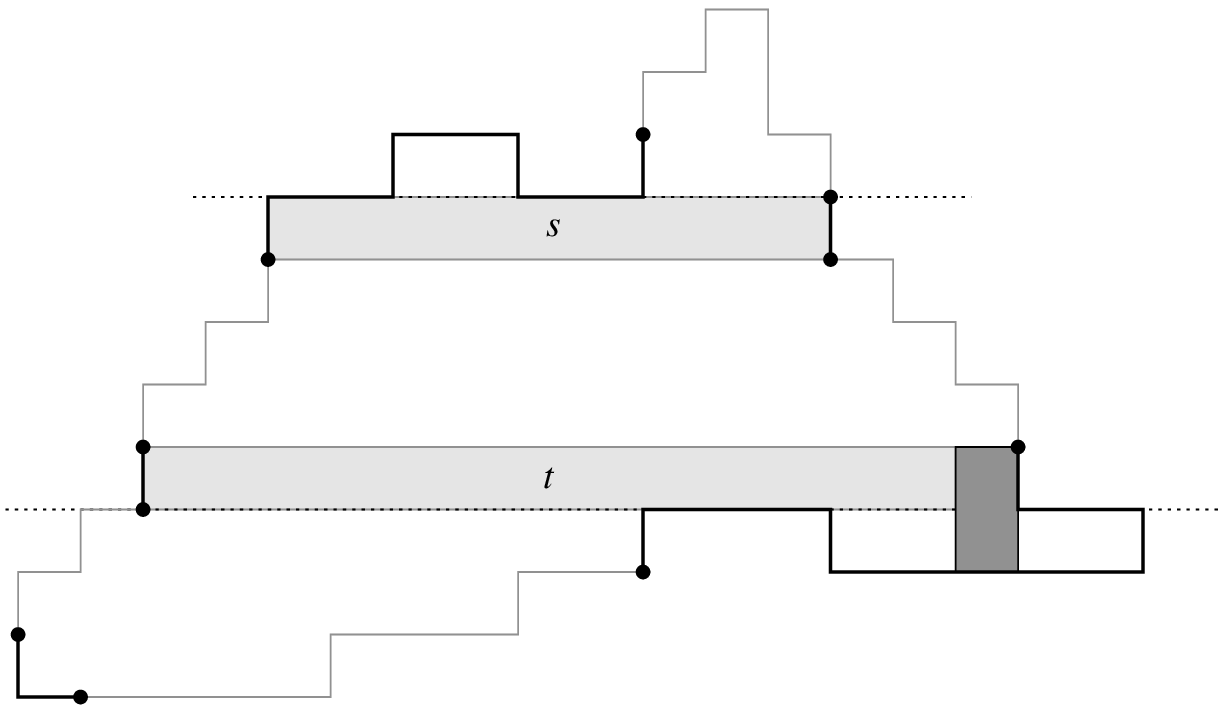}}&&
		\subfigure[The bottom indent is level with the top indent.]
		{\raisebox{12pt}
		{\includegraphics[scale=0.5]{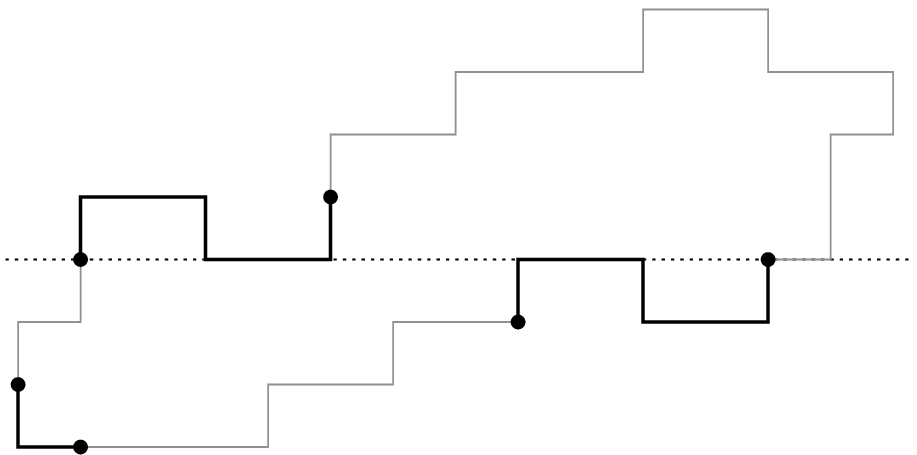}}}\\
		\subfigure[The level indents are interweaved.]
		{\qquad\includegraphics[scale=0.5]{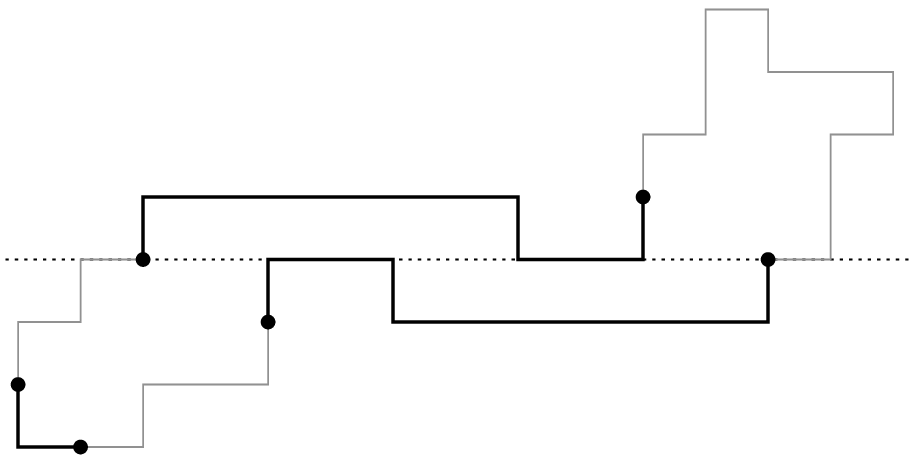}}&&
		\subfigure[Interweaved indents have a pyramid top factor.]
		{\includegraphics[scale=0.45]{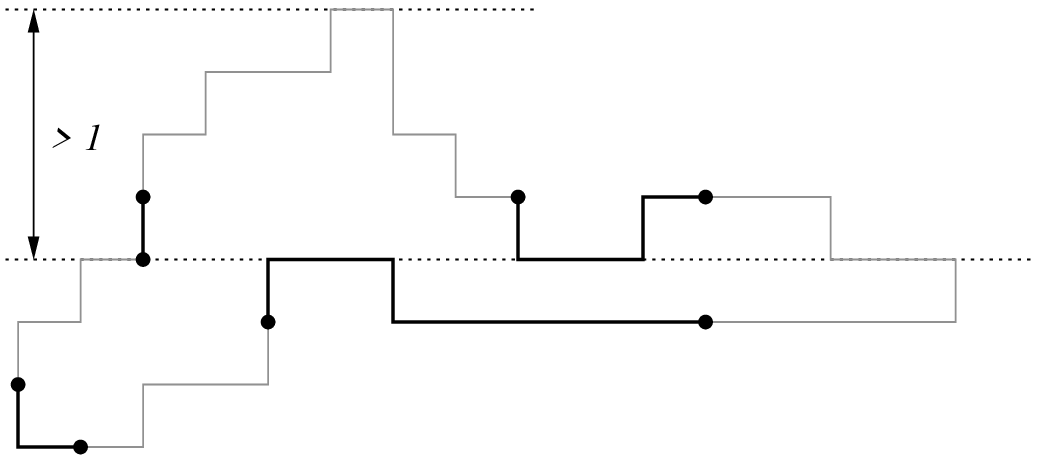}}\\
		\subfigure[The bottom indent is above the top indent.]
		{\includegraphics[scale=0.45]{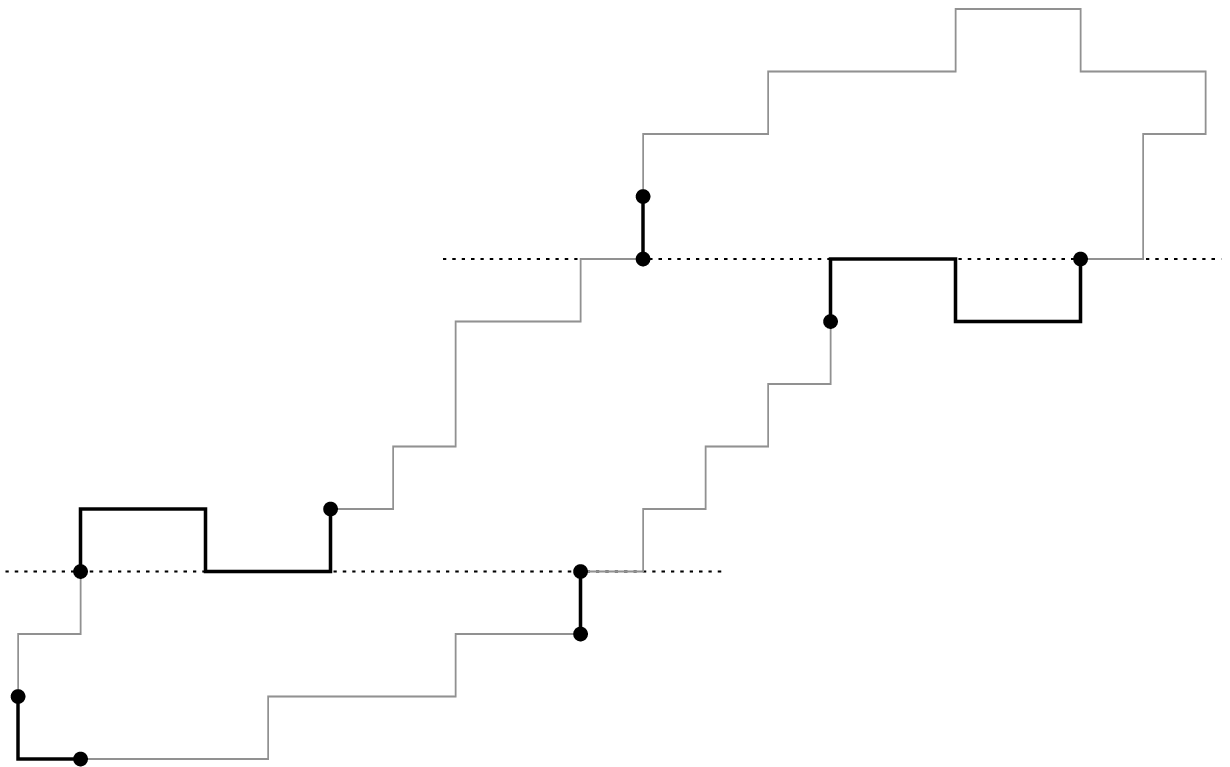}}&&
		\subfigure[The top indent is in the corner.]
		{\includegraphics[scale=0.45]{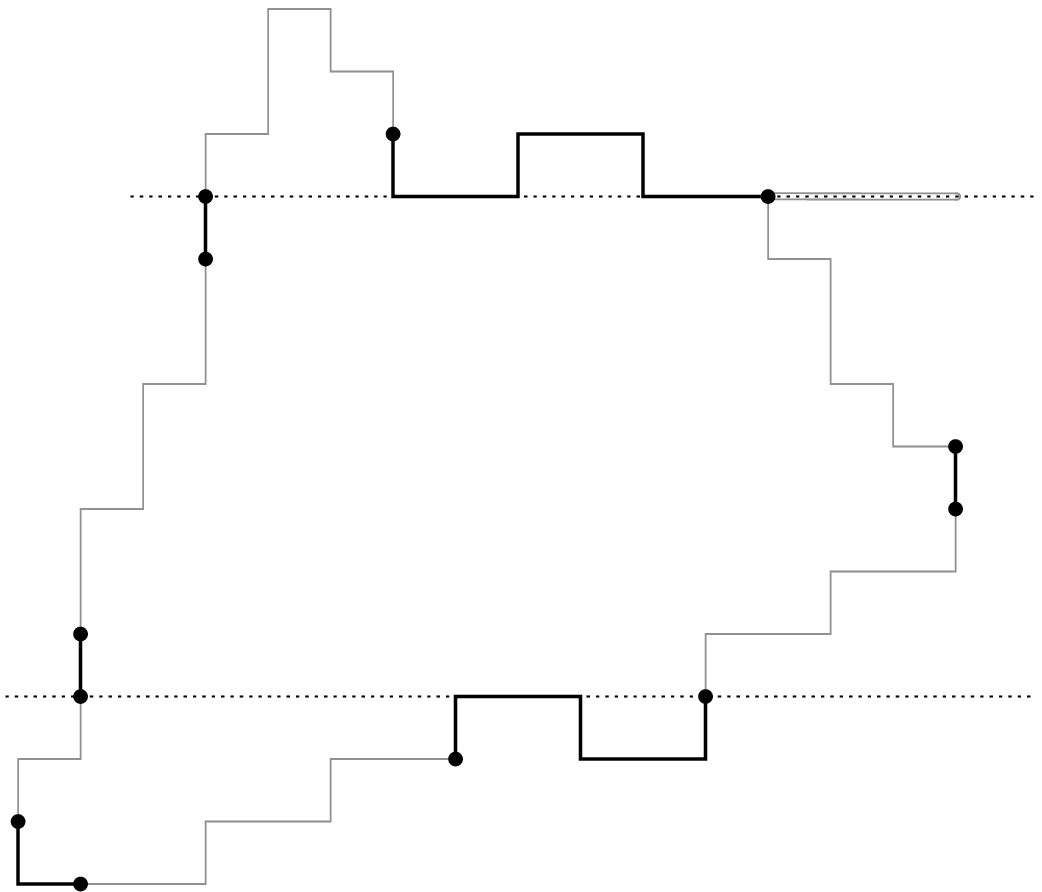}}
	\end{tabular}
	\end{center}
	\caption{The form of a 2-unimodal polygon with indents on opposite edges.}
	\label{Fig:2-u.opp}
\end{figure}

When the indents are both vertical and on different edges of the polygon,
the top indent is either on the left edge, or in the corner. They therefore
belong to one of the cases depicted in Figure~\ref{Fig:2-u.opp}.

\subsubsection{Top indent on the left}

Let us first consider the case where the top indent is on the left edge. This
case is given in Figure~\ref{Fig:2-u.opp}(a), which one can construct by joining
the middle factor, enumerated by $\bar{T}(s,t)$, with the top one. It is clear,
from Figure~\ref{Fig:2-u.opp}(b), that when the middle factor is farthest to the
right, such that fixed, joined steps of the top factor form a 1-dimensional loop
and the middle factor is wrapped, we may well form a 1-dimensional loop with the
steps along its base. Now, as was shown in Section~\ref{s_wrapping}, moving
the last vertical step down to the right-most edge preserves the self-avoiding
condition. However, we will now need to enumerate the base according to a
certain parameter, say $t$, and when the bottom horizontal steps change direction
and form a 1-dimensional loop along the bottom, the length of the base is
double-counted. We must therefore enumerate these configurations according to
the length of this 1-dimensional loop, which becomes the height one row shown in
the figure, in order to adjust for this. This can be done by taking the
derivative of the pyramid generating function to distinguish a height at which
to place the indent. We then remove the right-most height one row that has
weight $xt$ for each cell, and replace it by a row of at least length one that
weights its cells by $xrt$, such that $r$ counts the right-end length.
The generating function is therefore
\begin{equation}
	\label{Eq:2-u.opp.bottom.adj}
	\frac{xrt(1-xt)}{1-xrt} \lr{\frac{xty}{1-xt}}^2 \dd y \frac{P(xt,y)}{y}
\end{equation}
We therefore generate the configurations shown in Figure~\ref{Fig:2-u.opp}(a),
using the inclusion-exclusion principle, with the function
\begin{multline} \label{Eq:2-u.opp.below}
	\bigg(\lr{\EE{xy^2}{1-x-y}{\cdot\frac{s}{x-s}} - \frac{sv}{1-s-v} \cdot
	\frac{2xy}{\Delta} } \frac{s^2}{(1-s)^3} \odot_s \bar{T}(s,t) \\
-	\frac{x^4t^4y^2}{((1-xt)^2-y)^2} \lr{\frac{t}{1-xt^2} - \frac{1}{1-xt}} \bigg)
	\odot_t\ \frac{t^3v}{(1-t)^3(1-t-v)}.
\end{multline}

Now, if the indent factor is farthest to the right, it is of the form shown in
Figure~\ref{Fig:2-u.opp}(c). This form is very similar to the adjustment we had
to make in the previous enumeration and its generating function follows \mm.

The case in which the indents are side by side is shown in
Figure~\ref{Fig:2-u.opp}(d), (e) and (f). These are straight-forward to
enumerate using our standard techniques, as is part (g), which depicts the case
where the bottom indent is higher than the top indent.

\subsubsection{Top indent in the corner}

Let us now consider the case where the top indent is in the corner. Its
form is shown in Figure~\ref{Fig:2-u.opp}(h), where it can be seen to be similar
to the form of part (a) of the figure. Due to this
symmetry, we must ensure that the top, pyramid factor is of at least height two.
Its generating function then follows from the enumeration of parts (a) to (c) of
the figure \mm.

%
%%
%%% Case 3
%%
%

\subsection{Case 3: indents in different directions, on the same edge}

We now consider 2-unimodal polygons that have indents in different directions
on the same edge of the polygon. We note that this means that they are either on
the same side, or one is in the corner and the other is not. Now, the cases
where one indent is on the left are symmetrical, up to an interchange of
variables, with cases with one indent on the bottom. We therefore only consider
those cases where there is an indent on the left, or they are both in the
corner.

\begin{figure}
	\begin{center}
	\subfigure[The indents form a concave region on the left edge.]
		{\quad \includegraphics[scale=0.4]{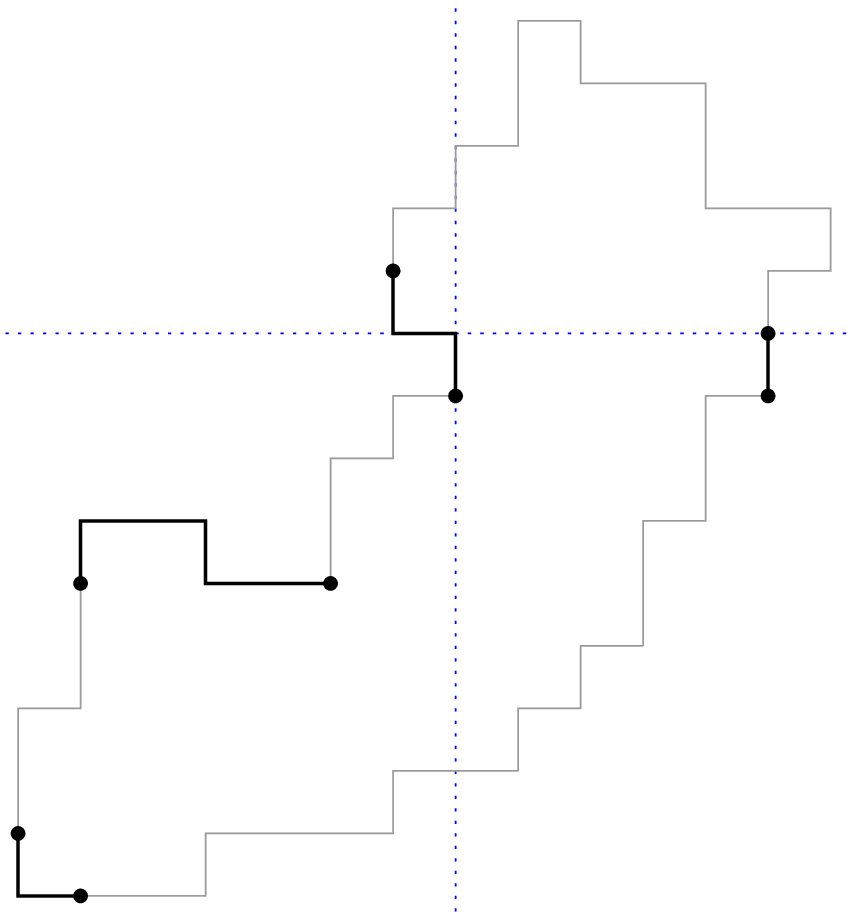}
		 \quad \qquad
		 \includegraphics[scale=0.5]{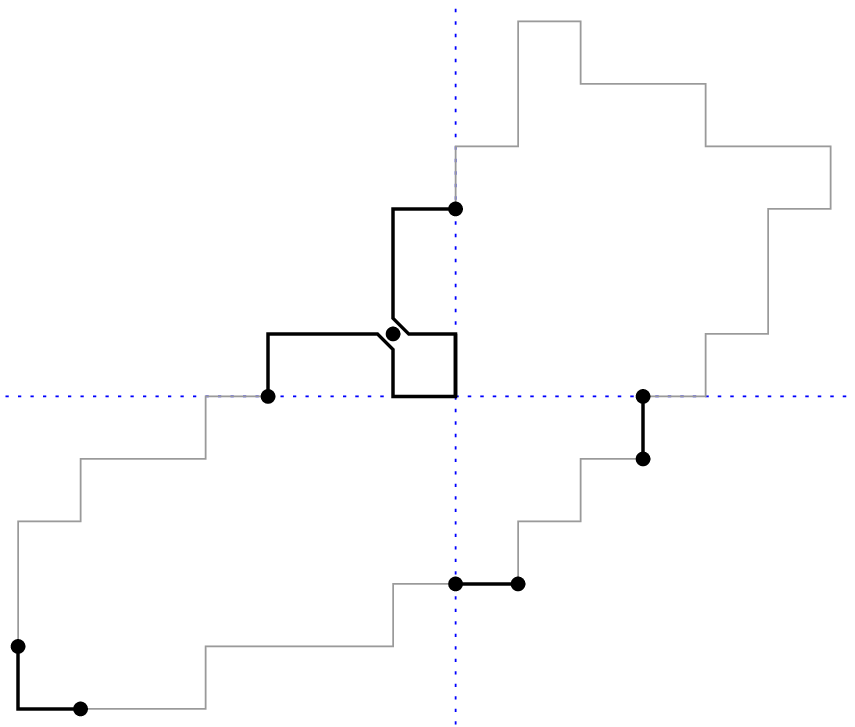}}
	\subfigure[The indents form a convex region on the left edge.]
		{\includegraphics[scale=0.4]{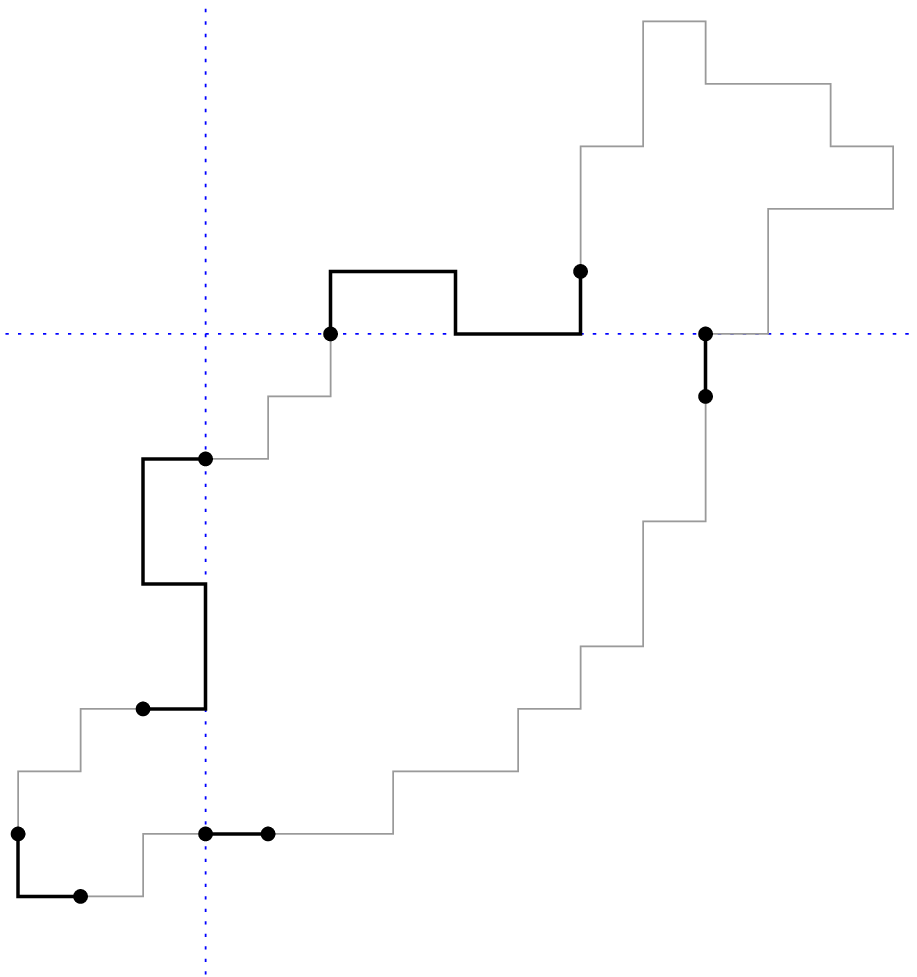}
		\qquad \qquad
		 \includegraphics[scale=0.5]{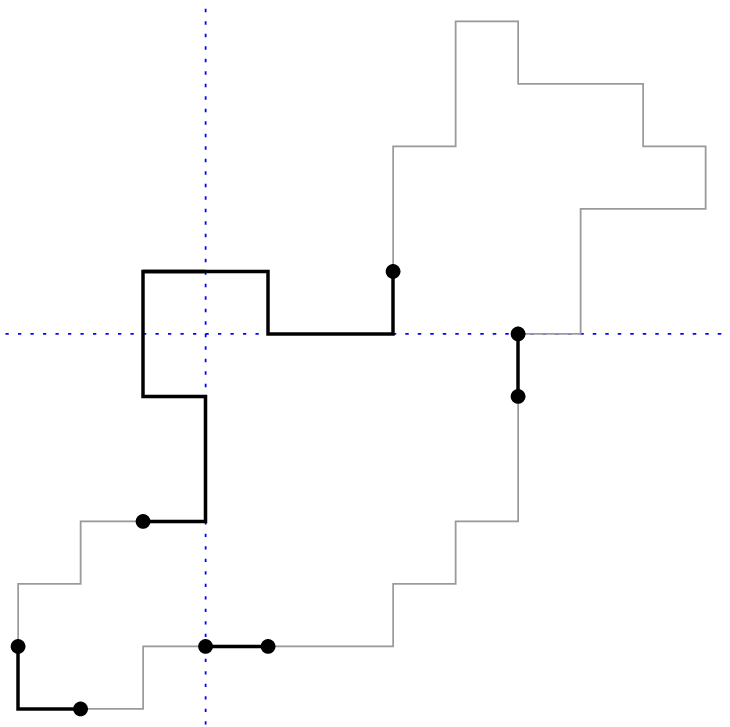}}
	\subfigure[The indents form a concave region in the corner.]
		{\includegraphics[scale=0.4]{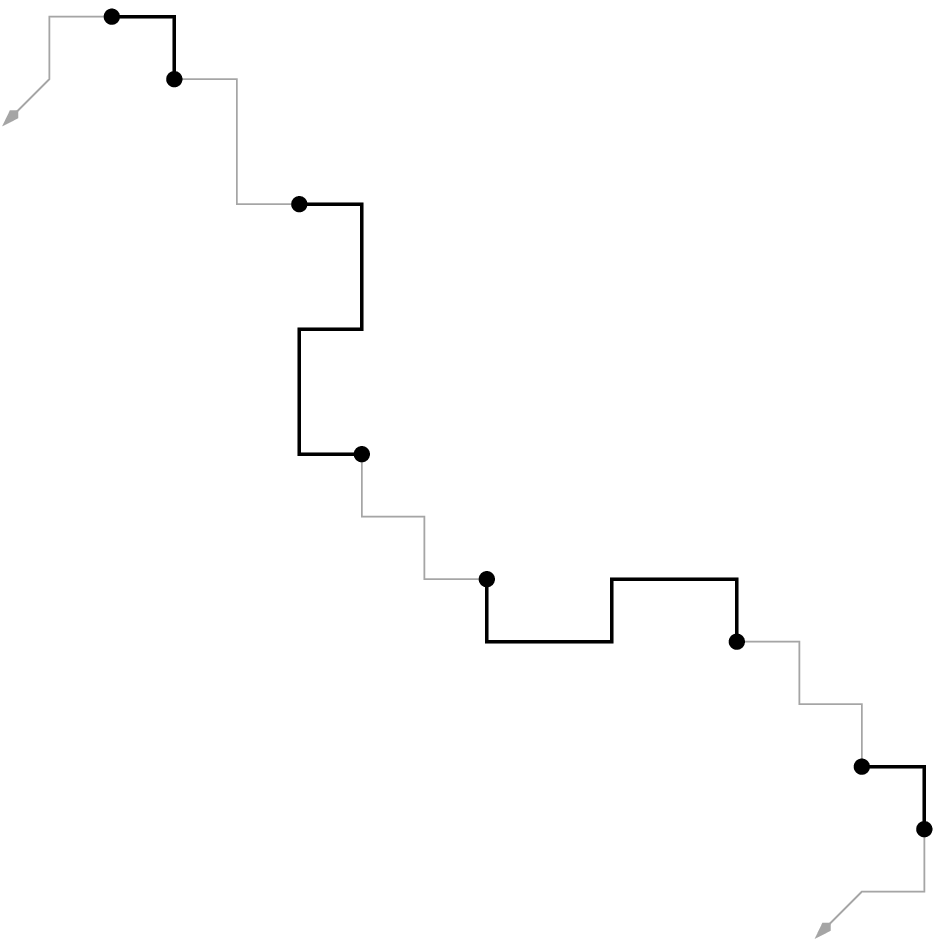} \qquad
		 \includegraphics[scale=0.5]{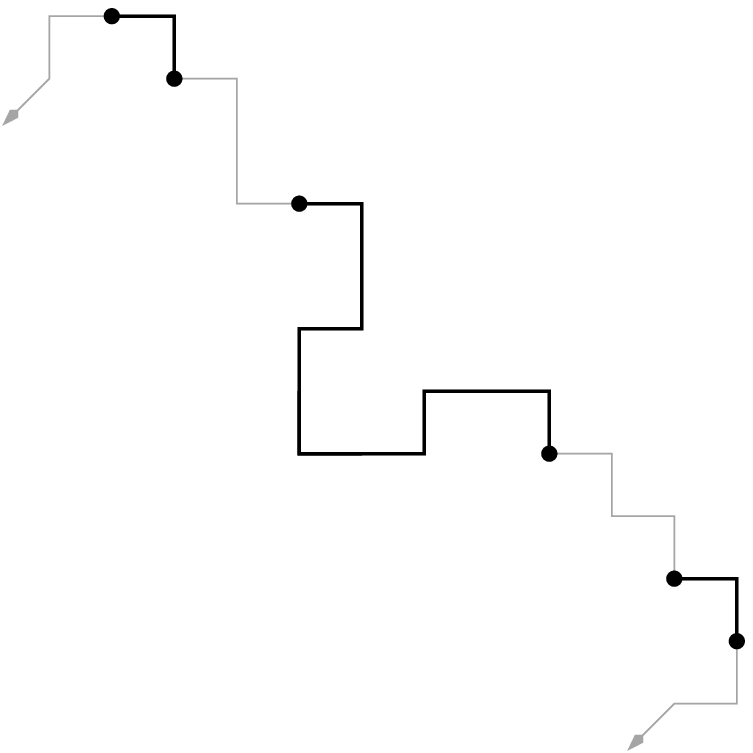}}
	\subfigure[The indents form a convex region in the corner.]
		{\includegraphics[scale=0.4]{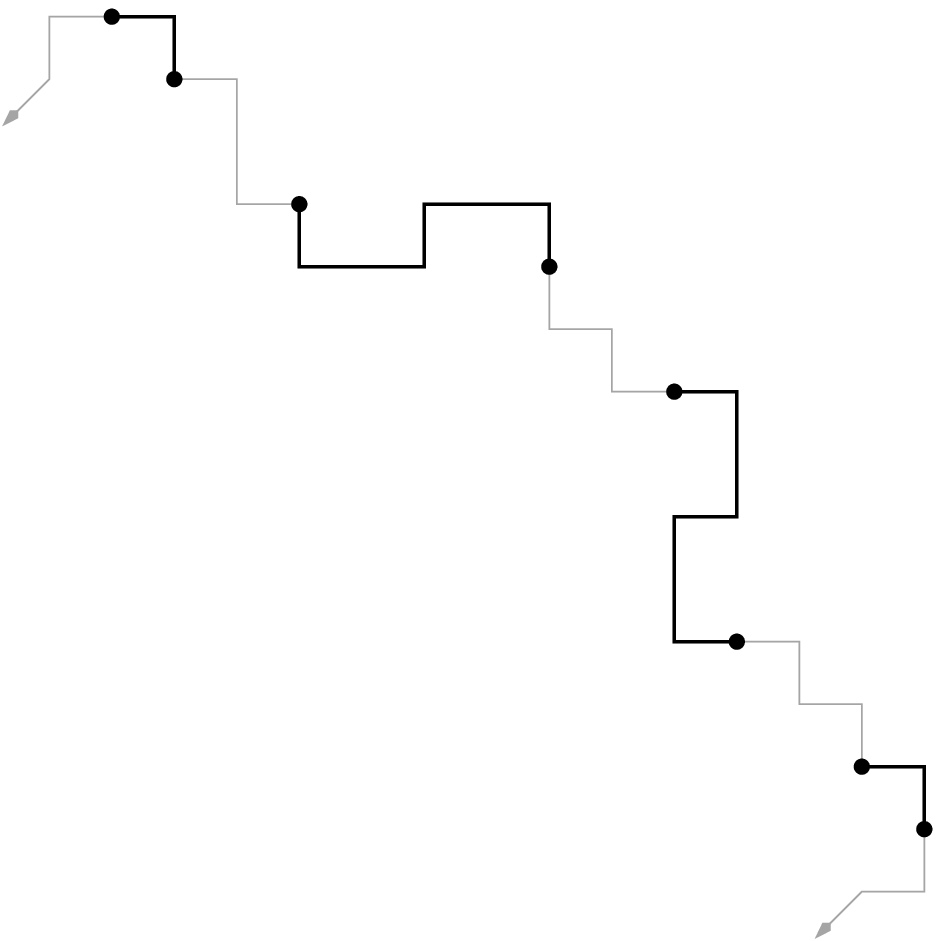} \qquad
		 \includegraphics[scale=0.5]{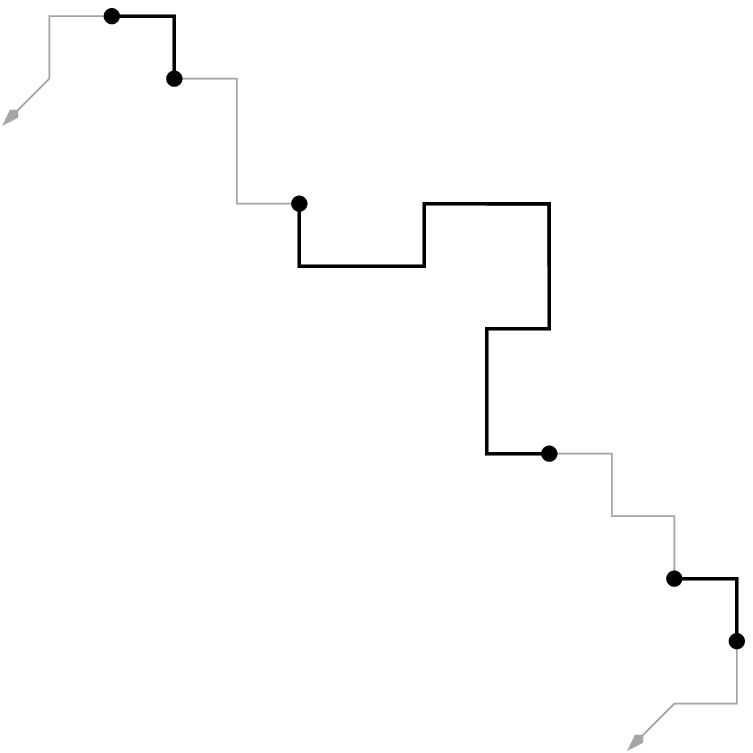}}
	\end{center}
	\caption{The form of a 2-unimodal polygon with indents in opposite directions.}
	\label{Fig:2-u.adj}
\end{figure}

We start by considering the two cases where they are on the same side.
There are two distinct cases: when both indents are on the left edge and when
they are in the corner. In either case, when the vertical indent is below the
horizontal one, the polygon is locally concave around the indent, in the
Euclidean sense. For example, see Figures~\ref{Fig:2-u.adj}(a) and (c).
Otherwise the polygon is locally convex around the indents, as in
Figure~\ref{Fig:2-u.adj}(b) and (d).

The remaining polygons have indents on adjacent sides. The only cases we have to
 consider therefore have one indent on the left edge and one in the corner. These
are: when the left indent is vertical and when it is horizontal.

We note that if a vertical (resp. horizontal) indent touches the top (resp.
right) edge of the minimum bounding rectangle, it could be considered to be on
the left (resp. bottom) edge \emph{or} in the corner.  We arbitrarily choose to
include such polygons in the case where the indents are on different sides.

\subsubsection{Indents on the left edge}

\paragraph{Locally concave.}
When the two indents form a concave edge, they are of the form shown in the
first diagram of Figure~\ref{Fig:2-u.adj}(a). We note that they must have a top
factor of width at least two, giving the generating function
\[
y \SP^2 Z \lr{ \EE{x^2 u^*}{(1-y)(1-x-y)(x-u^*)^2}{} - \frac{2x\SP}{\Delta^2} }.
\]
There is one sub-class of polygon whose indents intersect, whose form is shown
in the second diagram. These are enumerated by
\[
yu\SP \EE{x}{(1-y)(1-x-y)}{\lr{\frac{u^*}{(x-u^*)}-\frac{x}{1-x}}} -2\SP^4Z^3.
\]

\paragraph{Locally convex.}

When the two indents form a convex edge, they are of the form shown in
Figure~\ref{Fig:2-u.opp}(b).  Again, there is a subclass whose indents intersect
given in the second diagram.  We note that the vertical indent may not lie along
the top of the minimum bounding rectangle, and so the top factor must be of
height greater than one. The generating function is given by standard methods.

\subsubsection{Indents in the corner}

\paragraph{Locally concave.}

In this case, the horizontal indent is above the vertical one. And so, if the
indents are distinct, as in the left diagram of Figure~\ref{Fig:2-u.adj}(c),
the top factor is a pyramid with a horizontal indent. This is therefore
enumerated by $P'(u,y)$, given by equation \eqref{Eq:P'}, and the overall
generating function is
\begin{equation} \label{Eq:adj.corner.concave}
	\T \cdot P'(u,y) \cdot \lr{\SP Z + \frac{v}{1-x}},
\end{equation}
which can be re-expressed in the following form that is symmetric in $x$ and $y$
to reflect its geometrical symmetry:
\(
	{(1+u)(1+v)\SP^4 Z^3}/{(1-x)(1-y)}.
\)
Otherwise, there is a 2-dimensional indent, as in the right diagram in the
figure. Noting that the vertical (resp. horizontal) indent factor can extend
further to the right (resp.  top) than the bottom (resp. top) factor, the
generating function for this case is
\[
	\T \cdot \bar{\T} (1-(1-u)(1-v)) \lr{\SP Z + \frac{u}{1-y}}
	\lr{\SP Z + \frac{v}{1-x}}.
\]

\paragraph{Locally convex.}

There are no complications in this case. We again note that the indents may not
lie along the edge of the minimum bounding rectangle. Thus, the generating
function can be given on inspection of Figure~\ref{Fig:2-u.adj}(d), and
simplified to $(u+v)\SP^5 Z^3 /xy$.

%
%% Adjacent
%

\subsubsection{Horizontal indent on left, vertical indent in corner}

\begin{figure}
	\begin{center}
	\subfigure[The left indent is above the one in the corner.]
	{\includegraphics[scale=0.55]{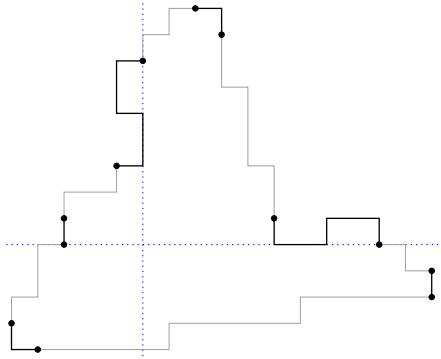}} \qquad
	\subfigure[The left indent is next to the vertical one.]{\raisebox{9pt}
	{\includegraphics[scale=0.6]{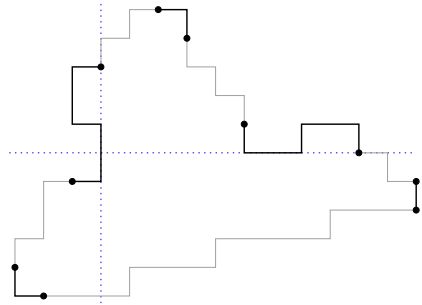}}}
	\subfigure[The indent factor is adjacent to the vertical indent.]
		{\includegraphics[scale=0.6]{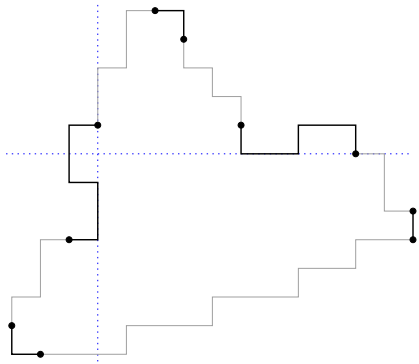}
		 \qquad \qquad
		 \includegraphics[scale=0.6]{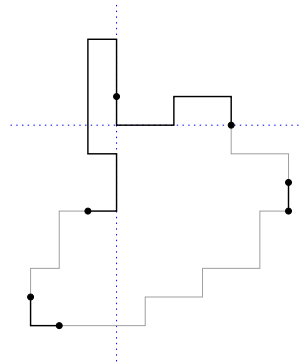}}
	\subfigure[The left indent is below the other.]
	{\quad\includegraphics[scale=0.6]{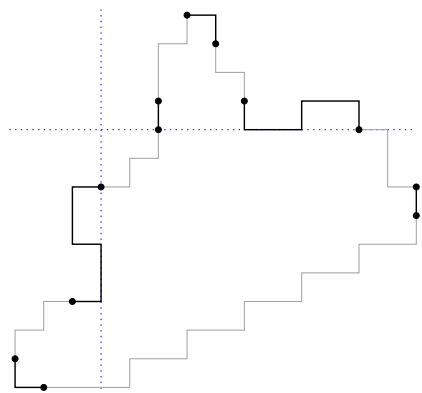}\quad}
	\end{center}
	\caption{The form of 2-unimodal polygons with a horizontal indent on the left,
			 and a vertical one in the corner.}
	\label{Fig:2-u.xyo}
\end{figure}

When there is a horizontal indent on the left edge, and a vertical one in the
corner, the polygons are of the form depicted in Figure~\ref{Fig:2-u.xyo}. They
are divided into four classes, according to the position of the horizontal
indent relative to the vertical indent: above, next-to, adjacent and below.
We say that the indents are {\em next to} each other if the {\em indentations}
formed overlap in height. If the {\em hump} of the indent factor (which is the
non-convex part of the interior of the polygon) is at the same height as the
vertical indent, we say that they are {\em adjacent}.

When the horizontal indent is above the one in the corner, we see that
reflecting the top, pyramid factor would give the case where both indents are
in the corner, such that the polygon is locally concave around the indents.
This implies that the indents are distinct, as in the first example of
Figure~\ref{Fig:2-u.adj}(c), which is enumerated by
\(
	\T\bar\T y\SP Z (1/(1-x)+uZ) (\SP/(1-y)+uZ).
\)

When the horizontal indent is next to or below the one in the corner, the
generating function is obtained simply by standard methods. This is also true
for when the horizontal indent factor is adjacent to the vertical indent, except
that we note that the pyramid factor in the first quadrant can be empty, as
shown in the second diagram of Figure~\ref{Fig:2-u.xyo}(c). In this case, the
generating function is given by
\begin{multline*}
	\lr{\frac s{1-s}}^2 \lr{\frac 1{1-x}+\frac{s^2}{x-s^2}}
	\frac{sy}{1-s-y} \lr{\frac y{1-y}+\frac{(1-s)^2}{(1-s)^2-y}} \\
	\ \odot_s\ \frac{xsty}{1-xs-ty}\ \odot_t\ \frac{ut^2}{(1-t)^2(1-u-t)}\\
+	\lr{\frac s{1-s}}^2 \lr{\frac 1{1-x}+\frac{s^2}{x-s^2}} \frac y{1-y} \\
	\ \odot_s\ \frac{xsty}{1-xs-ty}\ \odot_t\ \frac{ut^3}{(1-t)^2(1-u-t)}.
\end{multline*}

\subsubsection{Vertical indent on left, horizontal indent in corner}

\begin{figure}
	\begin{center}
	\subfigure[A horizontal indent on the bottom.]
		{\includegraphics[scale=0.6]{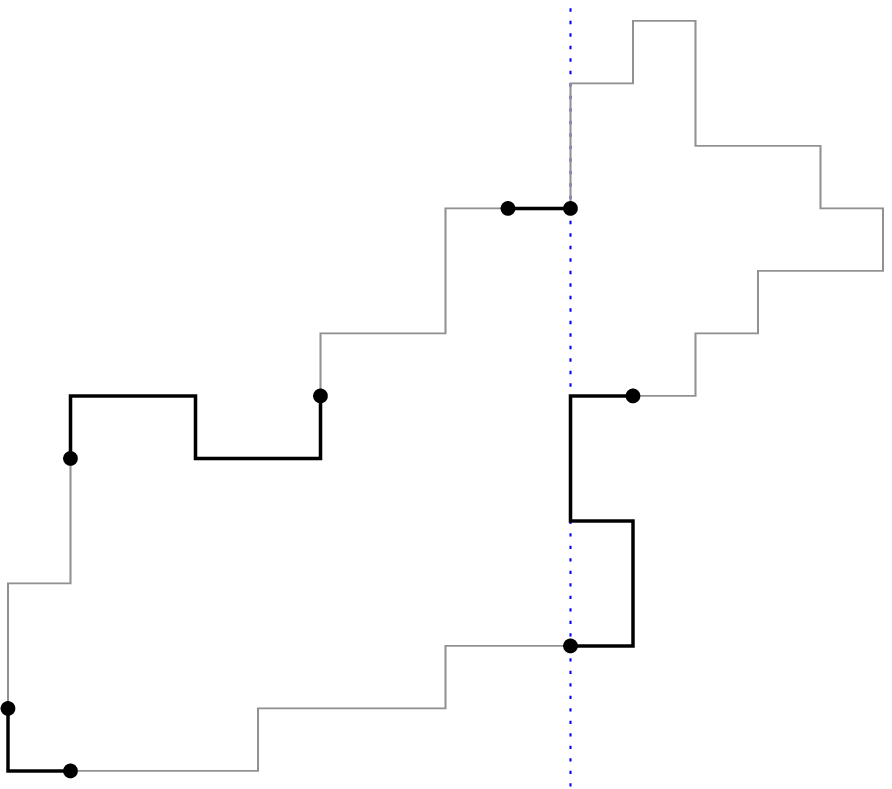} \qquad
		 \includegraphics[scale=0.6]{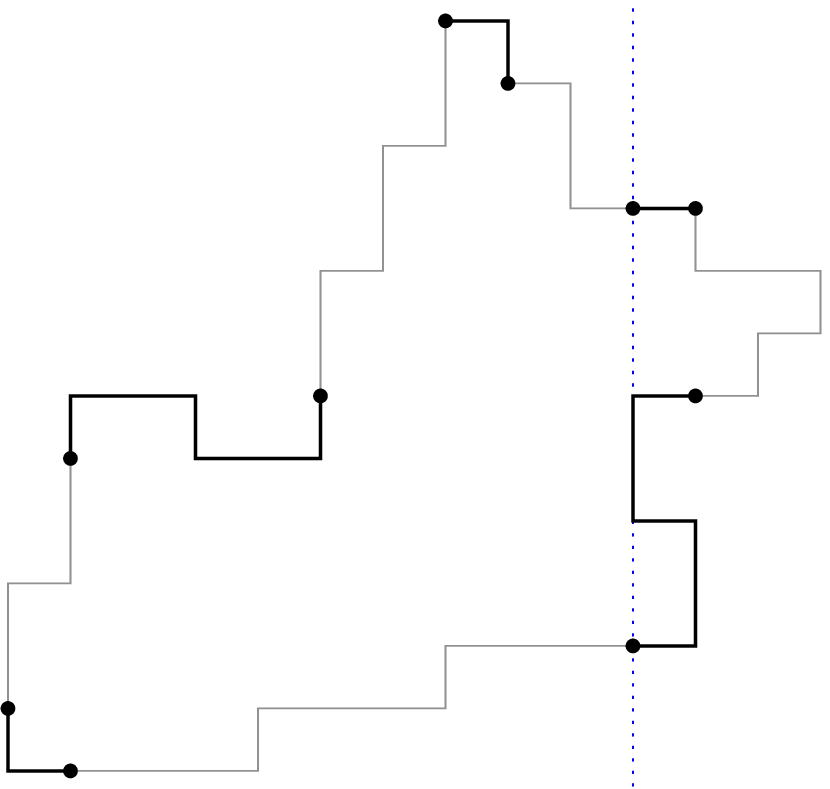}}
	\subfigure[The symmetrical class.]
		{\includegraphics[scale=0.6]{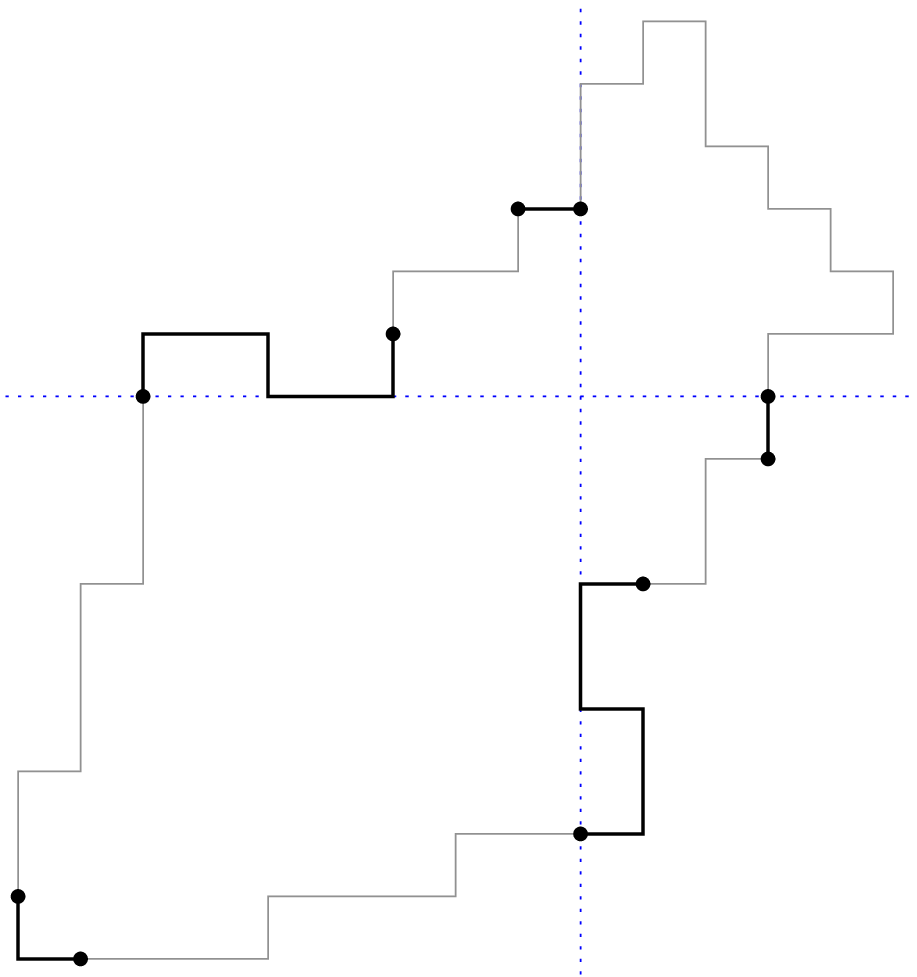}} \qquad
	\subfigure[The exclusion case of the symmetrical class.]
		{\includegraphics[scale=0.55]{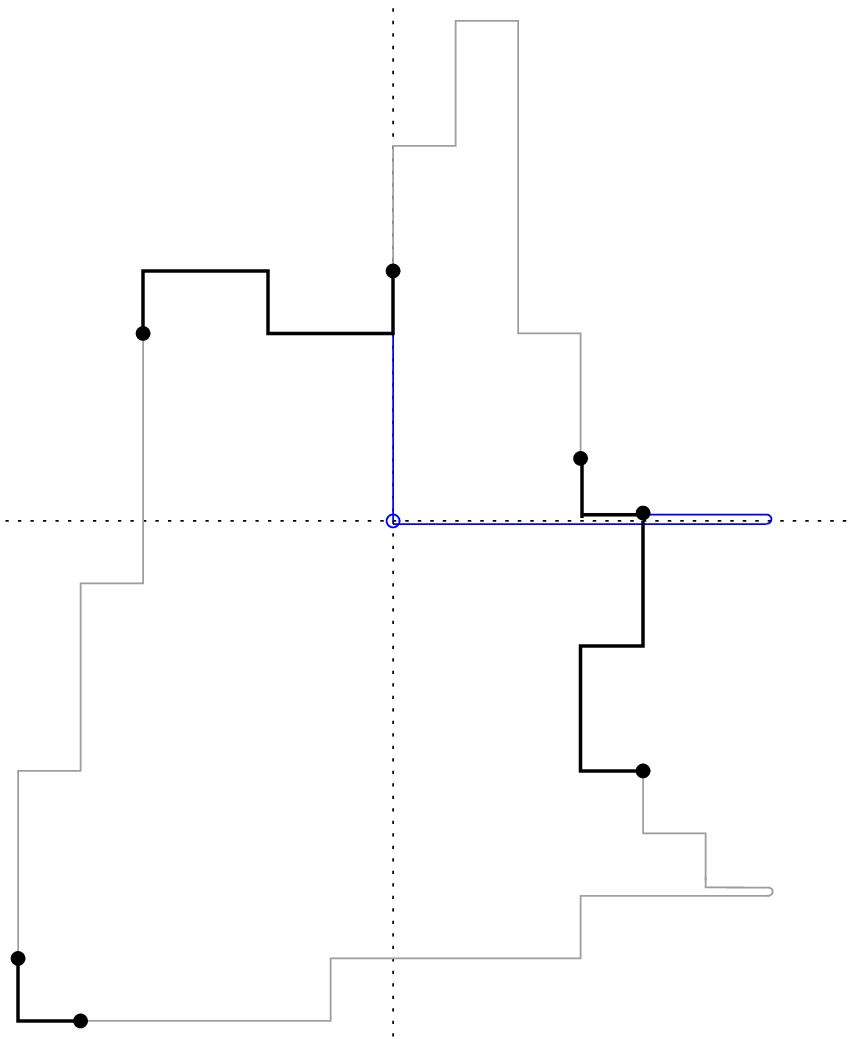}}
	\subfigure[A horizontal indent in the corner.]
		{\includegraphics[scale=0.6]{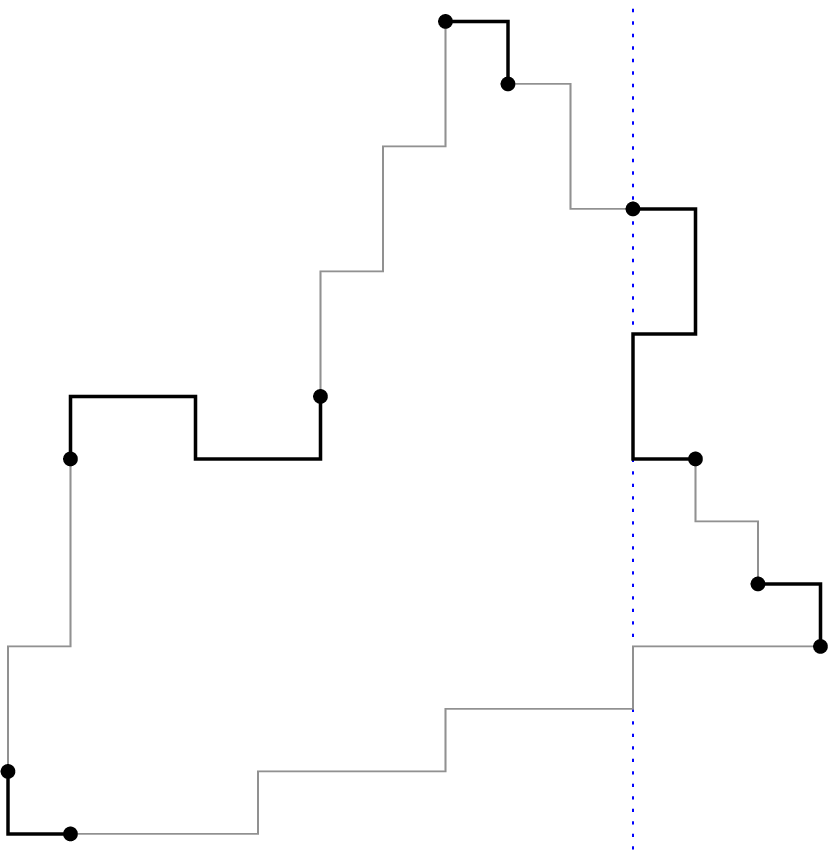}} \qquad
	\subfigure[The indent factor intersects itself.]
		{\qquad\quad
		\includegraphics[scale=0.58]{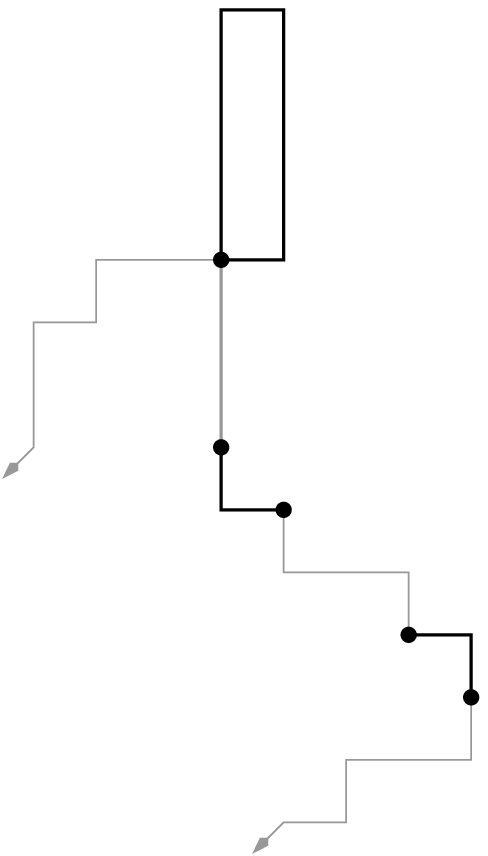}\qquad\quad}
	\end{center}
	\caption{The form of 2-unimodal polygons with a vertical indent on
			 the left, and a horizontal one on a different edge.}
	\label{Fig:2-u.against.y}
\end{figure}

The class of polygons that have a vertical
indent on the left edge and a horizontal one in the corner, as shown in
Figure~\ref{Fig:2-u.against.y}(d), is equivalent to the
second case of Figure~\ref{Fig:2-u.against.y}(a), after one flips the right
factor and the horizontal indent. However, there are three differences:
\begin{itemize}
\item we must allow the possibility of the indent touching the top of the
minimum bounding rectangle, such that it is at least as high as the left factor;
\item when the indent is the highest factor and intersects itself, as in
Figure~\ref{Fig:2-u.against.y}(e), a vertical indent can be placed next to the
width one indent factor (the ``hump'') in the top right, and the resulting
polygon will be self-avoiding; and
\item the right, pyramid factor must be of at least height two, for when the
indent lies along the right edge, we have chosen to consider this to be on the
bottom edge, rather than in the corner.
\end{itemize}
The first point adds a $1/(1-y)$ term to the expression, and the last adds a
multiplicative factor of $u/(1-v)$ to the contribution of the pyramid. The
second case is separate. Taking the derivate of the indent factor gives
$y^2 \dd y \lr{ \frac 1y \cdot \frac y{(1-y)} } = y^2/(1-y)^2$. This leaves us
with the following expression for the overall generating function:
\[
	\T \cdot y^2 \dd y \lr{
	\frac{\SP^3 Z}{xy} \cdot \frac{u}{1-v} \cdot \lr{\frac{1}{1-y}+vZ} }
+	\frac{u}{1-v} \lrf y{1-y}^2 v\, \T \SP Z.
\]

%
%%
%%% Case 4
%%
%

\subsection{Case 4: indents in different directions, on opposite edges}

When the indents of a 2-unimodal polygon are on opposite edges, there is one on
the left and one on the bottom. There are therefore two cases to consider:
when the left indent is vertical and when it is horizontal.

Before starting, we note that when we enumerate those cases where the indent
lies along the MBR there is one case that has already been enumerated in the
previous case: when the vertical indent is on the left perimeter, and the
horizontal indent lies on the right perimeter.

\subsubsection{Vertical indent on left, horizontal indent on bottom}

When we have a vertical indent on the left and a horizontal indent on the
bottom, an exchange of the variables $x$ and $y$ (corresponding to a reflection
in the $x=y$ axis) is an automorphism, giving the same class of polygons.  Now,
if the vertical indent is to the left of the horizontal one, as depicted in
Figure~\ref{Fig:2-u.against.y}(a), there is a staircase factor to the
bottom-left, allowing us to place the indent by taking the derivative. We can
obtain those with their horizontal indent below the vertical one via an exchange
of $x$ and $y$ in the generating function.
We can then obtain those polygons whose vertical indent is \emph{not} to the left
by excluding the double-counted, symmetrical case, shown in
Figure~\ref{Fig:2-u.against.y}(b), whose vertical indent is also to the left.
And so, we only have to enumerate these two classes.

Considering those polygons of the form given in Figure~\ref{Fig:2-u.against.y}(a),
we can place the indent distinguishing a vertical step in the left factor. We
must therefore separate the two cases shown (according
to whether the left or right factor is highest), as the usual bijection which
allows us to rejoin these cases does not distinguish between the contribution
to the height from the two factors. The first case, when the right factor is
highest, is enumerated by
\[
	\T y^2 \dd {y^*} \lr{ \frac{\bar{\T}^*}{y^*} \Ex x{u^*v^*}{1-x-v^*}
	{\cdot\frac{x^2(1-x)^2}{(1-x)^2-y}} - \frac{2xy}{\Delta}\cdot
	\frac{uv^*}{1-u-v^*}},
\]
where the derivative is taken to act only on those factors that are asterisked,
and afterward $y^*$ is equated with $y$. The case where the right indent is
higher is simpler, as there is no contribution to the height from the right
factor. We note that one can flip the right factor to obtain the form shown in
Figure~\ref{Fig:2-u.against.y}(d), whose generating function one can write on
inspection.
Finally, the symmetric case, shown in Figure~\ref{Fig:2-u.against.y}(b), is
enumerated using the inclusion-exclusion principle.
We then note that this enumerates polygons whose fixed steps along the base of
the top-right factor form a 1-dimensional loop, as shown in
Figure~\ref{Fig:2-u.against.y}(c), which we must exclude.  Similarly, the fixed
steps along the left edge may form a vertical 1-dimensional loop. 
Surprisingly, the generating function simpliies to $xy(1-x-y)/\Delta^{3/2}$, which is
symmetrical in $x$ and $y$. Thus, the exclusion case for both the horizontal and
vertical 1-dimensional loop have the same generating function. Finally, the
doubly-excluded case, of both a horizontal and vertical 1-dimensional loop, is
re-included to give the required result.

\subsubsection{Horizontal indent on left, vertical indent on bottom}

\begin{figure}
	\begin{center}
	\subfigure[The left indent is above the one on the bottom.]
		{\includegraphics[scale=0.6]{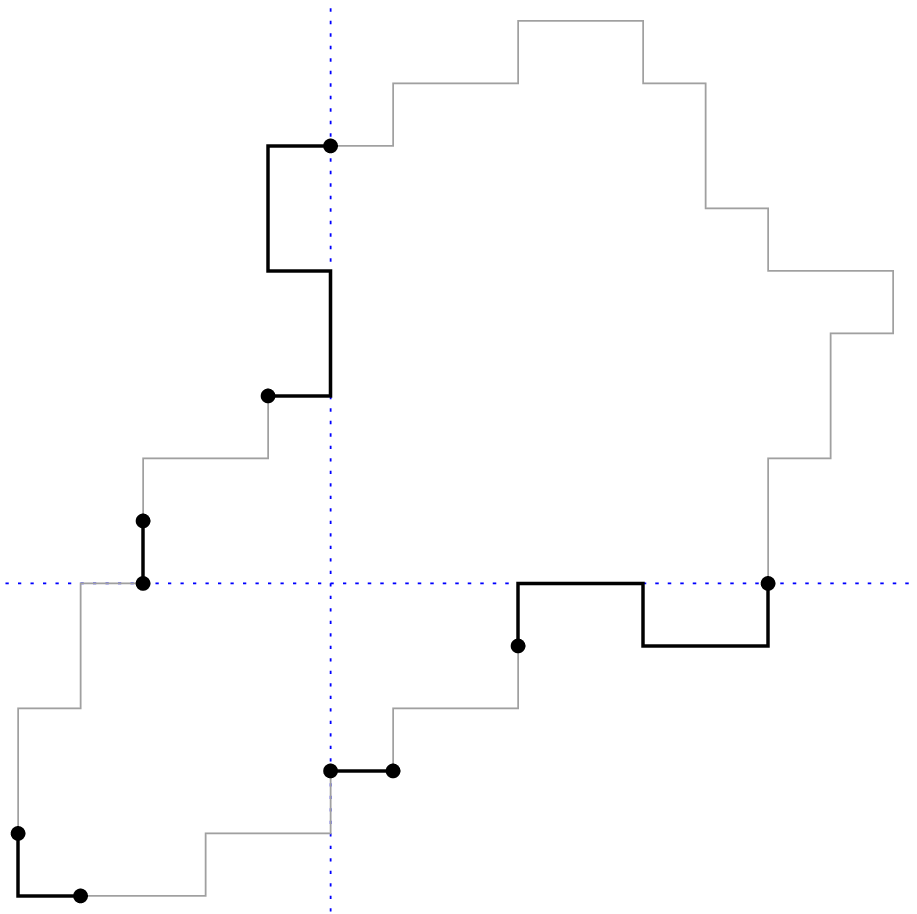}\qquad
		 \includegraphics[scale=0.6]{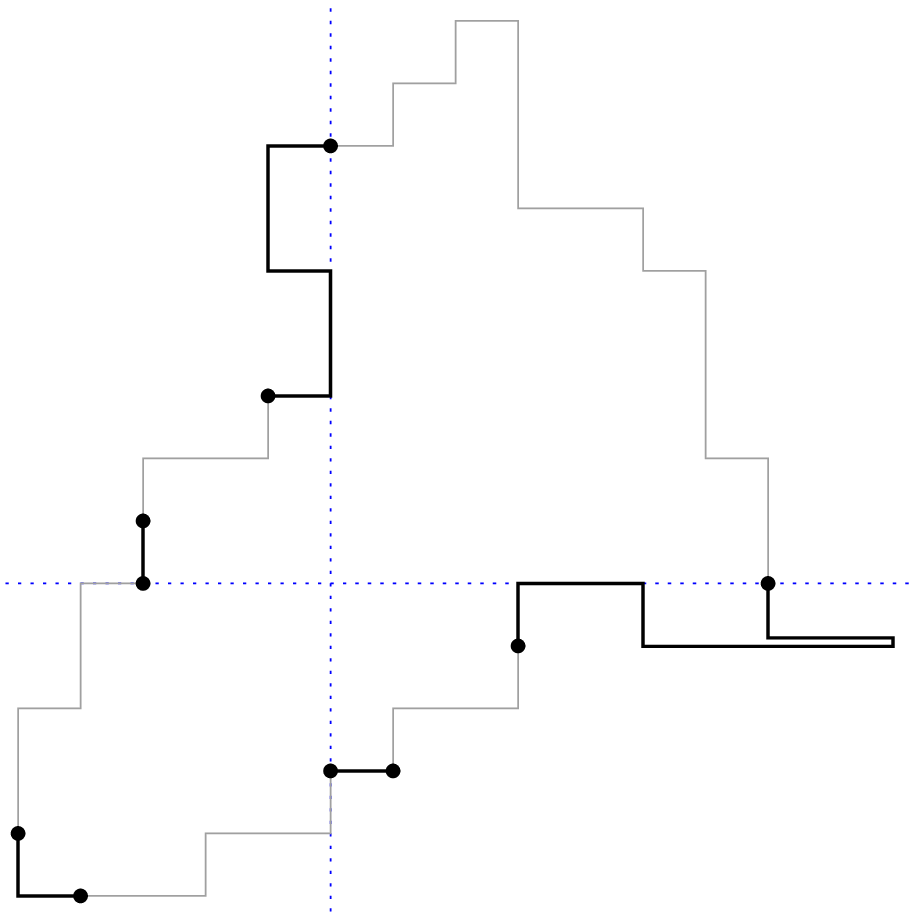}}
	\subfigure[The left indent is next to the vertical one.]
		{\includegraphics[scale=0.6]{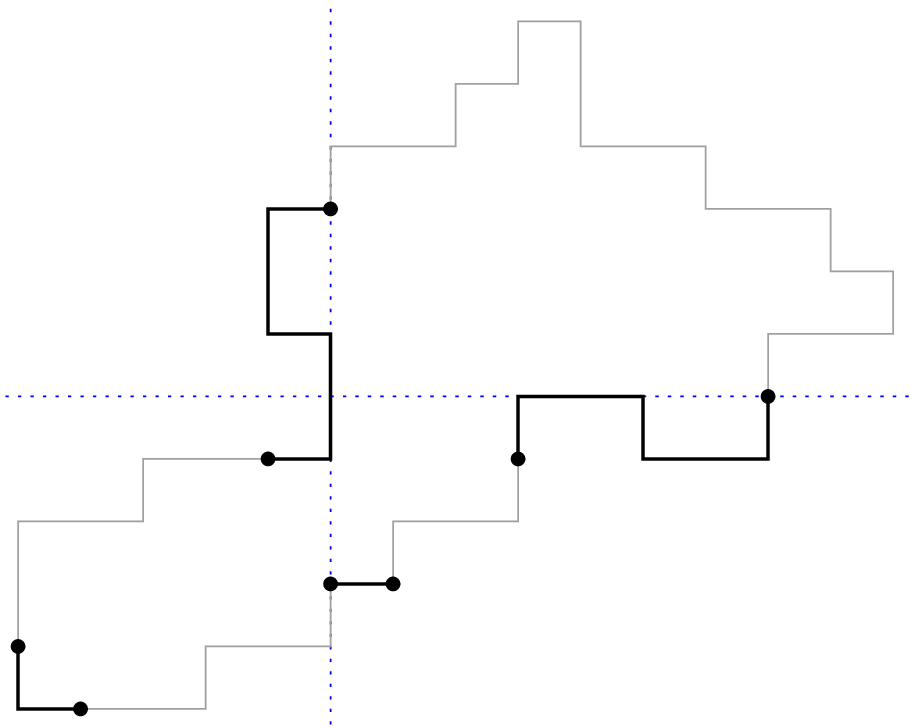}}\qquad
	\subfigure[The indent factor is adjacent to the vertical indent.]
		{\includegraphics[scale=0.6]{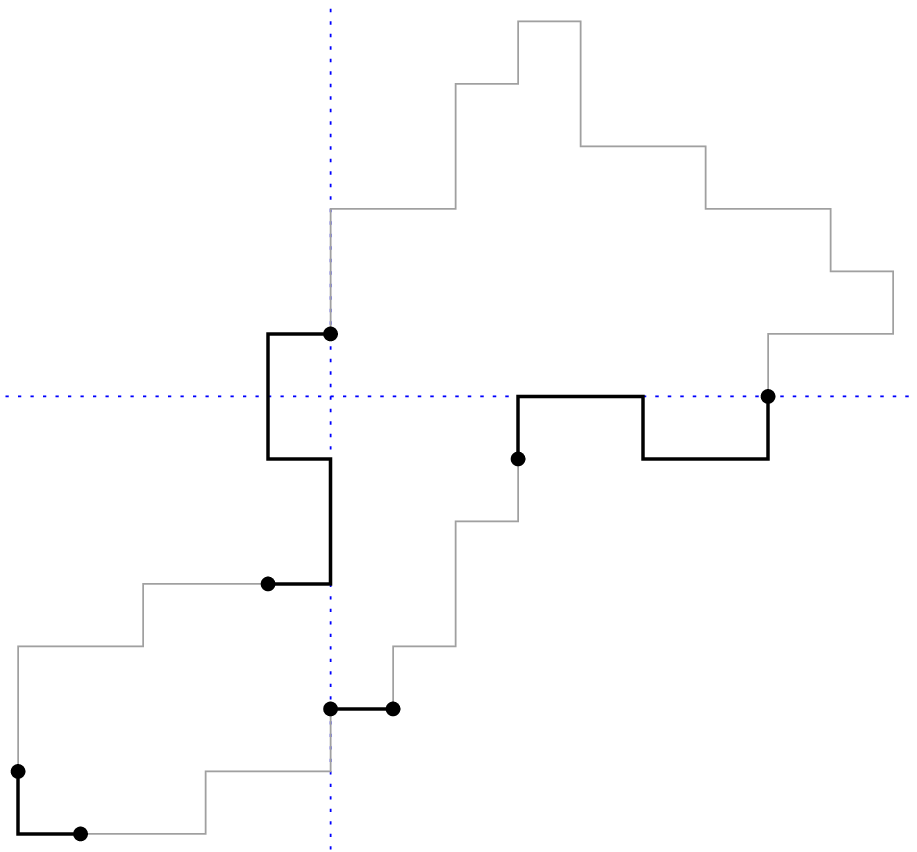}}
	\subfigure[The left indent is below the other.]
		{\includegraphics[scale=0.6]{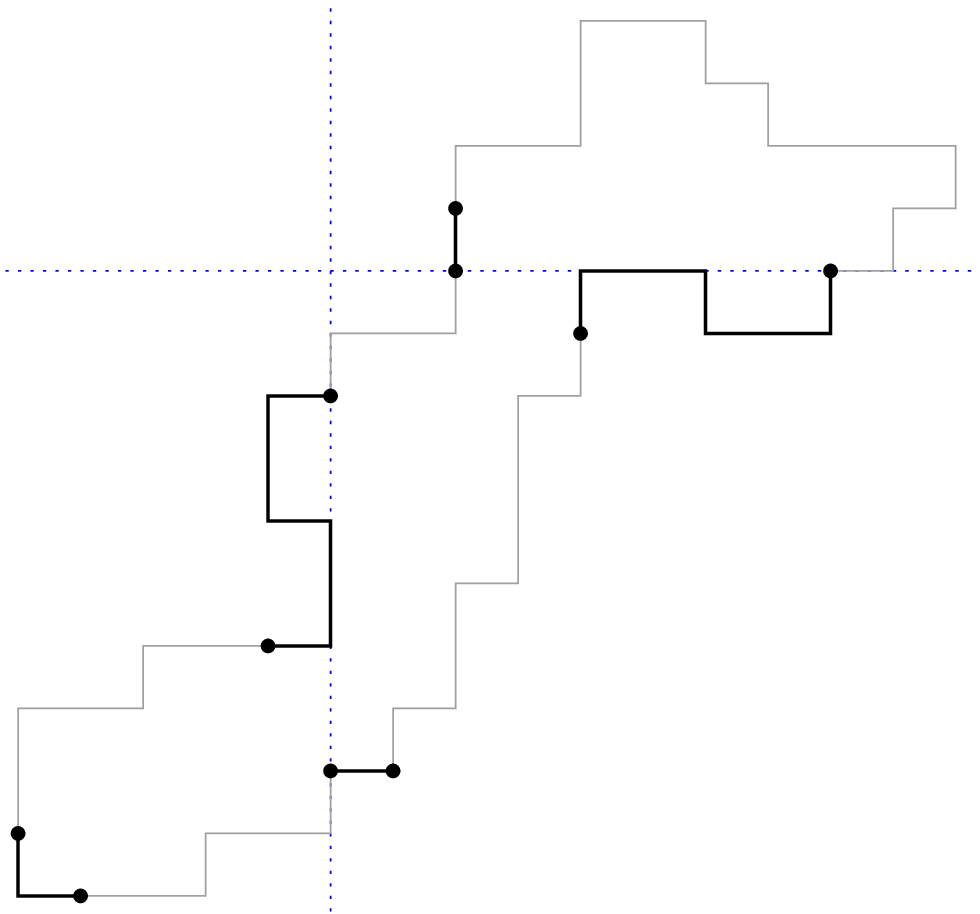}}
	\end{center}
	\caption{The form of 2-unimodal polygons with a horizontal indent on the left,
			 and a vertical one on the bottom.}
	\label{Fig:2-u.xoy}
\end{figure}

When the indents of a 2-unimodal polygon consist of a horizontal one on the left
edge and a vertical one on the bottom, the polygon is of a form depicted in
Figure~\ref{Fig:2-u.xoy}. Again, these forms divide the polygons into four
parts according to the height of the left indent, relative to that of the
vertical indent.

When the horizontal indent is above the vertical one, the latter indent could
be to the right, under, adjacent to, or to the left of the horizontal
indent above it. However, similarly to the previous case, the latter three
sub-cases are symmetric with other cases, and we need not enumerate them.

All of the following four cases can be enumerated by an inclusion-exclusion
argument, together with our other standard techniques.
We note that we have to exclude the polygons whose fixed steps along the
factorisation lines form a  1-dimensional loop, as in the second diagram of part
(a) of the figure. We also note that we must enumerate the cases
where the indent factors are farthest to the right or the top. 

\subsection{The 2-unimodal generating function}

Combining all of the above cases gives the following generating function for
2-unimodal polygons:
\begin{multline} \label{Eq:2-unimodal}
	\frac{A^{(u)}_2}{2x^3y^3(1-x)^3(1-y)^3\Delta^{5/2}} \\
	- \frac{B^{(u)}_2}{2x^3y^3(1-x)^5(1-y)^5((1-x)^2-y)^3((1-y)^2-x)^3(1-x-y)},
\end{multline}
{\small
\begin{minipage}[t]{\linewidth}$
\mbox{\normalsize{where }} A^{(u)}_2 = 
\ 2 (1-x)^{11} x^4
- 2 (1-x)^9 x^4 (13-10 x+3 x^2) y
+ 2 (1-x)^7 x^2 (4-10 x+77 x^2  -97 x^3+40 x^4-23 x^5+3 x^6) y^2
- 2 (1-x)^5 x^2 (38-109 x+328 x^2-501 x^3  +314 x^4-113 x^5-5 x^6-15 x^7+x^8) y^3
+ (1-x)^4 (2-18 x+378 x^2-878 x^3  +1656 x^4-2024 x^5+1145 x^6-294 x^7+45 x^8+122 x^9+2 x^{10}) y^4
- (1-x)^3 (22-188 x+1342 x^2-2610 x^3 +
  3402 x^4-3350 x^5+1832 x^6-337 x^7+192 x^8  -159 x^9-40 x^{10}+2 x^{11}) y^5
+ (1-x) (110-1012 x+4640 x^2-10118 x^3+12949 x^4  -11749 x^5+7989 x^6-3394 x^7+177 x^8+213 x^9+175 x^{10}-120 x^{11}+12 x^{12}) y^6 
+ (-330+2958 x-10938 x^2+21112 x^3-24520 x^4+19285 x^5-11383 x^6+5556 x^7  -1584 x^8-611 x^9+555 x^{10}-160 x^{11}+28 x^{12}) y^7
+ (660-5172 x+14664 x^2  -20918 x^3+17843 x^4-10049 x^5+3571 x^6-1584 x^7+1342 x^8-295 x^9+82 x^{10}  -32 x^{11}) y^8
+ (-924+6300 x-13980 x^2+14186 x^3-8426 x^4+3578 x^5+36 x^6  -611 x^7-295 x^8-24 x^9+14 x^{10}) y^9
+ (924-5460 x+9708 x^2-6454 x^3+2105 x^4  -1350 x^5-38 x^6+555 x^7+82 x^8+14 x^9) y^{10}
+ (-660+3372 x-5012 x^2+1968 x^3  +250 x^4+547 x^5-295 x^6-160 x^7-32 x^8) y^{11}
+ (330-1458 x+1930 x^2-496 x^3 -431 x^4-33 x^5+132 x^6+28 x^7) y^{12}
- 2 (1-x)^2 (55-101 x+7 x^2+35 x^3+6 x^4) y^{13} 
- 2 (1-x)^3 (-11+4 x+x^2) y^{14}
- 2 (1-x)^3 y^{15}
$\end{minipage}\\[0.5\baselineskip]

\noindent
\begin{minipage}[t]{\linewidth}$
\mbox{\normalsize{and }} B^{(u)}_2 = 
-2 (1-x)^{18} x^4
+ 2 (1-x)^{16} x^4 (20 - 26 x + 5 x^2) y
- 2 (1-x)^{14} x^2 (4 - 10 x + 189 x^2 - 464 x^3 + 388 x^4 - 122 x^5 + 10 x^6) y^2
+ 2 (1-x)^{12} x^2 (66 - 243 x + 1343 x^2 - 4004 x^3 + 5623 x^4 - 4025 x^5 +
   1459 x^6 - 230 x^7 + 10 x^8) y^3
- (1-x)^{10} (2 - 20 x + 1096 x^2 - 5130 x^3 + 17888 x^4 - 49004 x^5 + 86111 x^6
  - 92423 x^7 + 60313 x^8 - 23283 x^9 + 4864 x^{10} - 440 x^{11} + 10 x^{12}) y^4
+ (1-x)^8 (36 - 404 x + 6712 x^2 - 34332 x^3 + 107222 x^4 - 257436 x^5 + 471559 x^6 - 618724 x^7 + 559981 x^8 - 341505 x^9 + 135916 x^{10} - 33115 x^{11} + 4316 x^{12} - 220 x^{13} + 2 x^{14}) y^5
+ (1-x)^6 (-306 + 3806 x - 36028 x^2 + 177500 x^3 - 546201 x^4 + 1230135 x^5 -
2183847 x^6 + 3061206 x^7 - 3288667 x^8 + 2624163 x^9 - 1512420 x^{10} + 609513
x^{11} - 163267 x^{12} + 26485 x^{13} - 2144 x^{14} + 52 x^{15}) y^6
- (1-x)^4 (-1632 + 22272 x - 170304 x^2 + 786188 x^3 - 2405909 x^4 + 5350296 x^5
  - 9271113 x^6 + 13015854 x^7 - 14902848 x^8 + 13692922 x^9 - 9840757 x^{10} +
	5379124 x^{11} - 2166934 x^{12} + 616413 x^{13} - 115521 x^{14} + 12519
	x^{15} - 582 x^{16} + 4 x^{17}) y^7 +
$\end{minipage}
\noindent
\begin{minipage}[t]{\linewidth}$
(1-x)^3 (-6120 + 84760 x - 581646 x^2
	+ 2434626 x^3 - 6888332 x^4 + 14198657 x^5 - 22614758 x^6 + 
    29071458 x^7 - 30896121 x^8 + 27177494 x^9 - 19436717 x^{10} + 10987113 x^{11} - 4750238 x^{12} + 1510618 x^{13} - 
    334368 x^{14} + 46862 x^{15} - 3416 x^{16} + 72 x^{17}) y^8 - 
  (1-x)^2 (-17136  + 240688 x
 - 1579636 x^2 + 6285816 x^3 - 16975856 x^4 + 33355749 x^5 - 50170572 x^6 + 
    60137293 x^7 - 59249831 x^8 + 48829663 x^9 - 33627193 x^{10} + 18983318 x^{11} - 8514558 x^{12} + 2920248 x^{13} - 
    729963 x^{14} + 123709 x^{15} - 12407 x^{16} + 524 x^{17}) y^9
-  (1-x) (37128 - 528528 x + 3424538 x^2 - 13385812 x^3 + 35510647 x^4 - 68394451 x^5 + 100101788 x^6 - 
    115262274 x^7 + 107466746 x^8 - 83069604 x^9 + 53987064 x^{10} - 29441819 x^{11} + 13202121 x^{12} - 4701667 x^{13} + 
    1271236 x^{14} - 245731 x^{15} + 30752 x^{16} - 2002 x^{17} + 28 x^{18}) y^{10} + 
  (63648 - 917696 x + 5965232 x^2 - 23360456 x^3 + 62122407 x^4 - 119859819 x^5 + 175071125 x^6 - 199526930 x^7 + 
    181725867 x^8 - 135067367 x^9 + 83428883 x^{10} - 43315180 x^{11} + 18849138 x^{12} - 6726629 x^{13} + 1892384 x^{14} - 
    397178 x^{15} + 57372 x^{16} - 4893 x^{17} + 152 x^{18} + 4 x^{19}) y^{11} + 
  (-87516 + 1189760 x - 7206914 x^2 + 26101372 x^3 - 63794722 x^4 + 112416740 x^5 - 148870993 x^6 + 
    152402508 x^7 - 123199222 x^8 + 80108387 x^9 - 42643940 x^{10} + 18849138 x^{11} - 6931650 x^{12} + 2077871 x^{13} - 
    484525 x^{14} + 81618 x^{15} - 8776 x^{16} + 446 x^{17} + 2 x^{18}) y^{12} + 
  (97240 - 1241240 x + 7002476 x^2 - 23476728 x^3 + 52823226 x^4 - 85191730 x^5 + 102552835 x^6 - 94614843 x^7 + 
    68159388 x^8 - 38932682 x^9 + 17903788 x^{10} - 6726629 x^{11} + 2077871 x^{12} - 519020 x^{13} + 99499 x^{14} - 
    13256 x^{15} + 991 x^{16} - 16 x^{17} - 2 x^{18}) y^{13} + 
  (-87516 + 1043900 x - 5472038 x^2 + 16957822 x^3 - 35081419 x^4 + 51716269 x^5 - 56520264 x^6 + 46940030 x^7 - 
    30104049 x^8 + 15085017 x^9 - 5972903 x^{10} + 1892384 x^{11} - 484525 x^{12} + 99499 x^{13} - 15652 x^{14} + 
    1683 x^{15} - 86 x^{16}) y^{14} + (63648 - 705536 x + 3423376 x^2 - 9772052 x^3 + 18515467 x^4 - 24837930 x^5 + 
    24516909 x^6 - 18219941 x^7 + 10332058 x^8 - 4503883 x^9 + 1516967 x^{10} - 397178 x^{11} + 81618 x^{12} - 
    13256 x^{13} + 1683 x^{14} - 156 x^{15} + 6 x^{16}) y^{15} + 
  (-37128 + 380016 x - 1697802 x^2 + 4439228 x^3 - 7652806 x^4 +
$\end{minipage}
\noindent
\begin{minipage}[t]{\linewidth}$
  9266605 x^5 - 8180675 x^6 + 5376370 x^7 - 
    2657724 x^8 + 989788 x^9 - 276483 x^{10} + 57372 x^{11} - 8776 x^{12} + 991 x^{13} - 86 x^{14} + 6 x^{15}) y^{16} + 
  (17136 - 160720 x + 656508 x^2 - 1559760 x^3 + 2421888 x^4 - 2613261 x^5 + 2030050 x^6 - 1155943 x^7 + 
    485274 x^8 - 149047 x^9 + 32754 x^{10} - 4893 x^{11} + 446 x^{12} - 16 x^{13}) y^{17} + 
  (-6120 + 52120 x - 192822 x^2 + 411534 x^3 - 566923 x^4 + 534144 x^5 - 355377 x^6 + 169105 x^7 - 57326 x^8 + 
    13455 x^9 - 2030 x^{10} + 152 x^{11} + 2 x^{12} - 2 x^{13}) y^{18} + 
  (1632 - 12480 x + 41276 x^2 - 77826 x^3 + 92923 x^4 - 73915 x^5 + 40129 x^6 - 14871 x^7 + 3632 x^8 - 524 x^9 + 
    28 x^{10} + 4 x^{11}) y^{19} - 2 (1-x)^4 (153 - 424 x + 392 x^2 - 155 x^3 + 36 x^4) y^{20}
 + 4 (1-x)^5 (9 - 8 x + x^2) y^{21} - 2 (1-x)^5 y^{22}
$\end{minipage}
}

\section{Enumerating 2-convex polygons} \label{s_2-convex}

There is a choice of side when the indent occurs at the edge, that is, touching
the MBR. We choose to define this according to the relative position of the
indents, such that they occur either on the same side, or the closest possible
sides. This is relevant to Cases 3 and 4, but is in contrast to the choice made
in Cases 3 and 4 of the previous section. However, the choice is
arbitrary and does not make a difference to the methodology.

In our `divide-and-conquer' approach, we factor the polygons along lines that
run along the base of the indentations. When the indents are in different
directions, we can therefore divide the plane into four quadrants, by extending
these factorisation lines.  The first quadrant is in the top-right, and they
are ordered in an anti-clockwise fashion.

%
%%
%%% Case 1
%%
%
\subsection{Case 1: indents in the same direction, on the same edge} 

If the indents of a 2-convex polygon are in the same direction and are on the
same edge, then they can be further classified according to whether they are on
the same side or not. We will assume without loss of generality that they are both
vertical and on the top edge.

If one or both of the indents lies along the top edge of the MBR, then they
could be considered to be on either side of the polygon. We choose by default
that the indents be classified as on different sides, where possible, as this
will force us to make an adjustment to the easier of the two cases. Moreover,
when the indents lie on the same side, there is a symmetrical case (e.g. when
the indents are on the top-left and top-right sides) for which this choice
avoids double-counting.

\subsubsection{Indents on the same side} 

\begin{figure}
	\begin{center}
		\includegraphics[scale=0.55]{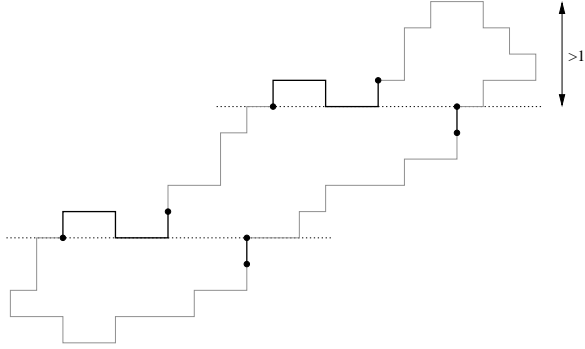}
	\end{center}
	\caption{The form of a 2-convex polygon with both indents on the same side.}
	\label{Fig:2-c.same}
\end{figure}

If the indents of a 2-convex polygon are both on the same side, then they are
either at the same height, or they are not. Due to our choice that the indents
be considered on different sides where possible, the top factor must be of at
least height two.
We refer to Figure~\ref{Fig:2-c.same} for the form of these polygons, which is
the same as the unimodal case, except that the bottom factor can now be
unimodal, not just staircase.
We enumerate the top factor as unimodal, and use wrapping to
include the case where the middle or bottom factor touches the right edge of
the MBR. This is a perfect example of how versatile wrapping is -- there are
even exclusion cases that intersect in the bottom \emph{right} corner, with a
2-unimodal self-avoiding factor that has both indents in the corner.

We adopt an inclusion-exclusion argument, allowing the unimodal factors to
intersect, and then exclude the ones that factor as a 2-unimodal intersecting
polygon and a unimodal SAP. Noting that we must now have a top factor of at
least height two, we obtain these generating functions from
Section~\ref{ss_case1} \mm.

\subsubsection{Indents on adjacent sides} 

When the two vertical indents of a 2-convex polygon are both on top, but on
adjacent sides, their form is very similar to the previous case, where the
indents are on the same side.  As outlined in the unimodal case 1 (c.f.
Figure~\ref{Fig:2-u.same}), once we wrap the polygon such that the middle factor
is the one that touches the right edge of the MBR and the top factor is a
pyramid, we can take thse polygons and flip top factor and indent (that is,
along the top join) such that
the indents are on adjacent sides (c.f. Figure~\ref{Fig:2-u.same}(d)). In this
case the top-right indent is above the top-left one. Similarly, if the flip is
done along the bottom join, we obtain the case where the top-right indent is
below the top-left one (c.f. Figure~\ref{Fig:2-u.same}(a)). Obviously, now that
the polygon is not rooted in one corner, these two cases are symmetrical.
The one difference in this case is that now the top indent can extend furthest
to the right, and may touch the right edge of the MBR. We therefore need to
include this case as well.

Proceeding as above, \mm, the generating function when the indents are at
different heights is therefore
\begin{multline*}
	\lr{ \EE{sy^2}{1-x-y}{\lr{\frac x{x-s} - \frac{1-y}{1-y-s}}}
	+ \frac{P(s,y)}{1-x} } \frac{s^2}{(1-s)^3} \\
	\odot_s\ \bar T(s,t)\ \odot_t\ 	\frac 1{1-t} \Ex y{ty^2}{1-t-y}
	{\lr{\frac{(1-y)^2}{(1-y)^2-x}+\frac x{1-x}} \lr{\frac x{1-x-y}}^2\,} \\
- \frac{2xy}\Delta \cdot \T y^3 \dd y \frac {\T \SP Z}y \lr{u Z + \frac 1{1-x}},
\end{multline*}
and the case where the indents are at the same height is obtained similarly.

\subsubsection{The generating function.}
The above two generating functions combine, with a multiple of two for the
symmetry of the indents being on the same side, to give the generating
function for 2-convex polygons with both vertical indents on top. This result
then gives the overall generating function
\begin{multline} \label{Eq:2-c.same}
	\frac{2A(x,y)}{x^2y^2(1-x)^2(1-y)^2((1-y)^2-x)((1-x)^2-y)\Delta^{7/2}} \\
  + \frac{2B(x,y)}{x^2y^2(1-x)^2(1-y)^2((1-y)^2-x)^2((1-x)^2-y)^2\Delta^{4}},
\end{multline}
where $A(x,y)$ and $B(x,y)$ can be found in \cite{thesis}. We note that
$A(x,y)$ is $O(x^{17}y^{17})$, $A(x,x)$ is $O(x^{20})$,
$B(x,y)$ is $O(x^{21}y^{21})$, $B(x,x)$ is $O(x^{26})$.

%
%%
%%% Case 2
%%
%
\subsection{Case 2: indents in the same direction, on opposite edges} 

The case where the indents are in the same direction, but on opposite edges, is
equivalent to Case 2 of the 2-unimodal polygon enumeration of the previous
section, whose form is depicted in Figure~\ref{Fig:2-u.opp}. Now the polygon is
no longer rooted in the bottom left corner, and so we break up the enumeration
into two parts that are equivalent to those of the last section.

\subsubsection{Indents on opposite sides} 

\begin{figure}
	\begin{center}
	\begin{tabular}{ccc}
		\subfigure[The middle and top factors touch the left and right edges respectively.]
		{\includegraphics[scale=0.5]{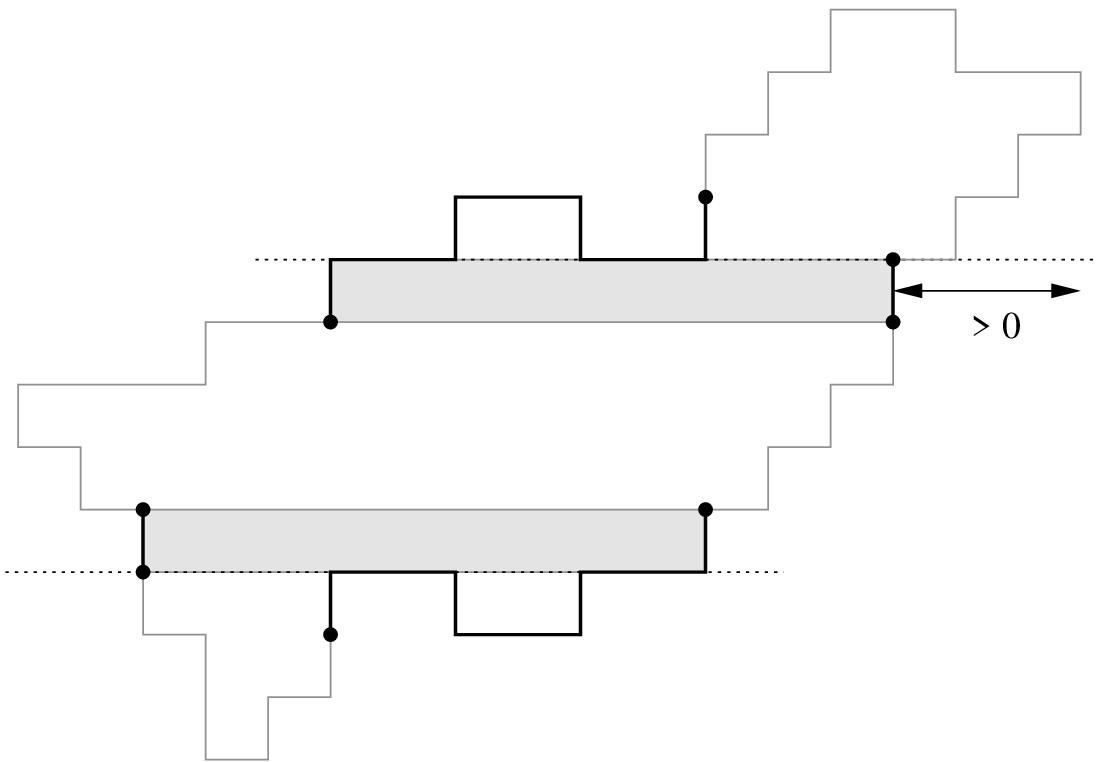}} &\qquad&
		\subfigure[The middle factor touches both the left and the right edges.]
		{\raisebox{6pt}
		 {\includegraphics[scale=0.5]{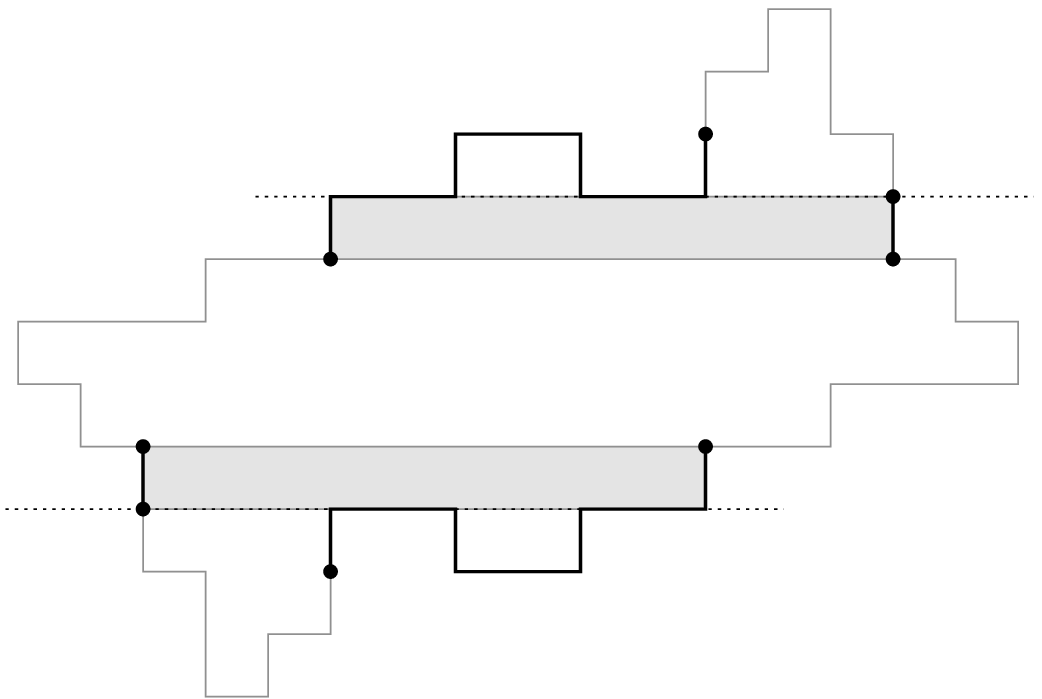}}}\\
		\subfigure[The top indent touches the left edge and the top factor
			touches the right.]
		{\includegraphics[scale=0.5]{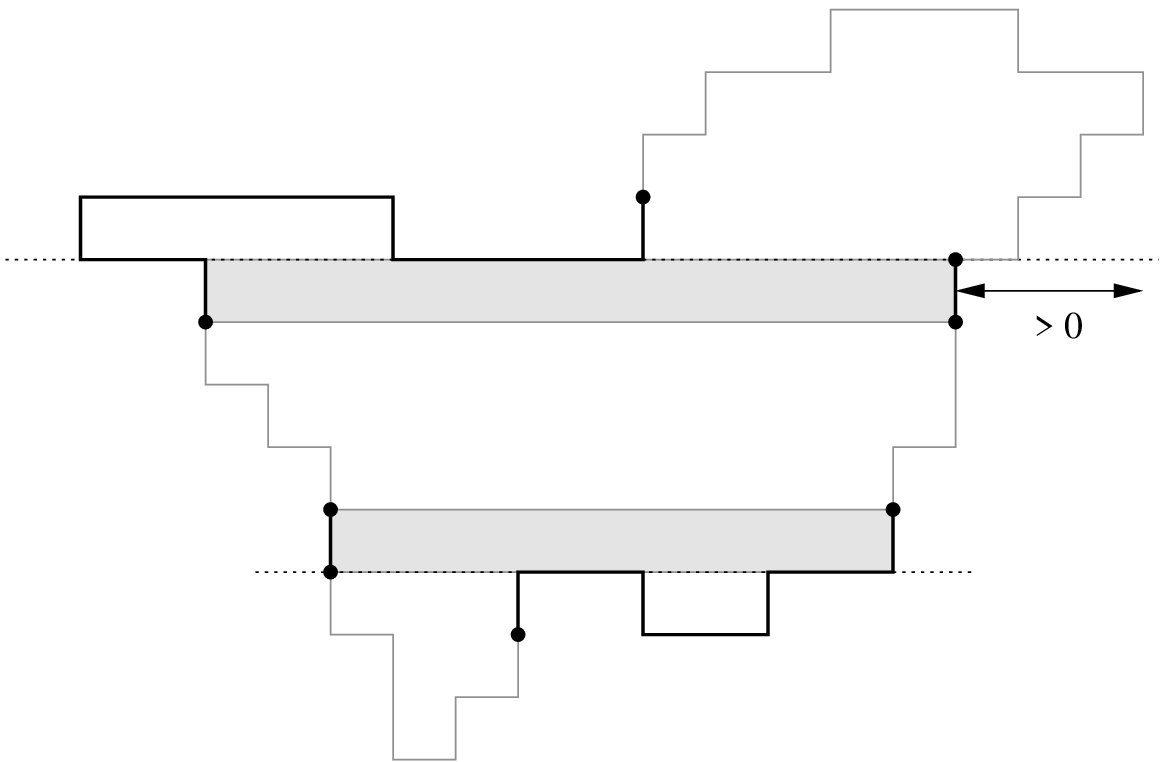}}&&
		\subfigure[The top indent touches the left edge and the middle factor
		touches the left.]
		{\raisebox{0pt}
		{\includegraphics[scale=0.5]{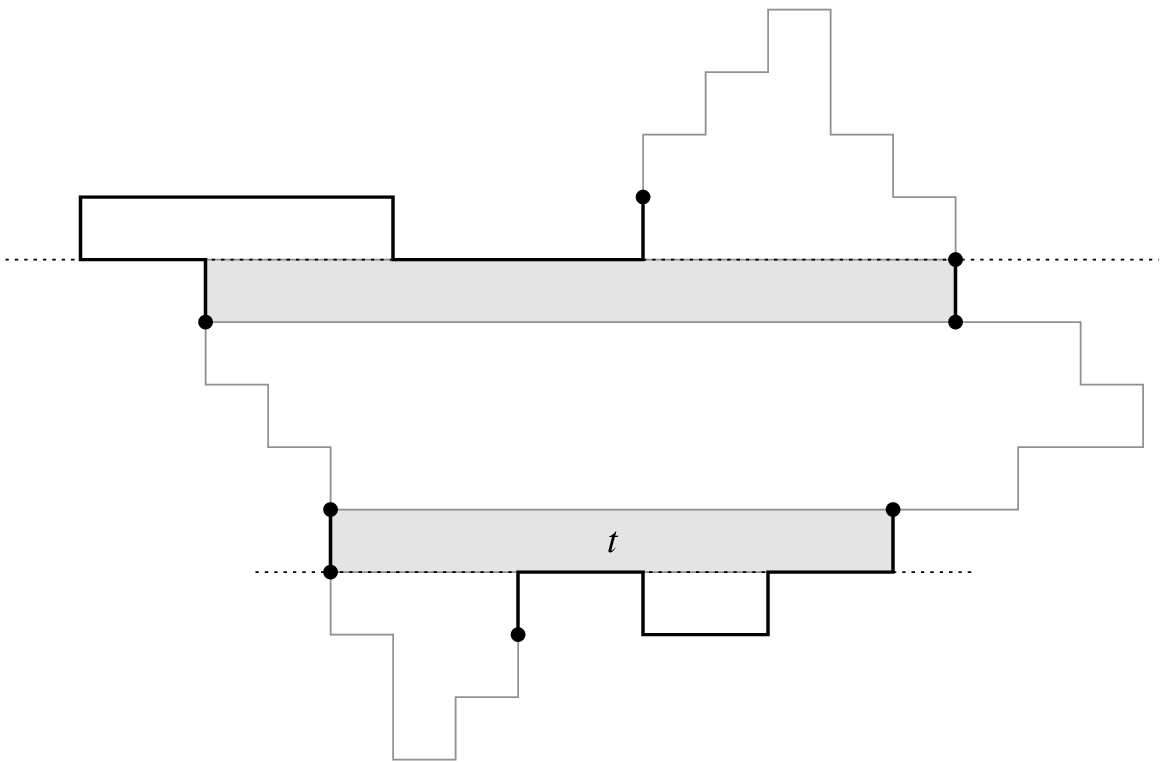}}}\\
	\end{tabular}
	\end{center}
	\caption{The form of 2-convex polygons with indents on opposite sides.}
	\label{Fig:2-c.opp}
\end{figure}

When the indents of a 2-convex polygon are on opposite sides, its form is as
depicted in Figures~\ref{Fig:2-u.opp}(a)--(e) and (g), except that the bottom left
factor is unimodal.
We note that the form shown in part (f) of the figure has the indents on
opposite edges, but adjacent sides, and does not belong to this case.
When the bottom indent is below the top indent, the form is
as shown in part (a) of the figure. Proceeding as per Section~\ref{ss_2-unimodal_opp},
we join the top unimodal factor and indent to the middle staircase factor,
wrapping the resulting polygon such that the middle factor may touch the
right edge of the MBR. However, again, when the base of the staircase factor is
folded, as in part (b) of the figure, an adjustment must be made. The top two
thirds of the polygon is therefore enumerated as in equation
\eqref{Eq:2-u.opp.below}.

What is different about the convex case is that we cannot simply apply the same
logic to join the bottom factor as before, when we could adjust
for the folding of the base \emph{before} applying the join. Now, if we intend
adding a unimodal factor to the bottom that can be wrapped, such that it becomes
a pyramid, we must then remove the cases where the top indent is folded.
Consideration of this case will reveal that the top join requires that the
top-right unimodal factor must not be folded, or the join will have a factor
$s/(x-s)$
on one side, and $s^2x/(1-s^2x)$ on the other, which causes the result to
diverge. Enumerating the cases that were enumerated by the folded factor
separately makes the enumeration just as complicated as if we never folded the
bottom factor in the first place! And so we don't, but instead enumerate
the two cases determined by whether or not the middle factor touches the left
edge of the MBR.

We therefore proceed by enumerating the cases where the bottom-left unimodal
factor touches the left edge of the MBR by following the unimodal cases in
Section~\ref{ss_2-unimodal_opp} \mm.
The case where it is the middle factor that touches the left edge can be
divided into two parts, as shown in Figure~\ref{Fig:2-c.opp}(a) and (b), which
are also generated by standard techniques. We note that the latter case is the
first instance where we require the convex generating function from
Section~\ref{ss_convex_top_bottom} that is enumerated according to both its top
and bottom perimeter.
We also need to generate the symmetrical cases where either indent (but not both)
touches the side. Their form is shown in Figures~\ref{Fig:2-c.opp}(c) and (d).

\subsubsection{Indents on adjacent sides} 

The second part of the enumeration of 2-convex polygons with their indents on
opposite edges is when they are on adjacent sides. This is equivalent to the
unimodal cases depicted in Figures~\ref{Fig:2-u.opp}(f) and (h). We enumerate
them by taking the previous classes with indents on opposite sides, and flipping
the bottom factor and indent.

Now, let us consider the opposite-side cases shown in Figure~\ref{Fig:2-c.opp}.
It is evident that if the bottom factor and indent are reflected such that the
indent is to the left, it will now be on the adjacent side. Furthermore, we
notice that parts (c) and (d) of the figure show when the top indent touches the
left-edge of the MBR. If we have swapped the bottom factor and indent in these
cases, it could be the bottom factor that touched the right edge of the MBR.
There is therefore an extra case to enumerate.

There is another fact that we must note -- that if the top or bottom factors are
of height one, then they were enumerated above, as the indent could be
considered to be on either of two sides. We must therefore re-do the enumeration
of the polygons of Figure~\ref{Fig:2-c.opp}, but with factors of minimum height
two.
(We note that it is helpful to first enumerate the extra cases with factors of
minimum height one and then combine them with the results from above to verify
the total against series which include the double-counted polygons.)

We finish by noting that the interweaved case, equivalent to the unimodal
polygons shown in Figure~\ref{Fig:2-u.opp}(f), does not follow as simply from
the unimodal enumeration as the other cases. By making a unimodal factor wrap
a pyramid bottom factor, we obtain the generating function
\begin{multline}
	\frac{x^2(1+x)}{(1-x)^3} \bigg(
		\lr{\EE{xy}{1-x-y}{\cdot\frac y{1-x}\cdot\frac s{x-s}}
			- \frac{2xy}\Delta \cdot \frac{sv}{1-s-v}} \times
\\		\lr{\frac s{1-s}}^2 \ \odot_s\ 
		\lr{P(xs,y) - \frac{xsy}{1-xs}}
		- \frac{2xy}\Delta \cdot v \T \lr{P(u,y)-\SP}
	\bigg).
\end{multline}

\subsubsection{The generating function.}

Summing the above two results gives the generating function for polygons with
two vertical indents on opposite edges, with one being in the top-left corner.
To obtain the generating function for all those with two vertical indents on
opposite edge we multiply by two, due to symmetry. However, we are now
double-counting the case where the top and bottom indents touch the top and
bottom of the MBR respectively. Denoting \(
	\Phi_x^n = (\delta/\delta x)^n / n!,
\) this case is generated by
\begin{multline}
	\Phi_s^4 \Phi_t^4 \lr{st C(s,t)}
	+ \frac{4x}{1-x} \Phi_s^3 \Phi_t^4 \lr{t \bar U(s,t)} \\
	+ 2\lr{\frac{x}{1-x}}^2 \lr{
		\Phi_s^3 \Phi_t^3 \, \bar T(s,t)
		+ \Phi_s^2 \Phi_t^4 \lr{t P(xs,t)/s}
	}.
\end{multline}
Then, a change of variables gives the generating function for two horizontal
indents, which in turn gives the following generating function for 2-convex
polygons with indents on opposite edges:
\begin{multline} \label{Eq:2-c.opp}
	\frac{-2 A(x,y)}{x^2y^2 (1-x)^3 (1-y)^3 \Delta^{7/2}} \\
  + \frac{B(x,y)}{x^2y^2 (1-x)^7 (1-y)^7 ((1-x)^2-y)^2 ((1-y)^2-x)^2 \Delta^{4}},
\end{multline}
where $A(x,y)$ and $B(x,y)$ can be found in \cite{thesis}. We note that
$A(x,y)$ is $O(x^{15}y^{15})$, $A(x,x)$ is $O(x^{19})$,
$B(x,y)$ is $O(x^{26}y^{26})$, $B(x,x)$ is $O(x^{36})$.

%
%%
%%% Case 3
%%
%
\subsection{Case 3: indents in different directions, on the same edge}

We now consider two out of three cases which enumerate 2-convex polygons
with indents in different directions: where the indents are on the same side and
where they are on adjacent sides.
This is equivalent to Case 3 of the 2-unimodal polygon enumeration of the
previous section. When they are on the same side, they form a locally
concave or convex region, as per Figures~\ref{Fig:2-u.adj}(c) and (d)
respectively. When they are on adjacent sides, their form is equivalent to that
shown in Figure~\ref{Fig:2-u.xyo}.

\subsubsection{Locally concave}

If the polygon's indents form a locally concave region, they are of the form
depicted in Figure~\ref{Fig:2-u.adj}(c). Aside from the walk between the two
indents (which is in
the first quadrant by definition) the polygon can either enter the first or the
third quadrant, but not both. Those polygons which do not enter either quadrant,
passing through the origin, we will consider to be passing through the first
quadrant.

\paragraph{Third quadrant.}

Let us first consider those polygons that enter the third (bottom-left) quadrant.
The indents are both in the first quadrant, within which the polygon is
enumerated by
\(
	Q_1(s,t) = s^2t^2(s+t)/(1-s)(1-t)(1-s-t).
\)
Considering the factor in the second quadrant, we define
\[
	Q_2(x,y;s,t) = \lr{\frac{y}{(1-x)^2-y} + \frac 1{1-y}} \frac{ty}{1-x-ty}
					\cdot \frac{x^2s}{1-xs},
\]
and, by symmetry, $Q_4(x,y;s,t) = Q_2(y,x;t,s)$. And so, as the factor in the
third quadrant is simply a directed walk, we have the following inclusion case:
\[ E\bigg[
		Q_2(x,y^*;p,t) Q_2(y,x^*;q,s) \odot_{p,q,s,t}\ Q_1(s,t)
		\cdot \frac{pq}{1-p-q}
	\bigg].
\]
Wrapping then forms all of the required polygons, including those cases where
the polygon is unimodal in the thrid quadrant. The polygon may therefore
intersect in any of the second, third or fourth quadrant, forming unimodal
loops. Completing the inclusion-exclusion argument by removing those intersecting
polygons, enumerated as per the incusion case \mm, gives the desired result.

\paragraph{First quadrant.}

We will now consider those polygons that pass through the origin, or enter the
first quadrant. Aside from the indents, the factor in the first quadrant is
staircase. We must join unimodal factors to the bottom and to the left, unless
the indent factors touch the MBR, in which case, the second and fourth quadrant
factors must be pyramids. As the unimodal factors can intersect themselves in
the second and fourth quadrants, we adopt the usual inclusion-exclusion
argument, and the generating function follows via standard techniques.

\subsubsection{Locally convex}

When the indents form a convex region, the `humps' of the indents are either
distinct, in which case the polygon does not enter the first quadrant, or they
are joined by a common cell contained in the first quadrant. In the first case,
the third quadrant factor is staircase, and otherwise it is formed simply by a
directed walk. To this we join two unimodal factors in the second and fourth
quadrant. The unimodal factors may intersect, and we use a standard
inclusion-exclusion approach. `Wrapping' the factors together ensures that all
possible configurations of the polygons are formed. This means that despite the
fact the factor in the third quadrant is staircase, when it is wrapped, it may
intersect itself, forming a unimodal loop in the bottom-left corner. we
therefore also need to exclude this case.

To help enumerate the unimodal factors of the second and fourth
quadrant, let
\[
	Q_2(s; n) = \frac{s^n}{(1-s)^{n+1}} \cdot \frac{sP(x,y^*)}{x-s}, \quad
	Q_2^*(s; n) = \frac{s^n}{(1-s)^{n+1}} \cdot \frac{s v^*}{1-s-v^*},
\]
and we define $Q_4$ and $Q_4^*$ similarly, with $x$ and $y$ interchanged.
The inclusion case is therefore enumerated by
\begin{multline*}
	E\bigg[ Q_2(s;2) Q_4(t;2) \bigg] \odot_{s,t} F(s,t) \\
+	E\bigg[ Q_2(s;1) Q_4(t;1) - Q_2(s;0) Q_4(t;0) \bigg]
	\odot_{s,t} \frac{xsty}{1-xs-ty},
\end{multline*}
and the other terms follow \mm.

\subsubsection{The generating function for indents on the same side.}
Summing the above two results for the polygons that are locally concave and
convex around their indents gives the following generating function for polygons
with indents in different directions on the same side:
\begin{equation} \label{Eq:2-c.with}
	\frac{x y A(x,y)}{2 (1-x)   (1-y)   \Delta^{7/2}}
  + \frac{x y B(x,y)}{2 (1-x)^2 (1-y)^2 \Delta^{4}},
\end{equation}
where $A(x,y)$ and $B(x,y)$ can be found in \cite{thesis}. We note that
$A(x,y)$ is $O(x^{8}y^{8})$, $A(x,x)$ is $O(x^{9})$,
$B(x,y)$ is $O(x^{10}y^{10})$, $B(x,x)$ is $O(x^{12})$.

\subsubsection{Indents on adjacent sides}

\begin{figure}
	\begin{center}
		 \includegraphics[scale=0.5]{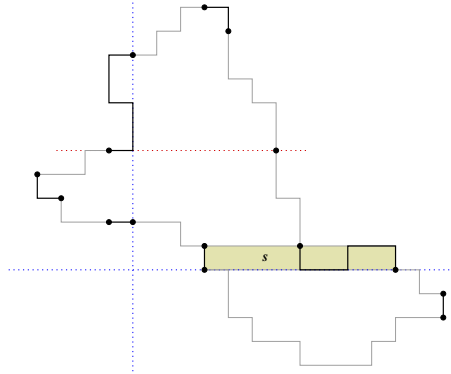}
	\end{center}
	\caption{The form of 2-convex polygons with indents in different directions
			 on adjacent sides of the polygon and which touch the left edge of
			 the MBR in the second quadrant.}
	\label{Fig:2-c.adj.above}
\end{figure}

We enumerate this case by braking it up into parts that are
classified by the relative height of the indent.
As per the previous section, we say that the indents are {\em next to} each
other if the {\em indentations} formed overlap in height. If the {\em hump} of
the indent factor, rather than the indentation, is at the same height as the
vertical indent, we say that they are {\em adjacent}.  And so, the horizontal
indent is either above, next to, adjacent to or below the vertical indent. The
form of the equivalent unimodal polygons are shown in Figure~\ref{Fig:2-u.xoy}.
However, there is one case missing that was not relevant for the unimodal
situation; when the horizontal indent is above the vertical one, and the polygon
touches the MBR in the second quadrant. This case is depicted in
Figure~\ref{Fig:2-c.adj.above}.

When the horizontal indent is higher, the polygons are of one of the forms shown
in Figures~\ref{Fig:2-u.xyo}(a) or \ref{Fig:2-c.adj.above}. Let us first
enumerate the latter case, which has a significantly different form to the other
polygons of this class. We factorise as usual, which divides the polygon into
two main factors: the top and bottom sections.

We begin by enumerating the top section. We need to enumerate it according the
its base. We have not already needed to enumerate 1-unimodal polygons whose
indent is in the corner according to their base. And so, we break this section
of the polygon into two parts again by extending a line from the bottom of the
indentation (shown in red in the figure). We therefore join an almost-pyramid
polygon to a unimodal polygon and obtain the generating function
\begin{equation} \label{Eqn:1-u.corner.base}
	U_1^{c}(s) = \lrf y{1-t-y}^2 \lrf t{1-t}^2
		\lr{\frac{y}{(1-t)^2-y} + \frac{1}{1-y}} \odot_t \bar{U}(s,t),
\end{equation}
where $s$ counts the base. We then join the indent factor by multiplying by
$(xs/(1-xs))^2$ and join it to the bottom section.

Now the bottom section is unimodal, but if it is pyramid, the indent factor can
extend further to the right. We therefore add the pyramid and non-pyramid
factors separately. We do the join in the most straight-forward manner by
overlapping one row of the top and bottom factors (indicated by the shaded
region in the figure), which requires a division of the generating function by
$y$. In this case, we run into the case where the polygon intersects itself. We
must therefore exclude this intersecting case.
The generating function is therefore given by
\begin{equation}
	\frac 1y \lr{U_1^{c}(s) - U_1^{c}(1)} \lrf{xs}{1-xs}^2 \odot_s
	\bar U(s/x) + P(s,y) \lr{\frac{1-s}{1-x} - 1},
\end{equation}
where $\bar U(s)$ is the unimodal polygon generating function, with $s$ counting
the bottom perimeter, and $P(x,y)$ is the pyramid generating function.

If the horizontal indent occurs above the vertical one, but the polygon touches
the left edge of the MBR in the third quadrant, its form is as depicted in
Figure~\ref{Fig:2-u.xyo}(a), except that the bottom factor is convex. The case
where the left perimeter is in the second quadrant is obtained via wrapping.
The remaining three cases are all obtained by a similar extension of the
unimodal case \mm.

\subsubsection{The generating function.}
In the above derivation, we have taken the position and direction of one indent,
and then assumed that the other indent is one the next side, following the
polygon in one of the two possible orientations. However, those polygons whose
indent is on the other adjacent side also belong to this class. We obtain their
generating function by swapping the variables $x$ and $y$. And so, the total
generating function is
\begin{equation} \label{Eq:2-c.against.adj}
	\frac{A(x,y)}{(1-x)(1-y)\Delta^{5/2}}
  - \frac{B(x,y)}{(1-x)(1-y) ((1-x)^2-y) ((1-y)^2-x) \Delta^4},
\end{equation}
where $A(x,y)$ and $B(x,y)$ can be found in \cite{thesis}. We note that
$A(x,y)$ is $O(x^{8}y^{8})$, $A(x,x)$ is $O(x^{9})$,
$B(x,y)$ is $O(x^{14}y^{14})$, $B(x,x)$ is $O(x^{16})$.

%
%%
%%% Case 4
%%
%
\subsection{Case 4: indents in different directions, on opposite edges} 

We note that if the indents are on opposite edges, they must be on opposite
sides, due to symmetry.  We again enumerate this case by braking it up into
parts that are classified by the relative height of the indent, adopting the
same terminology.  The form of the equivalent unimodal polygons are shown in
Figure~\ref{Fig:2-u.xoy}.
The only difference is that the factor in the third quadrant can now be unimodal.
We note, however, that the indents must not lie along the edge of the MBR,
as these polygons have already been enumerated, in contrast to the unimodal case.

\begin{figure}
	\begin{center}
		 \includegraphics[scale=0.5]{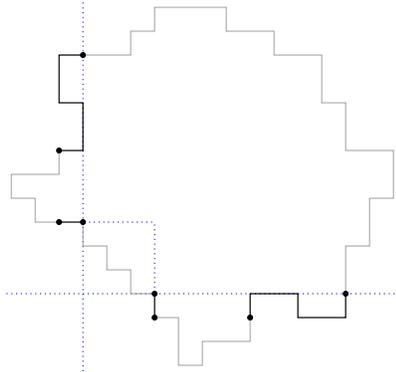}
	\end{center}
	\caption{The form of 2-convex polygons with indents in different directions
			 on opposite sides of the polygon that do not enter the third
			 quadrant.}
	\label{Fig:2-c.opp.above}
\end{figure}

As was the case when the indents were on adjacent sides, there is one new case
to enumerate. This is when the polygon does not enter the thrid quadrant, as
shown in Figure~\ref{Fig:2-c.opp.above}.
This is an interesting case, as we cannot proceed as normal, which is because
the horizontal and vertical joins (over $s$ and $t$ respectively) must be done
{\em simultaneously}\footnote{For an illustration of why this is, see
\cite{thesis}.}.
We are forced to re-evaluate our approach and come up with a new way to evaluate
such constructions. 
Using wrapping, we can make a unimodal factor on top by making the top factor
staircase and add double-bonds to the bottom factor, which wrap the staircase
factor. Each bond has two steps,
contributing $s^2$ to the weight of the join, as well as a weight of $1/x$.
This leaves us with the following expression for the bottom factor in the join:
\[
\frac{1}{1-s} \EE{sy^2}{1-s-y}{\frac{y}{1-s}} \left(\frac{s}{1-s}\right)^2
\left( \frac{1}{1-x} + \frac{s^2}{x-s^2} \right),
\]
where the two terms in the last factor correspond respectively to the cases where
the indent touches the right-edge and where only the convex factor touches it.

Up until now it has been possible to do joins with such factors in them by
simplifying other terms and applying the equivalent operations to them. Now we
will have two distinct factors that both have expressions of this form, and so
we cannot avoid evaluating it. We therefore need to re-express these terms
in a form which we know how to evaluate, which we achieve by noting that
\begin{equation}
	\frac{s^2}{x-s^2} = \Ex{x}{s}{x-s}{}.
\end{equation}
This therefore achieves the required wrapping by only applying the $E$ operator
to the 1-dimensional loop part of the join, and not the indent factor.

We now apply this approach to our current problem by joining the indents and
pyramid factors to a staircase polygon that is enumerated by its bottom and side
perimeters. `Folding' the staircase factor twice creates a convex polygon,
gives the possibility of intersecting itself. We will therefore exclude this
case, which is straight-forward to enumerate. Keeping in mind that the pyramid
factors must be of at least height two, we obtain the following generating
function:
\begin{multline}
  \EE{sy^3}{1-s-y}{\cdot\frac{s^2}{(1-s)^4}\lr{\frac{s}{x-s}-\frac{1}{1-x^*}}}
  \odot_s F(s,t) \ \odot_t \\
  \EE{x^3t}{1-x-t}{\cdot\frac{t^2}{(1-t)^4}\lr{\frac{t}{y-t}-\frac{1}{1-y^*}}}
  - (\SP Z)^3\cdot \T \cdot \bar\T \cdot \frac{2xy}{\Delta}.
\end{multline}

The reamining four cases are equivalent to those depicted in
Figure~\ref{Fig:2-u.xoy}. Using the above application of wrapping next to
the indents, we obtain the generating functions by enumerating the first
quadrant as a directed walk.
We note that the unimodal factor in the third quadrant must be of at least width
two, and of at least height two when the top indent is above the bottom one.
The latter case is generated by
\begin{multline}
	E_{w,x,y,z}\bigg[\frac{(wz)^2}{1-w-z}
	\lrf{s}{1-s}^2 \frac{s}{1-s-z} \lr{\frac{s}{x-s} + \frac{1}{1-x^*}} \\
		\odot_s \frac{x^*y^*}{1-x^*s-ty^*} \ \odot_t
	\lrf{t}{1-t}^2 \frac{t}{1-w-t} \lr{\frac{t}{y-t} + \frac{1}{1-y^*}}
	\bigg] \bigg|_{w=x, z=y} \\
  - \T\cdot{\bar\T}\cdot\frac{2xy}{\Delta} \cdot
    \Ex{w,z}{(wz)^2}{1-w-z}{\cdot\frac{v^*}{1-w-v^*}\cdot\frac{u^*}{1-u^*-z}},
\end{multline}
for minimum height and width one, and the required result follows the exclusion
of the height or width one cases. The remaining three results follow \mm, and we
refer to \cite{thesis} for the details of the calculations.

\subsubsection{The generating function.}
The same symmetries in the $x=y$ axis exist for the 2-convex case as for the
2-unimodal one, and it is necessary to reflect the last three, non-symmetric parts
and add them to the other classes to get the all possible 2-convex polygons with
a vertical indent on the top-left side of the polygon and a horizontal
one on the opposite side. Summation gives the following generating function:
\begin{multline} \label{Eq:2-c.against.opp}
	\frac{A(x,y)}{2(1-x)(1-y)\Delta^{7/2}} \\
  - \frac{B(x,y)}{2 (1-x-y) (1-x)^2(1-y)^2 ((1-x)^2-y)^2 ((1-y)^2-x)^2 \Delta^4},
\end{multline}
where $A(x,y)$ and $B(x,y)$ can be found in \cite{thesis}. We note that
$A(x,y)$ is $O(x^{10}y^{10})$, $A(x,x)$ is $O(x^{11})$,
$B(x,y)$ is $O(x^{19}y^{19})$, $B(x,x)$ is $O(x^{23})$.

\subsection{The 2-convex generating function}

In the above derivations, the direction and position of one of the indents was
chosen arbitrarily, such that the direction and position of the other determined
which sub-class the polygon belongs to. Therefore, when the two indents are in
different directions, the cardinality of the set of possible combinations of
direction and location for the fixed indent is four. When the indents are in the
same direction, both directions are enumerated by the generating functions.  If
the indents are on the same edge, the cardinality is two. If they are opposite,
it is one. We note that we obtained the bimodal 2-convex generating function in
Section~\ref{s_bimodal}.
We therefore obtain the generating function for 2-convex
polygons by summing the results, multiplying each term by the cardinality of the
class it counts.  This gives the following generating function for 2-convex
polygons:
\begin{multline} \label{Eq:2-convex}
	\frac{-4A^{(c)}_2}{(1-x)^3 x^2 (1-y)^3 y^2 \Delta^{7/2}} \\
	- \frac{B^{(c)}_2}{(1-x)^7 x^2 (1-y)^7 y^2 ((1-x)^2-y)^3 ((1-y)^2-x)^3 (1-x-y) \Delta^4},
\end{multline}
{\small
\begin{minipage}[t]{\linewidth}$
\mbox{\normalsize{where }} A^{(c)}_2 = 
\ (1-x)^{11} x^4
- 3 (1-x)^9 x^4 (5 - 2 x + x^2) y
+ (1-x)^7 x^2 (4 - 12 x + 103 x^2 - 79 x^3 + 31 x^4 - 11 x^5 + 3 x^6) y^2
- (1-x)^5 x^2 (40 - 124 x + 455 x^2 - 533 x^3 + 264 x^4 - 69 x^5 - 21 x^6 - 3 x^7 + x^8) y^3
+ (1-x)^3 (1 - 12 x + 232 x^2 - 742 x^3 + 1696 x^4 - 2297 x^5 + 1626 x^6 - 539 x^7 + 133 x^8 + 3 x^9 - 56 x^{10} + 3 x^{11}) y^4
+ (1-x)^2 (-11 + 119 x - 943 x^2 + 2443 x^3 - 4014 x^4 + 4513 x^5 - 3054 x^6 + 867 x^7 - 58 x^8 + 221 x^9 - 137 x^{10} - 11 x^{11} + x^{12}) y^5
+ (1-x) (55 - 542 x + 2765 x^2 - 6154 x^3 + 8193 x^4 - 7901 x^5 + 5521 x^6 - 2140 x^7 - 284 x^8 + 430 x^9 + 19 x^{10} - 81 x^{11} + 7 x^{12}) y^6
+ (-165 + 1503 x - 5996 x^2 + 11929 x^3 - 14004 x^4 + 11488 x^5 -
$\end{minipage}
\noindent
\begin{minipage}[t]{\linewidth}$
7661 x^6 + 4474 x^7 - 1456 x^8 - 506 x^9 + 504 x^{10} - 160 x^{11} + 18 x^{12}) y^7
+ (330 - 2502 x + 7381 x^2 - 10693 x^3 + 8925 x^4 - 4846 x^5 + 1856 x^6 - 1456 x^7 + 1364 x^8 - 328 x^9 + 103 x^{10} - 22 x^{11}) y^8
+ (-462 + 2898 x - 6353 x^2 + 6348 x^3 - 3639 x^4 + 1204 x^5 + 714 x^6 - 506 x^7 - 328 x^8 - 32 x^9 + 10 x^{10}) y^9
+ (462 - 2394 x + 3925 x^2 - 2349 x^3 + 873 x^4 - 637 x^5 - 411 x^6 + 504 x^7 + 103 x^8 + 10 x^9) y^{10}
+ (-330 + 1422 x - 1797 x^2 + 434 x^3 + 47 x^4 + 484 x^5 - 100 x^6 - 160 x^7 - 22 x^8) y^{11}
+ (165 - 603 x + 632 x^2 - 4 x^3 - 180 x^4 - 114 x^5 + 88 x^6 + 18 x^7) y^{12}
- (1-x)^2 (55 - 67 x - 18 x^2 + 27 x^3 + 7 x^4) y^{13}
- (1-x)^3 (-11 + x^2) y^{14}
- (1-x)^3 y^{15}
$\end{minipage}\\[0.5\baselineskip]
\noindent
\begin{minipage}[t]{\linewidth}$
\mbox{\normalsize{and }} B^{(c)}_2 = 
	\ 4 (1-x)^{26} x^4
	- 4 (1-x)^{24} x^4 (30 - 20 x + 7 x^2) y
	+ 4 (1-x)^{22} x^2 (4 - 12 x + 433 x^2 - 580 x^3 + 368 x^4 - 128 x^5 + 21 x^6) y^2
	- 4 (1-x)^{20} x^2 (100 - 360 x + 4214 x^2 - 8192 x^3 + 7616 x^4 - 4297 x^5 + 1528 x^6 - 350 x^7 + 35 x^8) y^3
	+ 4 (1-x)^{18} (1 - 12 x + 1252 x^2 - 5226 x^3 + 32426 x^4 - 76437 x^5 + 93156 x^6 - 70591 x^7 + 36061 x^8 - 12917 x^9 + 3217 x^{10} - 
		    532 x^{11} + 35 x^{12}) y^4
	- (1-x)^{16} (104 - 1296 x + 42936 x^2 - 198752 x^3 + 859020 x^4 - 2164230 x^5 + 3180634 x^6 - 
		    3007804 x^7 + 1964426 x^8 - 927123 x^9 + 322834 x^{10} - 83390 x^{11} + 15842 x^{12} - 1961 x^{13} + 84 x^{14}) y^5
	+ (1-x)^{14} (1300 - 16868 x + 289112 x^2 - 1412308 x^3 + 5065708 x^4 - 12664938 x^5 + 20907814 x^6 - 23411754 x^7 + 18540432 x^8 - 
		    10767204 x^9 + 4716295 x^{10} - 1577652 x^{11} + 403402 x^{12} - 79358 x^{13} + 12073 x^{14} - 1114 x^{15} + 28 x^{16}) y^6
	- (1-x)^{12} (10400 - 140832 x +
$\end{minipage}
\noindent
\begin{minipage}[t]{\linewidth}$
	 1632176 x^2 - 8077604 x^3 + 26714752 x^4 - 64771716 x^5 + 113604188 x^6 - 143655868 x^7 + 
		    132836186 x^8 - 91595522 x^9 + 48102042 x^{10} - 19641626 x^{11} + 6301753 x^{12} - 1574826 x^{13} + 304364 x^{14} - 47486 x^{15} + 
		    6143 x^{16} - 400 x^{17} + 4 x^{18}) y^7
	- (1-x)^{10} (-59800 + 846768 x -
$\end{minipage}
\noindent
\begin{minipage}[t]{\linewidth}$
	7892920 x^2 + 38780872 x^3 - 125239164 x^4 + 297723354 x^5 - 
		    537008734 x^6 + 736108516 x^7 - 767512896 x^8 + 611776002 x^9 - 375779397 x^{10} + 180183183 x^{11} - 68757620 x^{12} + 21225217 x^{13} - 
		    5236463 x^{14} + 983039 x^{15} - 136742 x^{16} + 16579 x^{17} - 1952 x^{18} + 78 x^{19}) y^8
	+ (1-x)^8 (-263120 + 3902272 x - 32692880 x^2 + 159495728 x^3 - 518222904 x^4 + 1235851899 x^5 - 2273834564 x^6 + 3286427610 x^7 - 
		    3744105010 x^8 + 3360169186 x^9 - 2372247082 x^{10} + 1316981454 x^{11} - 578123374 x^{12} + 204616292 x^{13} - 60356681 x^{14} + 
		    15025753 x^{15} - 2968864 x^{16} + 408498 x^{17} - 36191 x^{18} + 3653 x^{19} - 450 x^{20} + 7 x^{21}) y^9
	+ (1-x)^6 (920920 - 14323848 x + 115445704 x^2 - 565821908 x^3 + 1881114592 x^4 - 4595055622 x^5 + 8689038821 x^6 - 
		    13096803662 x^7 + 15935786206 x^8 - 15697620474 x^9 + 12485332431 x^{10} - 7965917640 x^{11} + 4040636749 x^{12} - 1619517767 x^{13} + 
		    518745262 x^{14} - 140030433 x^{15} + 34601374 x^{16} - 7726101 x^{17} + 1290999 x^{18} - 115181 x^{19} + 1171 x^{20} - 68 x^{21} + 43 x^{22}) y^{10}
	+ (1-x)^4 (-2631200 + 42965472 x - 346781248 x^2 + 1732729172 x^3 - 5960460908 x^4 + 15133786878 x^5 - 29768776966 x^6 + 
		    46871962702 x^7 - 60257832295 x^8 + 63871376160 x^9 - 55950160750 x^{10} + 40346497476 x^{11} -
	23700589523 x^{12} + 11139364197 x^{13} - 
		    4087427426 x^{14} + 1145559467 x^{15} - 251575315 x^{16} + 52651169 x^{17} - 13478667 x^{18} + 3310337 x^{19} - 518486 x^{20} + 32936 x^{21} + 
		    909 x^{22} + 2 x^{23} + 3 x^{24}) y^{11}
	- (1-x)^3 (-6249100 + 100966228 x - 786008036 x^2 + 3790621220 x^3 - 12634194004 x^4 + 31120468774 x^5 - 59333283779 x^6 + 90481426644 x^7 - 112836337344 x^8 + 
			116619494723 x^9 - 100526269572 x^{10} + 72257250633 x^{11} -  42935196678 x^{12} + 20646926051 x^{13} - 7717297655 x^{14} + 2081066620 x^{15} - 344551302 x^{16} + 19708982 x^{17} - 190245 x^{18} + 
		    2604586 x^{19} - 1089925 x^{20} + 167990 x^{21} - 5325 x^{22} - 819 x^{23} + 5 x^{24}) y^{12}
	- (1-x)^2 (12498200 - 200432672 x + 1523814896 x^2 - 7173568104 x^3 + 23394202064 x^4 - 56432397385 x^5 + 105258925336 x^6 - 
		    156756419659 x^7 + 190737577732 x^8 - 192647570536 x^9 + 163201860685 x^{10} - 116641434399 x^{11} + 70265808073 x^{12} - 
		    35189781755 x^{13} + 14129858643 x^{14} - 4206300270 x^{15} + 752961328 x^{16} + 2762280 x^{17} -
$\end{minipage}
\noindent
\begin{minipage}[t]{\linewidth}$
			41077210 x^{18} + 8781488 x^{19} + 
		    172960 x^{20} - 372435 x^{21} + 59206 x^{22} - 2116 x^{23} - 89 x^{24} + 11 x^{25}) y^{13}
	- (1-x) (-21246940 + 339297288 x - 2543202516 x^2 + 11796257160 x^3 - 37956243060 x^4 + 90389633588 x^5 - 166278246916 x^6 + 243693863948 x^7 -
			291099521233 x^8 + 288163147560 x^9 - 239402103854 x^{10} + 168782468152 x^{11} - 101938448323 x^{12} + 
		    52836781903 x^{13} - 23089813175 x^{14} + 8076150924 x^{15} - 2006509792 x^{16} + 229306444 x^{17} + 52941397 x^{18} - 29443771 x^{19} + 
		    5233574 x^{20} - 137308 x^{21} - 
$\end{minipage}
\noindent
\begin{minipage}[t]{\linewidth}$
			96084 x^{22} + 13023 x^{23} - 102 x^{24} + 33 x^{25}) y^{14}
	+ (-30904640 + 493047872 x - 3671274848 x^2 + 16909765768 x^3 - 54076869168 x^4 + 128034624944 x^5 - 233957204524 x^6 + 
		    339860131708 x^7 - 401161475162 x^8 + 390960587591 x^9 - 318456383150 x^{10} + 219480309921 x^{11} - 130108934416 x^{12} + 
		    67655799311 x^{13} - 31165964099 x^{14} + 12455845878 x^{15} - 4042728520 x^{16} + 929084638 x^{17} - 93195098 x^{18} - 24059794 x^{19} + 
		    11705293 x^{20} - 1983652 x^{21} + 86582 x^{22} + 23547 x^{23} - 3151 x^{24} - 68 x^{25} + 5 x^{26}) y^{15}
	+ (38630800 - 579105408 x + 4018988224 x^2 - 17153382704 x^3 + 50579891852 x^4 - 109880934764 x^5 + 183275392200 x^6 - 
		    241692656260 x^7 + 257400864147 x^8 - 224539648706 x^9 + 161936707748 x^{10} - 97616424050 x^{11} + 50386570178 x^{12} - 
		    23295420511 x^{13} + 10082660716 x^{14} - 4042728520 x^{15} + 1381926138 x^{16} - 352763968 x^{17} +
$\end{minipage}
\noindent
\begin{minipage}[t]{\linewidth}$
			51556603 x^{18} + 1631397 x^{19} - 
		    2559687 x^{20} + 547928 x^{21} - 54439 x^{22} - 611 x^{23} + 432 x^{24} - 31 x^{25}) y^{16}
	+ (-41602400 + 585761792*x - 3791089760*x^2 + 15002579968*x^3 - 40793586584*x^4 + 81262920502*x^5 - 123537868000*x^6 + 
		    147559428776 x^7 - 141331365760 x^8 + 109736530568 x^9 - 69193200719 x^{10} + 35421383132 x^{11} - 15013860403 x^{12} + 
		    5709460646 x^{13} - 2235816236 x^{14} + 929084638 x^{15} - 352763968 x^{16} + 103409188 x^{17} - 19757856 x^{18} + 1570918 x^{19} + 
		    189358 x^{20} - 71908 x^{21} + 11883 x^{22} - 268 x^{23} + 77 x^{24}) y^{17}
	+ (38630800 - 510758608 x + 3084598320 x^2 - 11324232712 x^3 + 28395037220 x^4 - 51826962892 x^5 + 71682098798 x^6 - 77353955608 x^7 + 
		    66429125667 x^8 - 45706907743 x^9 + 24906199156 x^{10} - 10403200527 x^{11} + 3174037717 x^{12} - 706359558 x^{13} + 176365047 x^{14} - 
		    93195098 x^{15} + 51556603 x^{16} - 19757856 x^{17} + 4723092 x^{18} - 605060 x^{19} + 49569 x^{20} - 3395 x^{21} - 1253 x^{22} - 143 x^{23}) y^{18} + (-30904640 + 383692992 x - 2163792656 x^2 + 
$\end{minipage}
\noindent
\begin{minipage}[t]{\linewidth}$
	7374341864 x^3 - 17053687520 x^4 + 28497776084 x^5 - 35790327136 x^6 + 34800116028 x^7 - 
		    26729224736 x^8 + 16280988977 x^9 - 7648226004 x^{10} + 2524992746 x^{11} - 406853569 x^{12} - 93698188 x^{13} + 82385168 x^{14} - 
		    24059794 x^{15} + 1631397 x^{16} + 1570918 x^{17} - 605060 x^{18} + 66554 x^{19} - 4454 x^{20} + 1819 x^{21} + 202 x^{22}) y^{19}
	+ (21246940 - 247772008 x + 1305971612 x^2 - 4135556872 x^3 + 8823936024 x^4 - 13492690928 x^5 + 15356870874 x^6 - 13416475166 x^7 + 
		    9196907145 x^8 - 4978222336 x^9 + 2052527732 x^{10} - 556811827 x^{11} + 29183400 x^{12} + 58467226 x^{13} - 34677345 x^{14} + 
		    11705293 x^{15} - 2559687 x^{16} + 189358 x^{17} + 49569 x^{18} - 4454 x^{19} - 1368 x^{20} - 110 x^{21}) y^{20}
	+ (-12498200 + 136980272 x - 675634696 x^2 + 1990453776 x^3 - 3921712444 x^4 + 5486790816 x^5 - 5649513210 x^6 + 4421089360 x^7 - 
		    2695886222 x^8 + 1305174655 x^9 - 491084981 x^{10} + 128534739 x^{11} - 11441768 x^{12} - 8063133 x^{13} + 5370882 x^{14} - 1983652 x^{15} + 
		    547928 x^{16} - 71908 x^{17} - 3395 x^{18} + 1819 x^{19} - 110 x^{20}) y^{21}
	+ (6249100 - 64433028 x + 297851620 x^2 - 817733112 x^3 + 1489860872 x^4 - 1908408554 x^5 + 1775572130 x^6 - 1240442976 x^7 + 
		    667928313 x^8 - 289317379 x^9 + 102644601 x^{10} - 29961766 x^{11} + 6383656 x^{12} - 977036 x^{13} - 41224 x^{14} + 86582 x^{15} - 
		    54439 x^{16} + 11883 x^{17} - 1253 x^{18} + 202 x^{19}) y^{22}
	+ (-2631200 + 25559072 x - 110940720 x^2 + 284446548 x^3 - 480234432 x^4 + 
		    564112724 x^5 - 474248936 x^6 + 294904852 x^7 - 138225516 x^8 + 52522411 x^9 - 17224508 x^{10} + 5578263 x^{11} - 1609051 x^{12} + 
		    492963 x^{13} - 109107 x^{14} + 23547 x^{15} - 611 x^{16} - 268 x^{17} - 143 x^{18}) y^{23}
	+ (920920 - 8444128 x + 34492144 x^2 - 82799336 x^3 + 129893340 x^4 - 140307634 x^5 + 106798702 x^6 - 59095890 x^7 + 23639299 x^8 - 
		    7486656 x^9 + 2001655 x^{10} - 644781 x^{11} + 181503 x^{12} - 63349 x^{13} + 13125 x^{14} - 3151 x^{15} + 432 x^{16} + 77 x^{17}) y^{24}
	+ (-263120 + 2283072 x - 8798288 x^2 + 19829008 x^3 - 28995448 x^4 + 28906719 x^5 - 20012600 x^6 + 9913898 x^7 - 3347134 x^8 + 
		    826000 x^9 - 124087 x^{10} + 29300 x^{11} - 2853 x^{12} + 1927 x^{13} - 135 x^{14} - 68 x^{15} - 31 x^{16}) y^{25}
	+ (59800 - 492568 x + 1795744 x^2 - 3809804 x^3 + 5207920 x^4 - 4809154 x^5 + 3046581 x^6 - 1369614 x^7 + 399732 x^8 - 78407 x^9 + 
		    2224 x^{10} + 919 x^{11} - 834 x^{12} + 111 x^{13} + 33 x^{14} + 5 x^{15}) y^{26}
	+ (-10400 + 81632 x - 282336 x^2 + 565028 x^3 - 723116 x^4 + 619222 x^5 - 359946 x^6 + 148482 x^7 - 39609 x^8 + 7449 x^9 - 326 x^{10} - 
$\end{minipage}
\noindent
\begin{minipage}[t]{\linewidth}$
		    10 x^{11} + 5 x^{12} - 11 x^{13}) y^{27}
	+ (1300 - 9768 x + 32140 x^2 - 60712 x^3 + 72592 x^4 - 57298 x^5 + 30217 x^6 - 11207 x^7 + 
		    2732 x^8 - 506 x^9 + 43 x^{10} + 3 x^{11}) y^{28}
	+ (-104 + 752 x - 2360 x^2 + 4200 x^3 - 4648 x^4 + 3305 x^5 - 1506 x^6 + 448 x^7 - 
		    78 x^8 + 7 x^9) y^{29}
	+ 4 (1-x)^7 y^{30}.
$\end{minipage}
}

\appendix

\section{Enumerating staircase and unimodal polygons with \mbox{$m$} vertical
indents on the left side} \label{s_m-same}

\begin{figure}
	\begin{center}
		\includegraphics[scale=0.5]{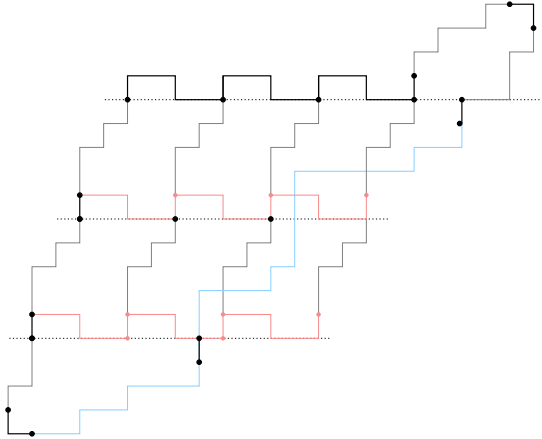}
	\end{center}
	\caption{The construction of 3-staircase polygons. Three vertical steps are
	selected as the heights of the indents. All indents are placed at the first
	choice.}
	\label{Fig:3-stair}
\end{figure}

\begin{prop}
	The generating function for $m$-staircase (resp. $m$-unimodal) polygons
	with left-side vertical indents at different heights is
	\begin{equation} \label{Eq:m-unimodal.same}
		\frac{y^{m+1}}{m!} \dd y \T^m \dd y \frac{\mathcal{Q}}{y},
	\end{equation}
	where $\T = (\SP/y)^2$ and $\mathcal{Q} = \SP/y$ (resp. $\UP/y$).
\end{prop}

\begin{proof}
	The form of 3-staircase polygons with only left-side
	vertical indents can be seen in Figure~\ref{Fig:3-stair}. This can be
	generalised to $m$-unimodal polygons. Staircase boxed regions are
	constructed as above, giving an $m-1\times m$ grid. By again factorising at
	the points of intersection, one can construct equivalence classes in terms
	of the boxes intersected by the bottom perimeter.
	
	For example, the generating function for all classes whose bottom perimeter
	intersects the \nth n box in the first row and does not intersect any box in
	the \nth {(n-1)} column, but rather follows the \nth n column down, out of the
	grid, can be written as
	\[
		\T^{m-1} \cdot \lr{\frac{\delta^{m-1} \T}{{\delta y}^{m-1}}}
				 \cdot \lr{\dd y \frac{\UP}{y}}.
	\]
	
	Hence, all classes of $m$-unimodal polygons that are defined by the
	intersections of boxes have a generating function that is determined by the
	relative position of the first box intersected in each column. For example,
	the path illustrated in Figure~\ref{Fig:3-stair} intersects the
	$3^{\mbox{\tiny rd}}$ box in the first row, and the \second box in the
	second row. This class could be uniquely described by the pair $(1,2)$ of
	the relative positions of these boxes. In general, one needs an
	$(m-1)$-tuple to classify classes as such. Noting that symmetric
	classifications also have the same generating function, one can show that
	the cardinalities of the sets of classes with the same generating function
	match that of the required set of classes of $m$-unimodal polygons.
\end{proof}

\begin{rem}
	The form of the generating function of $m$-unimodal polygons with only
	vertical indents on the left side is as follows:
	\[
		\frac{x(A+\sqrt{\Delta}B)}{y^{m+1}\Delta^{m+1/2}},
	\]
	where $A$ (resp. $B$) is a polynomial of order $4m-1$ (resp. $4m-2$) in both
	$x$ and $y$.
\end{rem}

\end{document}